\documentclass[10pt,reqno]{amsart}
\usepackage{amsmath,amstext,amsfonts,%
enumitem,%
extpfeil,graphicx,mathrsfs,%
mathtools,setspace,tikz,xparse}
\usepackage[alphabetic]{amsrefs}
\usepackage[unicode]{hyperref}
\usepackage[utf8]{inputenc}
\usepackage[T1]{fontenc}
\usepackage[normalem]{ulem}
\usepackage[cmtip,all]{xy}
\usepackage[light]{sourcesanspro}
\usepackage{lmodern}
\hypersetup{colorlinks=false,%
linkbordercolor={0.93 0.71 0.13},%
citebordercolor={1 0 .5},%
urlbordercolor={0.27 0.39 0.68}}%
\usetikzlibrary{calc,matrix,cd}

\usepackage{filecontents}

\DeclareMathOperator{\End}{End}
\DeclareMathOperator{\Hom}{Hom}

\DeclareMathOperator{\Ker}{ker}
\DeclareMathOperator{\Gal}{Gal}
\DeclareMathOperator{\Supp}{supp}
\DeclareMathOperator{\GL}{GL}
\DeclareMathOperator{\Sym}{Sym}
\DeclareMathOperator{\Img}{im}
\DeclareMathOperator{\Ind}{Ind}
\DeclareMathOperator{\ind}{ind}

\DeclareMathOperator{\Alg}{Alg}
\DeclareMathOperator{\Lisse}{Lisse}
\DeclareMathOperator{\Span}{span}

\newcommand{\m}[1]{\(\smash{#1}\)}
\newcommand{\algn}[1]{\begin{align*}\textstyle #1\end{align*}}
\newcommand{\Algn}[1]{\begin{align}\textstyle #1\end{align}}
\newcommand{\bb}[1]{\mathbb{#1}}
\newcommand{\QQ}{\bb{Q}}
\newcommand{\FF}{\bb{F}}
\newcommand{\ZZ}{\bb{Z}}
\newcommand{\CC}{\bb{C}}
\newcommand{\FFp}{\overline{\bb{F}}_p}
\newcommand{\QQp}{\smash{\overline{\bb{Q}}_p}}
\newcommand{\QQptimes}{\smash{\overline{\bb{Q}}_p^\times}}
\newcommand{\FFptimes}{\smash{\overline{\bb{F}}_p^\times}}
\newcommand{\ZZp}{\smash{\overline{\bb{Z}}_p}}

\newcommand{\defeq}{\vcentcolon=}
\newcommand{\SCR}{\mathscr}

\newcommand{\RHO}{\rho}
\newcommand{\HECKET}{\mathbf T}

\newcommand{\MU}{\mu}
\newcommand{\COVERA}[1]{\frac{C_{#1}}{a}}%
\newcommand{\COVERASQ}[1]{\frac{C_{#1}}{a^2}}%

\newcommand{\SIGMA}{\Sigma}
\newcommand{\ZETA}{\zeta}
\newcommand{\LAMBDA}{\lambda}

\newcommand{\LESS}[1]{%
{\uline{#1}}}
\newcommand{\MORE}[1]{%
{\overline{#1}}}
\newcommand{\BIGO}[1]{\smash{\mathsf{O}\!\left({#1}\right)}}
\newcommand{\NEBENPOWER}{\varrho}
\newcommand{\NEBENPOWEX}{\kappa}
\newcommand{\NEBENPOWEXX}{\kappa}
\newcommand{\NEBENPOWEY}{\kappa}

\newcommand{\MMM}[3]{M_{#1}({#3},{#2})}
\newcommand{\INZA}{I(n)Z}
\newcommand{\INZB}{n}
\newcommand{\SIZEEE}{4.25pt}
\newcommand{\VERIFYIF}[1]{[{#1}]}
\newcommand{\CONSTANT}[1]{C_{#1}}
\newcommand{\IRREDUCIBLE}[3]{\mathrm{Irr}_{#3}}

\newcommand{\BIRREDUCIBLE}[3]{\mathrm{Irr}_{\theta,#3}}
\newcommand{\ERRR}{\mathcal E}
\newcommand{\BRREDUCIBLE}[1]{\mathrm{Red}_{\theta,{#1}}}
\newcommand{\BPI}{\Pi_\theta}
\newcommand{\DIFF}{\hbar}
\newcommand{\TEMPORARYA}{\xi}
\newcommand{\TEMPORARYAA}{\vartheta}

\newcommand{\VECTORSP}{V_{\mathrm{coef}}}
\newcommand{\NEWPAGE}{}
\newcommand{\DOUBBR}[1]{[\![ {#1} ]\!]}
\newcommand{\aaa}{c}
\newcommand{\bbb}{d}
\newcommand{\ccc}{e}
\newcommand{\ddd}{f}
\newcommand{\GENN}{\overline{\GEN}}

\newcommand{\UNKNOWN}{\ast}
\newcommand{\TEMPORARYI}{\nu}

\newcommand{\absdet}{|\!\det\!|}
\newcommand{\IWASAWAVAR}{\QUARK}
\newcommand{\PADICVALUATION}{v_p}
\newcommand{\DELTAAA}{\delta}
\newcommand{\TADICVALUATION}{\smash{v_{\IWASAWAVAR}}}
\newcommand{\GEN}{\eta}

\newcommand{\AFIOFJ}{\\[-6pt]}
\newcommand{\AFIOFJFEWF}{\\[-13pt]}
\newcommand{\FUEWIFFEW}{\\[4pt]}
\newcommand{\LIFTA}[1]{{#1}}
\newcommand{\TEMPORARYAAAA}{\vartheta}
\newcommand{\LCLASS}{\mathscr{C}}
\newcommand{\LCLASSTEICH}{\mathscr{C}_{[\,]}}

\newcommand{\TILDE}{\widetilde}

\newcommand{\CIRCLES}{\SCR D}
\newcommand{\COMPAR}{\Psi}
\newcommand{\COMPARR}{\Sym(\Psi)}

\newcommand{\CONST}[2]{\gamma({#1},{#2})}
\newcommand{\GEQSLANT}{\geqslant}
\newcommand{\LEQSLANT}{\leqslant}
\newcommand{\QUARK}{\pi}

\newcommand{\QUOTIENT}{\mathfrak W}
\newcommand{\ZEROIDEAL}{\mathfrak I}
\newcommand{\EQUIVZEROIDEAL}{\equiv_{\mathfrak I}}
\newcommand{\EQUALZEROIDEAL}{=\hspace{0.51066em}}
\newcommand{\maxideal}{\mathfrak m}
\newcommand{\repnn}{\pi}
\newcommand{\NEWLINEA}{\\[2pt]}
\newcommand{\INAAA}{\SCR J}
\newcommand{\VAREPSILONp}{\varepsilon}

\newcommand{\SUBMODULE}{\mathrm{sub}}
\newcommand{\QUOTIENTI}{\mathrm{quot}}
\newcommand{\QUOTIENTMAP}{\Upsilon}

\newcommand{\MATRIX}{\mathrm{Van}}
\newcommand{\tempA}{u}
\newcommand{\tempB}{v}

\newcommand{\TEMPORARYMU}{\Phi}

\newcommand{\sqbrackets}[4]{{#1}\bullet_{#2,#3}{#4}}
\newcommand{\sqbracketss}[4]{{#1}\bullet_{#3}{#4}}
\newcommand{\Hypergeometric}[5]%
{\smash{\null_{#1}\mathrm{F}_{#2}\big(%
\smash{\genfrac{}{}{0pt}{2}{#3}{#4}};{#5}\big)}}

\newtheorem{theorem}{Theorem}
\newtheorem{lemma}[theorem]{Lemma}
\newtheorem{proposition}[theorem]{Proposition}
\newtheorem{corollary}[theorem]{Corollary}
\newtheorem{conjecture}{Conjecture}

\makeatletter
\newcommand{\dashover}[2][\mathop]{#1{\mathpalette\df@over{{\dashfill}{#2}}}}
\newcommand{\fillover}[2][\mathop]{#1{\mathpalette\df@over{{\solidfill}{#2}}}}
\newcommand{\df@over}[2]{\df@@over#1#2}
\newcommand\df@@over[3]{%
  \vbox{
    \offinterlineskip
    \ialign{##\cr
      #2{#1}\cr
      \noalign{\kern1pt}
      $\m@th#1#3$\cr
    }
  }%
}
\newcommand{\dashfill}[1]{%
  \kern-.5pt
  \xleaders\hbox{\kern.5pt\vrule height.4pt width \dash@width{#1}\kern.5pt}\hfill
  \kern-.5pt
}
\newcommand{\dash@width}[1]{%
  \ifx#1\displaystyle
    2pt
  \else
    \ifx#1\textstyle
      1.5pt
    \else
      \ifx#1\scriptstyle
        1.25pt
      \else
        \ifx#1\scriptscriptstyle
          1pt
        \fi
      \fi
    \fi
  \fi
}
\newcommand{\solidfill}[1]{\leaders\hrule\hfill}
\makeatother

\newenvironment{Smatrix}%
    {\tiny (\begin{smallmatrix}}%
    {\end{smallmatrix})}
\newenvironment{Lmatrix}%
    {\footnotesize (\begin{smallmatrix}}%
    {\end{smallmatrix})}
\newcommand{\BLACKSQUARE}{\hspace{\stretch{1}}\m{\mathbin{\vcenter{%
    \hbox{\rule{.9ex}{.9ex}}}}}}
\renewenvironment{proof}[1]%
    {\emph{{#1}}}%
    {\BLACKSQUARE}

\makeatletter
\newcommand{\xhookdoubleheadrightarrow}[2][]{%
  \lhook\joinrel
  \ext@arrow 0359\rightarrowfill@ {#1}{#2}%
  \mathrel{\mspace{-15mu}}\rightarrow}
\renewcommand{\xtwoheadrightarrow}[2][]{%
  \ext@arrow 0359\rightarrowfill@ {#1}{#2}%
  \mathrel{\mspace{-15mu}}\rightarrow}
\def\@settitle{\begin{center}%
  \baselineskip14\p@\relax\normalfont\LARGE
    \uppercasenonmath\@title
  \@title\end{center}%
}
\def\do#1{\@ifdefinable#1{\newdimen#1}}
\do\@TempWidthOne
\do\@TempWidthTwo
\newcommand*\WidthInMathUnitsOf[2]{%
    \settowidth\@TempWidthOne{$\m@th #1#2$}%
    \settowidth\@TempWidthTwo{$\m@th #1\mspace{1mu}$}%
    \strip@pt \dimexpr \@TempWidthOne*\p@/\@TempWidthTwo \relax
}
\newcommand*\SetToWidthInMathUnits[3]{%
    \settowidth\@TempWidthOne{$\m@th #2#3$}%
    \settowidth\@TempWidthTwo{$\m@th #2\mspace{1mu}$}%
    \edef #1{%
        \strip@pt \dimexpr \@TempWidthOne*\p@/\@TempWidthTwo \relax
    }%
}
\newcommand*\ExtractHexDigitOfFamilyOfSymbol[1]{%
    \expandafter\@ExtractHexDigitOfFamily \meaning #1%
}
\@ifdefinable\@ExtractHexDigitOfFamily{
    \def\@ExtractHexDigitOfFamily#1"#2#3#4#5{#3}
}
\makeatother

\newcommand*{\ReportValueOfOneEm}[1]{%
    &\text{in #1 size,}%
        &\quad\texttt{1em}%
            &=\texttt{\the \fontdimen 6 \csname #1font\endcsname 2}%
        &&\qquad\text{(from \texttt{\fontname \csname #1font\endcsname 2})}%
}
\newcommand*{\ReportSizeOfFontOfSymbol}[2]{%
    &\texttt{\pdffontsize\csname #2font\endcsname
                "\ExtractHexDigitOfFamilyOfSymbol{#1}}%
        &&\quad\text{in #2 size}%
        &\qquad(\texttt{\string\fontname}%
            &=\texttt{\fontname\csname #2font\endcsname
                "\ExtractHexDigitOfFamilyOfSymbol{#1}})%
}

\title%
[On the reductions of crystabelline representations]%
{On the reductions of certain  two-dimensional crystabelline representations}
\author{Bodan Arsovski}

\begin{document}
\maketitle

\begin{abstract}
Crystabelline representations are  representations of the absolute Galois group \m{G_{\QQ_p}} over \m{\QQp} that become crystalline on \m{G_{F}} for some abelian extension \m{F/\QQ_p}. Their relation to modular forms is that the %
 representation associated with a finite slope %
 newform of level divisible by \m{p^2} %
  is
 crystabelline. %
   In this article we study %
  the connection between the slopes of %
   two-dimensional crystabelline representations %
   and the reducibility of their modulo~\m{p} reductions. %
   This question is inspired by a theorem by Buzzard and Kilford which implies that the %
    slopes on the boundary of the \m{2}-adic eigencurve of tame level \m{1} are %
    integers (and in %
    arithmetic progression); an analogous theorem by Roe which says that the same is true for the \m{3}-adic eigencurve;  Coleman's halo conjecture and the ghost conjecture which  give predictions about the %
     slopes on the   \m{p}-adic eigencurve of general tame level; %
     and %
      Hodge theoretic 
       conjectures by Breuil, Buzzard, Emerton, and Gee which indicate that there is a connection between all of these  %
      and %
        the slopes of locally reducible two-dimensional crystabelline representations.
       We prove that the reductions %
       of certain two-dimensional crystabelline representations with slopes in  \m{(0,\frac{p-1}{2})\backslash \ZZ} are usually irreducible, with the exception of a small region where the slopes are half-integers and reducible representations do  %
       occur.
\end{abstract}

\section{Introduction}

Throughout this article we fix the odd prime number \m{p}, the  integer \m{n\in\ZZ_{\GEQSLANT2}}, the   primitive \m{p^n}th root of unity \m{\ZETA_F\in\QQp}, and the extension \m{F=\QQ_p(\ZETA_F)} of \m{\QQ_p}. The main objects we study are two-dimensional crystabelline representations of the absolute Galois group \m{G_{\QQ_p}} over \m{\QQp}. These  are defined by Berger and Breuil  in~\cite{BB10} as the potentially crystalline representations of \m{G_{\QQ_p}} that become crystalline on \m{G_{K}} for some abelian extension \m{K/\QQ_p}. %
 We restrict our discussion to the  extension \m{F/\QQ_p}. The crystabelline representations we consider are  parametrized by an integer  \m{k \in\ZZ_{\GEQSLANT 2}}, an element \m{a \in \ZZp} such that \m{\PADICVALUATION(a) \in (0,k-1)} and \m{a^2\ne p^{k-1}}, and %
  a character \m{\varepsilon:\ZZ_p^\times \to \QQptimes} which has conductor \m{p^n}. 
   For the remainder of this article we fix parameters \m{(k,a,\varepsilon)} with these properties, and we 
 denote the corresponding  crystabelline representation  by  \m{V_{k,a,\varepsilon}}. %
 We give a %
 construction of \m{V_{k,a,\varepsilon}} and the definition of its modulo \m{p} reduction \m{\overline{V}_{k,a,\varepsilon}} in section~\ref{lSecCrys}. 
 For \m{s\in\{1,\ldots,n-1\}}, we fix the primitive \m{p^{s}}th root of unity
     \m{\textstyle{}\ZETA_s = \ZETA_s (\varepsilon) =\varepsilon^{-1}(1-p^{n-s})}, and we write  \m{\ZETA=\ZETA_{n-1}}. We also fix a root %
    \m{\QUARK= p^{1/(p-1)p^{n-2}}\in\QQp} and define \m{\textstyle \TADICVALUATION =(p-1) p^{n-2}\cdot \PADICVALUATION}, so that \m{\TADICVALUATION} is the re-normalization of \m{\PADICVALUATION} such that \m{\TADICVALUATION(\QUARK)=\TADICVALUATION(1-\ZETA)=1}. %
    In particular, the \m{\TADICVALUATION}-valuation of any element of \m{\ZZ_p[\ZETA]} is an integer. When we refer to the ``slope'' of \m{V_{k,a,\varepsilon}} we mean the re-normalized valuation \m{\TADICVALUATION(a)}, and we define \m{\TEMPORARYI = \lceil \TADICVALUATION(a)\rceil \in \ZZ_{\GEQSLANT1}} as its ceiling. %
    The main question we consider in this article is whether reducibility of the modulo \m{p} representation \m{\overline{V}_{k,a,\varepsilon}} implies integrality of the slope \m{\TADICVALUATION(a)}. Let us phrase this question in the following form (see conjecture~4.2.1 in~\cite{BG16}).
\begin{conjecture}\label{MAINCONJ}
If \m{\TADICVALUATION(a) \not\in \ZZ} then \m{\overline{V}_{k,a,\varepsilon}} is irreducible.
\end{conjecture}

\subsection{The main theorem in this article} Let us fix an integer \m{\NEBENPOWER}  such that the reduction modulo \m{p} of \m{\VAREPSILONp} is the map \m{\UNKNOWN \mapsto {\UNKNOWN^{\NEBENPOWER}}}. In other words, \m{\VAREPSILONp} uniquely determines the congruence class \m{\NEBENPOWER + (p-1)\ZZ \in \ZZ/(p-1)\ZZ}, and we fix a  representative \m{\NEBENPOWER\in\ZZ} of this  class. Let also \m{\NEBENPOWEXX=k+\NEBENPOWER-2}, so that if \m{x^{k-2}} denotes the character \m{\ZZ_p^\times \to \QQptimes} which maps \m{z} to \m{z^{k-2}} then  the reduction modulo \m{p} of \m{x^{k-2}\VAREPSILONp} is the map \m{\UNKNOWN \mapsto {\UNKNOWN^{\NEBENPOWEXX}}}.
 For \m{\alpha\in \ZZ_{\GEQSLANT 1}}, %
 let us define \m{\CIRCLES_{\alpha}=\CIRCLES_{\alpha}(\VAREPSILONp)} as
    \algn{\CIRCLES_{\alpha} = \left\{x \in \QQp \,\big|\,\TADICVALUATION\left(x^2-\frac{1}{\alpha!(\alpha-1)!}(\VAREPSILONp(1+p)-1)^{2\alpha-1}\right)\GEQSLANT 2\alpha-\tfrac12\right\}.}
So \m{\CIRCLES_{\alpha}} is the union of the two closed disks that are centered at the square roots of \m{\frac{1}{\alpha!(\alpha-1)!}(\VAREPSILONp(1+p)-1)^{2\alpha-1}} and have radii \m{\QUARK^{-\alpha+1/4}}. Both centers of the disks have \m{\TADICVALUATION}-valuation \m{\alpha-\tfrac12}, and therefore all points of \m{\CIRCLES_{\alpha}} have \m{\TADICVALUATION}-valuation \m{\alpha-\tfrac12}. For \m{h\in\ZZ}, let  \m{\MORE{h}} be the number in \m{\{1,\ldots,p-1\}} that is congruent to \m{h} modulo \m{p-1}.  For \m{l \in \ZZ}, let us define \m{\IRREDUCIBLE{k}{\NEBENPOWER}{l}=\IRREDUCIBLE{k}{\NEBENPOWER}{l}(\NEBENPOWEXX)} as
    \algn{\textstyle \IRREDUCIBLE{k}{\NEBENPOWER}{l} = \ind(\omega_2^{\MORE{\kappa-2l+1}})\otimes \omega^{l}. }%
We prove the following theorem.
\begin{theorem}\label{MAINTHM1}
If \m{\TADICVALUATION(a) \in (0,\frac{p-1}{2})\backslash \ZZ}   then
    \algn{\overline V_{k,a,\varepsilon}\in \left(\{\Pi\} \cup\left\{\IRREDUCIBLE{k}{\NEBENPOWER}{\NEBENPOWEXX-j} %
    \right\}_{ 1\LEQSLANT j\LEQSLANT\TEMPORARYI-1}\right) \big\backslash \left\{\IRREDUCIBLE{k}{\NEBENPOWER}{j} 
    \right\}_{ 1\LEQSLANT j\LEQSLANT\TEMPORARYI-1},}
where
    \algn{\Pi \cong \left\{\begin{array}{rl}\mu_\lambda \omega^{2\TEMPORARYI-1} \oplus \mu_{\lambda^{-1}}\omega^{2\TEMPORARYI-1}%
    & \text{if } \NEBENPOWEXX\equiv 2\TEMPORARYI-1 \bmod {p-1} \text{ \& } a\in\CIRCLES_{\TEMPORARYI} \text{,} \\[5pt] \IRREDUCIBLE{k}{\NEBENPOWER}{\NEBENPOWEXX} & \text{otherwise,} \end{array}\right.}
and \m{\lambda + \lambda^{-1}} is the reduction modulo \m{\maxideal} of 
\algn{\textstyle  ((\TEMPORARYI-1)!)^{-1} (1-\VAREPSILONp(1+p))^{\TEMPORARYI-1} a^{-3} \left( \frac{1}{\TEMPORARYI!(\TEMPORARYI-1)!} (\VAREPSILONp(1+p)-1)^{2\TEMPORARYI-1}-a^2\right) .}
\end{theorem}
The following corollary is the particular case of theorem~\ref{MAINTHM1} when \m{\TADICVALUATION(a) \in (0,1)}. %

\begin{corollary}\label{MAINCOROLL1}
If \m{\TADICVALUATION(a) \in (0,1) } then
    \algn{\overline V_{k,a,\varepsilon} \cong \left\{\begin{array}{rl}\mu_\lambda \omega \oplus \mu_{\lambda^{-1}}\omega %
    & \text{if } \NEBENPOWEXX\equiv 1 \bmod {p-1} \text{ \& } \TADICVALUATION\left(a^2+1-\VAREPSILONp(1+p)\right)\GEQSLANT \tfrac32 \text{,} \\[5pt] \IRREDUCIBLE{k}{\NEBENPOWER}{\NEBENPOWEXX} & \text{otherwise,} \end{array}\right.}
 where \m{\lambda + \lambda^{-1}} is the reduction of \m{a^{-3}(\VAREPSILONp(1+p)-1-a^2)} modulo \m{\maxideal}.
\end{corollary}
So if \m{\TADICVALUATION(a) \in (0,\frac{p-1}{2})\backslash \ZZ}  then \m{\overline V_{k,a,\varepsilon}} is \emph{usually} irreducible,  all exceptions   satisfy \m{\TADICVALUATION(a) \in \frac12\ZZ}, and exceptions do exist.
Theorem~\ref{MAINTHM1} gives a list of at most \m{ \TEMPORARYI=\lceil\TADICVALUATION(a)\rceil} possibilities for %
\m{\overline V_{k,a,\varepsilon}}. %
 In~\cites{Ars1,Ars3} we prove that, if \m{\PADICVALUATION(a) \in (0,\frac{p-1}{2})\backslash \ZZ}  and some additional conditions are satisfied, there are exactly \m{\lceil \PADICVALUATION(a)\rceil} possibilities %
 for a certain %
 representation \m{\overline{V}_{k,a}} and all of them do  occur.
In addition, there are clear similarities between %
  theorem~\ref{MAINTHM1} and the  theorems about \m{\overline{V}_{k,a}}  proven in~\cites{BG09, BG13, BG15, GG19a}. %
 Based on this we expect that all of the possibilities for \m{\overline{V}_{k,a,\VAREPSILONp}} listed in theorem~\ref{MAINTHM1} %
  occur and that the equations determining exactly \emph{when} they occur are very complicated.

 \subsection{Relation to modular forms} The Hodge theoretic construction \m{V_{k,a,\varepsilon}} relates to the theory of modular forms because the representation associated with a finite slope %
 newform
  \algn{\textstyle f = \sum_{n=1}^\infty a_n q^n}
 ---which is normalized so that \m{a_1=1},  has
 weight \m{k}, level \m{\Gamma_1(Np^n)} with \m{(N,p)=1},  ramified character \m{\psi} which has conductor \m{p^n} and is such that \m{\psi(p)=1} (so that we can think of \m{\psi} as a character of \m{\ZZ_p^\times}), and is such that \m{a_p^2\ne p^{k-1}}---is crystabelline and equal to \m{V_{k,a_p,\psi}}. 
 
 Conjecture~\ref{MAINCONJ} %
  is motivated by %
   several known results %
 about the  slopes of modular forms. The main theorem in~\cite{BK05} implies that, for \m{m\in\ZZ_{\GEQSLANT2}}, the  \m{2}-adic slopes of \m{2}-adic classical modular forms of %
  level \m{2^m} and character which has conductor \m{2^m} are integer multiples of  \m{2^{3-m}} and in arithmetic progression. %
 The main theorems in~\cites{Roe14,Kil08,KM12} imply the \m{3}-adic analogue of this, %
 a similar result for \m{5}-adic modular forms of level \m{25}, and a similar result for \m{7}-adic modular forms of level \m{49}, respectively.  In fact, all of these apply more generally to overconvergent modular forms,  and are best phrased in terms of the slopes on the boundary of the %
  eigencurve (see~\cites{CM98,Buz07}). 
 The common theme is that the re-normalizations of the slopes appearing in these theorems are always integers. There are two general conjectures in this vein. 
 Coleman's  halo conjecture predicts, among other things, that the suitably normalized slopes on the boundary of the \m{l}-adic eigencurve are in a disjoint union of arithmetic progressions of rational numbers (see~\cite{LWX17}). %
 Bergdall and Pollack's ghost conjecture (see~\cites{BP1,BP2}) unifies, for \m{\Gamma_0(N)}-regular primes \m{l} with  \m{(N,l)=1},  Buzzard's conjecture about the slopes of modular forms of level \m{\Gamma_0(Nl)} (see Buzzard's algorithm in~\cite{Buz05}), and  predictions about the slopes of modular forms of level \m{\Gamma_1(Nl^m)} (\m{m\in\ZZ_{\GEQSLANT2}}) and character which has conductor \m{l^m}. Notably, all of the slopes predicted by Buzzard are integers, and \m{\Gamma_0(N)}-regularity relates to Galois representations because %
  an odd prime \m{l} is \m{\Gamma_0(N)}-regular if and only if, for all weights \m{t\GEQSLANT 2} and all \m{f \in S_t(\Gamma_0(N))}, the modulo \m{l} representation associated with \m{f} is reducible. Therefore these predictions naturally lead to a conjecture by Breuil, Buzzard, and Emerton about locally reducible crystalline representations (see conjecture~4.1.1 in~\cite{BG16}), %
  and to the analogous conjecture~\ref{MAINCONJ}.

   \begin{corollary}\label{MAINCOR3prime}
    If \m{f} is a %
    newform which has weight \m{k}, level \m{\Gamma_1(Np^n)} with \m{(N,p)=1}, and %
 \m{U_p}-eigenvalue \m{a_p} such that %
    \m{\TADICVALUATION(a_p)\in (0,p-1)\backslash\frac12\ZZ},  %
    then the residual representation associated with \m{f} is irreducible.
\end{corollary}
  
  \begin{proof}{Proof. }
  Note that
 \m{a_p^2\ne p^{k-1}}. We may without loss of generality assume that the character \m{\psi} of \m{f} satisfies \m{\psi(p)=1} (so that we can think of it as a character of \m{\ZZ_p^\times}). %
  Since \m{a_p\ne 0} it follows from part~(2b) of theorem~9.1.10 in~\cite{RS} that \m{\psi} has conductor \m{p^n}. It is shown in section~3.2 of~\cite{GM09} that the Galois representation associated with \m{f} is  \m{\textstyle V_{k,a_p,\psi}}, and theorem~\ref{MAINTHM1} implies that \m{\textstyle \overline{V}_{k,a_p,\psi}} is irreducible. %
\end{proof}

 For \m{N\in\ZZ_{\GEQSLANT1}} such that \m{(N,p)=1}, let \m{\SCR C_N} denote the %
 eigencurve of tame level \m{N}, and let %
 the ``boundary of \m{\SCR C_N}'' be the subset consisting of the points whose characters satisfy %
 \m{|\chi(1+p)-1| \GEQSLANT p^{-1/(p-1)}}   (see~\cite{LWX17} for more details). We can associate  with each point of \m{\SCR C_N} a residual representation, and by the ``locally reducible part of \m{\SCR C_N}'' we mean the part of \m{\SCR C_N} consisting of the points whose associated residual representations are reducible. The intersection of the boundary of \m{\SCR C_N} and the classical modular locus of \m{\SCR C_N} consists of points corresponding to finite slope newforms of level \m{\Gamma_1(Np^m)} (\m{m\in\ZZ_{\GEQSLANT2}}).  Therefore corollary~\ref{MAINCOR3prime}  implies the following.  %

  \begin{corollary}\label{MAINCOR3}
    The re-normalized slopes of classical points on the locally reducible part of the boundary of \m{\SCR C_N} %
    that are less than \m{\frac{p-1}{2}} are half-integers.
\end{corollary}

\subsection*{\texorpdfstring{Acknowledgements}{Acknowledgements}}\label{ackref}
I would like to thank professor Kevin Buzzard for his helpful remarks, suggestions, and support. This work was supported by Imperial College London and its President's PhD Scholarship.

\section{Crystabelline representations}
\label{lSecCrys}

Recall that we fix the integers \m{n,k\in\ZZ_{\GEQSLANT2}}, the   primitive \m{p^n}th root of unity \m{\ZETA_F\in\QQp},  the extension \m{F/\QQ_p=\QQ_p(\ZETA_F)/\QQ_p}, the element \m{a \in \ZZp} such that \m{\PADICVALUATION(a) \in (0,k-1)} and \m{a^2\ne p^{k-1}}, and %
  the character \m{\varepsilon:\ZZ_p^\times \to \QQptimes} with conductor \m{p^n}.  %
 Note that the maximal unramified subfield of \m{F} is \m{\QQ_p}. Let \m{\maxideal} be the maximal ideal of \m{\ZZp}. Let \m{D_{k,a,\varepsilon} %
 =  \overline{\QQ}_p e_1 \oplus \overline{\QQ}_p e_2} be the weakly admissible filtered \m{\varphi}-module with Hodge--Tate weights \m{(0,k-1)}, a filtration 
    \algn{\mathrm{Fil}^{k-1} (D_F) = \left\{\begin{array}{rl} F\otimes_{\QQ_p} D_{k,a,\VAREPSILONp} & \text{if } i \LEQSLANT 0\text{,} \\[2pt]  F\otimes_{\QQ_p} \QQp (e_1 + e_2) & \text{if } 0 < i < k\text{,}\\[2pt]
0 & \text{if } i \GEQSLANT k\text{,}  \end{array}\right.} %
a Frobenius map \m{\varphi} such that
\algn{\textstyle \left\{\begin{array}{ll} \varphi (e_1) = a^{-1}p^{k-1} e_1, &  \\[2pt] \varphi (e_2) = a e_2, &  \end{array}\right.}
 and an action of the Galois group \m{\Gal(F/\QQ_p)} such that
 \algn{\textstyle \left\{\begin{array}{ll} g(e_1) = \VAREPSILONp(g)^{-1}e_1, &  \\[2pt] g(e_2)=e_2. &  \end{array}\right.}
 Here we can view  \m{\VAREPSILONp} as a map \m{\Gal(F/\QQ_p) \to \QQptimes} since %
 \m{\Gal(F/\QQ_p)\cong(\ZZ/p^n\ZZ)^\times} (see subsection~1.2 of~\cite{BB10}). 
 By proposition~2.4.5 in~\cite{BB10}, there is a unique two-dimensional representation  \m{V_{k,a,\varepsilon}} that %
  becomes crystalline on \m{F}, has Hodge--Tate weights \m{(0,k-1)}, and is such that \m{\textstyle D_{cris}(V_{k,a,\varepsilon}^\ast|_{G_{F}}) = D_{k,a,\varepsilon}}, and this representation is absolutely irreducible.  Let \m{\overline{V}_{k,a,\varepsilon}} be the semi-simplification of the reduction modulo  \m{\maxideal} of a Galois stable %
  lattice in \m{{V_{k,a,\varepsilon}}}. This construction is independent of the choice of  lattice due to a theorem by Brauer and
Nesbitt (see~\cite{Ber10}). Berger and Breuil study  \m{{V_{k,a,\varepsilon}}} in detail in~\cite{BB10}, and the articles~\cites{Bre03a,Bre03b,Bre04} contain further background.

\section{The local Langlands correspondence}
\label{SECll}

In this section we %
 rephrase theorem~\ref{MAINTHM1} in terms of certain \m{\GL_2(\QQ_p)}-representations. %

 \subsection{\protect\boldmath The principal series representation associated with \texorpdfstring{\m{{V_{k,a,\varepsilon}}}}{V}}\label{subsecWe} For \m{x\in R^\times} with \m{R\in\{\QQp,\FFp\}}, let \m{\mu_x} be the unramified character of the Weil group that sends the geometric Frobenius to \m{x}. Let \m{|\ | : \QQp \to \overline\QQ} be the \m{p}-adic norm, and %
 let us fix embeddings \m{\overline{\QQ}\hookrightarrow \QQp} and \m{\overline{\QQ}\hookrightarrow \CC}. The characters in the  Weil--Deligne representation corresponding to \m{V_{k,a,\varepsilon}} are \m{\VAREPSILONp{\mu_{a p^{1-k}}}} and \m{\mu_{a^{-1}}}, %
 so the corresponding admissible representation of \m{G = \GL_2(\QQ_p)} via the  local Langlands correspondence is the principal series representation
    \algn{
    & \textstyle\mathrm{PS} (\VAREPSILONp{\mu_{a p^{1-k}}},{\mu_{a^{-1}}})
    \\ &\null\qquad\textstyle \defeq \{\text{smooth maps } G \to \CC \,|\, f (\begin{Lmatrix}
        x & y \\ 0 & z
        \end{Lmatrix}g) = \VAREPSILONp(x){\mu_{a  p^{1-k}}(x)}\mu_{a^{-1}}(z) |x/z|^{1/2} f(g)\}.}
Let \m{B} be the Borel subgroup of  \m{G=\GL_2(\QQ_p)}  consisting of those elements that are upper triangular, 
 let \m{K=\GL_2(\ZZ_p) \subset G}, let \m{Z} be the center of \m{G}, and, for \m{m\in\ZZ_{\GEQSLANT1}}, let   \m{I(m)} be the subgroup of \m{K} consisting of those elements that are upper triangular modulo \m{p^m}.   Let us define the character %
    \algn{  & \textstyle \RHO:\begin{Lmatrix} x&y\\0&z\end{Lmatrix} \in B\xmapsto{ \ \ } \VAREPSILONp(x)\mu_{a  p^{1-k}}(x)\mu_{a^{-1}}(z)|x/z|^{1/2} \in \QQptimes,} so that \m{\textstyle\Ind_B^G\RHO} %
    is the \m{\QQp}-representation that corresponds to \m{\mathrm{PS} (\VAREPSILONp{\mu_{a p^{1-k}}},{\mu_{a^{-1}}})}. Due to a theorem by  Bernstein and Zelevinsky, \m{\Ind_B^G\RHO} must be irreducible since the ratio of the characters \m{\VAREPSILONp\mu_{ap^{1-k}}} and \m{\mu_{a^{-1}}} is neither equal to the norm nor to its inverse.

 \subsection{Compact induction of locally algebraic representations}\label{loccomp}
  Let \m{W} be a finite-dimensional representation of a closed subgroup \m{H} of \m{G} over \m{\FF \in \{\QQp,\FFp\}}. By a locally algebraic   map \m{H\to W} we mean a map which on an open subgroup of \m{H} is the restriction of a rational map on (the algebraic group) \m{H}, and we say that \m{W} is locally algebraic if the map \m{h \in H \mapsto h w \in W} is locally algebraic for all \m{w \in W}.
 For a closed subgroup \m{H} of \m{G} and a locally algebraic finite-dimensional representation \m{W} of \m{H} we define the compact induction of \m{W} by
    \algn{\textstyle \ind_H^G W \defeq %
     \textstyle \big\{\text{locally algebraic } G \to W \,\,\big|\,\, & f (h g) = h f(g) \text{ for all } h \in H \\ &\textstyle{} \text{\& }\Supp f\text{ is compact in }H\backslash G \}.}
 Suppose that \m{H} is moreover open. %
  For \m{g\in G} and \m{w \in W}, let \m{\sqbrackets{g}{G,H}{\FF}{w}} be the unique element of \m{\ind_H^G W} that is supported on \m{Hg^{-1}} and maps \m{g^{-1}} to \m{w}. Since \m{H\backslash G} is discrete, %
  \m{\ind_H^G W} is the linear span of \m{\sqbrackets{g}{G,H}{\FF}{w}} for all \m{g\in G} and \m{w\in W}. %
    For open subgroups \m{H_1\subseteq H_2} of \m{G} such that \m{[H_2:H_1]} is finite, let us define
    \algn{\ind_{H_1}^{H_2} W \defeq \{\text{maps } H_2 \to W \,|\, f (h_1h_2) = h_1f(h_2)\}}
and let \m{\sqbrackets{g}{H_2,H_1}{\FF}{w}} be the unique element of \m{\ind_{H_1}^{H_2} W} that is supported on \m{H_1g^{-1}} and maps \m{g^{-1}} to \m{w}. %
 It is easy to check that \m{\sqbrackets{g}{G,H_2}{\FF}{\sqbrackets{h_2}{H_2,H_1}{\FF}{w}}} corresponds to \m{\sqbrackets{gh_2}{G,H_1}{\FF}{w}} under the isomorphism \m{\ind_{H_2}^{G}\ind_{H_1}^{H_2}W \cong \ind_{H_1}^GW}, and that \algn{\textstyle g_1(\sqbrackets{g_2}{G,H}{\FF}{(h w)}) = \sqbrackets{(g_1g_2h)}{G,H}{\FF}{w}} for all \m{g_1,g_2\in G}, \m{h \in H}, and \m{w\in W}. For \m{\xi \in \FF}, let \m{\textstyle{}\sqbrackets{\xi}{G,H}{\FF}{w} = \xi(\sqbrackets{\begin{Smatrix}
1 & 0 \\ 0 & 1
\end{Smatrix}}{G,H}{\FF}{w})}.
 For \m{g\in G}, \m{w\in W}, and \m{m\in \ZZ_{\GEQSLANT1}}, let us use the shorthands
 \algn{\sqbrackets{g}{G,I(m)Z}{\FF}{w} & \textstyle= \sqbrackets{g}{m}{\FF}{w}, \\ \sqbrackets{g}{G,KZ}{\FF}{w} & \textstyle{} = \sqbracketss{g}{KZ}{\FF}{w}.}
 Let us also define the character \m{\tau : I(n)Z\to\QQptimes} by %
    \algn{ & \textstyle \tau \left(\begin{Lmatrix}  \aaa&\bbb\\ \ccc p^n &\ddd \end{Lmatrix}\begin{Lmatrix} p^w & 0 \\ 0 & p^w \end{Lmatrix}\right) =\VAREPSILONp(\aaa)  \text{ for } \begin{Lmatrix}  \aaa&\bbb\\ \ccc p^n &\ddd \end{Lmatrix} \in I(n) \text{ and } w \in \ZZ.}
For \m{h \in \ZZ_{\GEQSLANT0}}, let %
\m{{\Sym}^{h}(\overline{\QQ}_p^2)} be the \m{G}-module of homogeneous polynomials in \m{x} and \m{y} of total degree~\m{h} with coefficients in \m{\QQp}, with \m{G} acting by \algn{\textstyle \begin{Lmatrix}g_1 & g_2 \\ g_3 & g_4
\end{Lmatrix} \cdot v(x,y) = %
v(g_1x+g_3y,g_2x+g_4y)} for \m{\begin{Smatrix}
g_1 & g_2 \\ g_3 & g_4
\end{Smatrix} \in G} and \m{v \in {{\Sym}^{h}(\overline{\QQ}_p^2)}}.
     In particular,  \m{\begin{Smatrix}
    p & 0 \\ 0 & p
    \end{Smatrix}v = p^{h}v} for all \m{v \in {\Sym}^{h}(\overline{\QQ}_p^2)}. %
     For \m{R \in \{\QQp,\ZZp,\FFp\}}, let 
     \m{\textstyle{}\uline{\Sym}^{h}(R^2)} be  the \m{KZ}-module %
     of homogeneous polynomials in \m{x} and \m{y} of total degree~\m{h} with coefficients in \m{R}, with \m{\begin{Smatrix}
    p & 0 \\ 0 & p
    \end{Smatrix}} acting trivially and \m{K} acting by \algn{\textstyle \begin{Lmatrix}g_1 & g_2 \\ g_3 & g_4
\end{Lmatrix} \cdot v(x,y) = v(g_1x+g_3y,g_2x+g_4y)} for \m{\begin{Smatrix}
g_1 & g_2 \\ g_3 & g_4
\end{Smatrix} \in K} and \m{v \in {\uline{\Sym}^{h}}({R^2)}}. Then \m{{\uline{\Sym}^{h}}({R^2)}} is also an \m{I(m)Z}-module for all \m{m\in\ZZ_{\GEQSLANT1}}.
 Let us define  the \m{I(n)Z}-modules   %
    \algn{\textstyle \TILDE\SIGMA_{k-2} & \textstyle = {\uline{\Sym}^{k-2}}({\overline{\QQ}_p^2)}\otimes \tau, \\
   \SIGMA_{k-2} & \textstyle = \text{the reduction of } \uline{\Sym}^{k-2}(\overline{\ZZ}_p^2) \otimes \tau \text{ modulo } \maxideal ,} and, for \m{h\in\ZZ_{\GEQSLANT0}}, the \m{KZ}-module  \algn{\sigma_h = \uline{\Sym}^h(\FFp^2).} For \m{t\in\{0,\ldots,p-1\}}, \m{\lambda \in\FFp}, and %
 \m{\psi:\QQ_p^\times \to \overline\FF_p^\times}, let %
    \algn{\textstyle \repnn(t,\lambda,\psi) = ({\ind_{KZ}^{G} \sigma_t}/{(T_{\sigma}-\lambda)})\otimes\psi,} where \m{T_\sigma\in \End_G(\ind_{KZ}^{G} \sigma_t)} is the Hecke operator corresponding to the double coset of \m{\begin{Smatrix} p&0 \\ 0 &1\end{Smatrix}}. Theorem~2.7.1 in~\cite{Bre03a}  classifies the representations \m{\repnn(t,\lambda,\psi)}. Let \m{\omega} be the modulo~\m{p} reduction of the cyclotomic character. Let \m{\ind(\omega_2^{t+1})} be the unique irreducible representation whose determinant is \m{{\omega^{t+1}}} and that is equal to \m{{\omega_2^{t+1}\oplus\omega_2^{p(t+1)}}} on inertia.  For \m{h\in\ZZ}, let \m{\MORE{h}\in\{1,\ldots,p-1\}}  and \m{\LESS{h}\in\{0,\ldots,p-2\}} be the numbers in the corresponding sets that are congruent to \m{h} modulo \m{p-1}.
    
   \subsection{\protect\boldmath A map  
 from %
 compact to parabolic induction} \label{subcompar}
  Let \algn{L : \absdet^{\frac{k-2}{2}}\otimes\tau \xrightarrow{\ \ } \absdet^{-\frac{1}{2}}\otimes\RHO} be a  linear isomorphism of one-dimensional vector spaces.
Let \algn{\textstyle\COMPAR\in\Hom_G(\absdet^{\frac{k-2}{2}}\otimes \ind_{\INZA}^{G}\tau, \absdet^{-\frac{1}{2}}\otimes \Ind_B^G \RHO)} be the map defined by linearly extending
    \algn{\textstyle \COMPAR\left( {\gamma}\sqbrackets{}{\INZB}{\QQp}{w}\right) = \left(g \mapsto \delta(g\gamma)L(w)\right) \text{ for } \gamma\in G \text{ and } w\in\absdet^{\frac{k-2}{2}}\otimes\tau ,}
 where, for \m{g\in G}, 
    \Algn{\label{eq3goigj45fjiugifweko3ig}\delta(g) = \left\{\begin{array}{rl}
    |\!\det(b)|^{-\frac12}\rho(b)\tau(u) & \text{if }g = bu \text{ with } b\in B \text{ 
     and } u \in I(n)\text{,}\\ 0 & \text{if } g \not \in BI(n)\text{.}
    \end{array}\right.}
Since \algn{\textstyle\delta(bg)L(w)=b\delta(g)L(w) &\vphantom{w\in \absdet^{\frac{k-2}{2}}\otimes\tau} \textstyle\text{ for } {b \in B}, {g \in G},  w\in\absdet^{\frac{k-2}{2}}\otimes\tau,\\
\textstyle \delta(u)=1 &\vphantom{w\in \absdet^{\frac{k-2}{2}}\otimes\tau} \textstyle\text{ for  } u\in 
\Gamma_1(p^n),\\
\textstyle \gamma_2(g\mapsto\delta(g\gamma_1)L(w)) = (g\mapsto\delta(g\gamma_2\gamma_1)L(w)) &\vphantom{w\in \absdet^{\frac{k-2}{2}}\otimes\tau} \textstyle \text{ for } {\gamma_1,\gamma_2\in G}, w\in \absdet^{\frac{k-2}{2}}\otimes\tau,\\
\textstyle \COMPAR \left({1}\sqbrackets{}{\INZB}{\QQp}{w}\right)\left(\begin{Lmatrix}
1 & 0 \\ 0 & 1
\end{Lmatrix}\right)=L(w) &\vphantom{w\in \absdet^{\frac{k-2}{2}}\otimes\tau}\textstyle \text{ for }w\in \absdet^{\frac{k-2}{2}}\otimes\tau,}
 \m{\COMPAR} is  a well-defined non-zero 
 homomorphism,    %
  which then must be surjective since \m{\Ind_B^G \RHO} is irreducible. Let  \m{\COMPARR} denote the induced surjective homomorphism   \Algn{\label{MAINEQN1}%
    \ind_{\INZA}^{G}\TILDE\SIGMA_{k-2} \xtwoheadrightarrow{\ \COMPARR \  %
    } {{\Sym}^{k-2}}({\overline{\QQ}_p^2)} \otimes \absdet^{-\frac{1}2} \otimes \Ind_B^G \RHO}
obtained from \m{\COMPAR} by tensoring with \m{{{\Sym}^{k-2}}({\overline{\QQ}_p^2)}}. To be more specific,  we use the general fact that if \m{H} is a closed subgroup of \m{G} and \m{U} is a \m{G}-representation and \m{W} is an \m{H}-representation, then there is an isomorphism \algn{\textstyle U \otimes \ind_{H}^G W \xrightarrow{ \ \cong  \ } \ind_{H}^G(%
U|_H\otimes W)} given by \m{u\otimes f \mapsto (g\mapsto gu\otimes f(g))}.
Let \m{U = {{\Sym}^{k-2}}({\overline{\QQ}_p^2)}}, %
 let us identify the %
 vector spaces underlying \m{U} and \m{\TILDE\SIGMA_{k-2} = (U\otimes \absdet^{\frac{k-2}2})|_{I(n)Z}\otimes\tau}, and write \m{(g,v)\mapsto g\cdot v} for the action of \m{G} on \m{U} and \m{(h,v)\mapsto hv} for the action of \m{I(n)Z} on \m{\TILDE \SIGMA_{k-2}}. %
  Then
\m{\COMPARR} is the map that fits in the commutative diagram\vspace{-4pt}
\begin{center}
    \begin{tikzpicture}[commutative diagrams/every diagram]
\node (P0) at (165:3cm) {$U%
\otimes \absdet^{\frac{k-2}{2}}\otimes \ind_{\INZA}^{G}\tau$};
\node (P1) at (15:3cm) {$\ind_{\INZA}^{G}\TILDE\SIGMA_{k-2}$} ;
\node (P2) at (270:0.8cm) {$U%
\otimes \absdet^{-\frac{1}2} \otimes \Ind_B^G \RHO $};
\path[commutative diagrams/.cd, every arrow, every label]
(P0) edge node {\m{\cong}} (P1)
(P0) edge node[left] {\m{U%
 \hspace{1pt}\otimes\hspace{1pt}\COMPAR}\rule{0pt}{0.15cm}\,\,\,\,\null} (P2)
(P1) edge node {\null\,\,\,\,\m{\COMPARR}} (P2);
\end{tikzpicture}
\end{center}\vspace{-3pt}
so that, for \m{\gamma\in G} and \m{v\in \TILDE\SIGMA_{k-2}%
},
\Algn{\label{eqf34ufh348gh348hf384h}
\COMPARR\left({\gamma}\sqbrackets{}{\INZB}{\QQp}{v}\right) & \textstyle= (\gamma\cdot v) \otimes \left(g \mapsto \delta(g\gamma)L(1)\right).%
}
 Let \m{\ZEROIDEAL} be the kernel of \m{\COMPARR}. %
  Let \m{\Theta_{k,a,\VAREPSILONp}} be image of the composition  %
 \Algn{\label{eq34g3g09g5i094gi409g}\ind_{\INZA}^{G}{\uline{\Sym}^{k-2}}({\overline{\ZZ}_p^2)}\otimes \tau \xhookrightarrow{\ \ } \ind_{\INZA}^{G}\TILDE\SIGMA_{k-2}\xtwoheadrightarrow{\ \ } \ind_{\INZA}^{G}\TILDE\SIGMA_{k-2}/\ZEROIDEAL.}
Let \m{\overline\Theta_{k,a,\varepsilon} = \Theta_{k,a,\varepsilon} \otimes \FFp}, and let \m{\INAAA  } be the kernel of the induced quotient map 
    \Algn{\label{eq34f3g3g3g45g445g45gg45}\ind_{\INZA}^G \SIGMA_{k-2} \xtwoheadrightarrow{\ \ } \overline\Theta_{k,a,\varepsilon}.}

   \subsection{\protect\boldmath The modulo~\texorpdfstring{\m{p}}{p} correspondence}\label{secBerger}
In~\cite{BB10}, Berger and Breuil prove that the locally algebraic representation
  \algn{\Pi(V_{k,a,\VAREPSILONp}) = \Alg(V_{k,a,\VAREPSILONp})\otimes \Lisse(V_{k,a,\VAREPSILONp}) = {\Sym^{k-2}}({\overline{\QQ}_p^2)} \otimes \absdet^{-\frac12} \otimes \Ind_B^G \RHO}
  is non-zero, irreducible, and admissible, by studying its \m{p}-adic completion \m{B(V_{k,a,\VAREPSILONp})}. 
  Let \m{\overline{\Pi}(V_{k,a,\VAREPSILONp})} be the semi-simplification of the reduction modulo  \m{\maxideal} of  a \m{G}-stable %
   lattice in \m{{{\Pi}(V_{k,a,\VAREPSILONp})}}. The following theorem follows from the main result in~\cite{Ber10}. %

          \begin{theorem}\label{BERGERX1}
          There are \m{t\in\{0,\ldots,p-1\}} and \m{\psi:\QQ_p^\times \to \overline\FF_p^\times} such that either 
          \algn{\textstyle \overline\Theta_{k,a,\varepsilon} \cong \repnn(t,0,\psi) %
          }
          or
          \algn{\textstyle \overline\Theta_{k,a,\varepsilon}^{\rm ss} \cong \left(\repnn(t,\lambda,\psi) \oplus\repnn(\LESS{p-3-t},\lambda^{-1},\omega^{t+1}\psi) \right)^{\rm ss}}
          for some \m{\lambda \in\FFptimes}.
          In the former case we have
          \algn{\textstyle \overline V_{k,a,\varepsilon} \cong \ind(\omega_2^{t+1})\otimes\psi,}
          and in the latter case we have
          \algn{\textstyle \overline V_{k,a,\varepsilon} \cong \left(\mu_\lambda \omega^{t+1} \oplus \mu_{\lambda^{-1}}\right)\otimes\psi.}
   \end{theorem}

\begin{proof}{Proof. }
 Galois representations that become semi-stable on an extension of \m{\QQ_p} by a \m{p}-power root of unity are trianguline (see section~2.3 of~\cite{Ber11}).
 So  theorem~3.1.1 in~\cite{Ber10} applies %
  and gives a correspondence between \m{\overline{V}_{k,a,\VAREPSILONp}} and \m{\overline{\Pi}(V_{k,a,\VAREPSILONp})}. 
Equations~\eqref{MAINEQN1} and~\eqref{eq34g3g09g5i094gi409g}
imply that \m{\Theta_{k,a,\VAREPSILONp}} is a lattice in \m{\Pi(V_{k,a,\VAREPSILONp})^{\mathrm{ss}}}. The correspondence we want follows from the fact that 
 \m{\overline{\Pi}(V_{k,a,\VAREPSILONp})} has finite length, %
 so  \m{\overline{\Pi}(V_{k,a,\VAREPSILONp})^{\rm ss} \cong \overline\Theta_{k,a,\VAREPSILONp}^{\rm ss}} %
  and therefore we can replace \m{\overline{\Pi}(V_{k,a,\VAREPSILONp})^{\rm ss}} with \m{\overline\Theta_{k,a,\VAREPSILONp}^{\rm ss}} in Berger's correspondence.
\end{proof}

\subsection{Reformulation of theorem~\ref{MAINTHM1}}
\label{subsRefor}
We can use theorem~\ref{BERGERX1}   to rephrase theorem~\ref{MAINTHM1} as an equivalent result about \m{G}-representations. For \m{l \in \ZZ}, let us define \m{\BIRREDUCIBLE{}{}{l}} as
    \algn{\BIRREDUCIBLE{}{}{l} = \repnn (\LESS{\NEBENPOWEXX-2l},0,\omega^{l}). %
    }
For \m{l \in \ZZ} and \m{\lambda\in \FFptimes}, let us define \m{\BRREDUCIBLE{l}(\lambda)} as
    \algn{\BRREDUCIBLE{l}(\lambda) = \left(\repnn(p-2,\lambda,\omega^{l}) \oplus\repnn(p-2,\lambda^{-1},\omega^{l}) \right)^{\rm ss}.}%
\begin{theorem}\label{MAINTHM2}
If \m{\TADICVALUATION(a) \in (0,\frac{p-1}{2})\backslash \ZZ}   then
    \Algn{\label{eqg30ig40it9304i90i}\overline {\Theta}_{k,a,\varepsilon}\in \left(\{\BPI\} \cup\left\{\BIRREDUCIBLE{k}{\NEBENPOWER}{\NEBENPOWEXX-j} %
    \right\}_{ 1\LEQSLANT j\LEQSLANT\TEMPORARYI-1}\right) \big\backslash \left\{\BIRREDUCIBLE{k}{\NEBENPOWER}{j} 
    \right\}_{ 1\LEQSLANT j\LEQSLANT\TEMPORARYI-1},}
where
    \algn{\BPI \cong \left\{\begin{array}{rl}\BRREDUCIBLE{2\TEMPORARYI-1}(\lambda)%
    & \text{if } \NEBENPOWEXX\equiv 2\TEMPORARYI-1 \bmod {p-1} \text{ \& } a\in\CIRCLES_{\TEMPORARYI} \text{,} \\[5pt] \BIRREDUCIBLE{k}{\NEBENPOWER}{\NEBENPOWEXX} & \text{otherwise,} \end{array}\right.}
and \m{\lambda+\lambda^{-1}} is the reduction modulo \m{\maxideal} of 
\algn{\textstyle  ((\TEMPORARYI-1)!)^{-1} (1-\VAREPSILONp(1+p))^{\TEMPORARYI-1} a^{-3} \left( \frac{1}{\TEMPORARYI!(\TEMPORARYI-1)!} (\VAREPSILONp(1+p)-1)^{2\TEMPORARYI-1}-a^2\right) .}
\end{theorem}

Therefore our goal is to prove theorem~\ref{MAINTHM2}. In the remainder of this article we make the assumption \m{\TADICVALUATION(a) \in (0,\frac{p-1}{2})\backslash\ZZ}, %
 so that \m{\TEMPORARYI\in\{1,\ldots,\frac{p-1}{2}\}} and \m{\TEMPORARYI-1<\PADICVALUATION(a)<\TEMPORARYI}. %
 In particular, \m{a^2 \ne p^{k-1}}.

    \section{Local constancy of \texorpdfstring{\m{V_{k,a,\varepsilon}}}{V}}

    In this section we show that it is enough to prove theorem~\ref{MAINTHM1} for \m{k>p^{100}}. We need the assumption \m{k>p^{100}} throughout sections~\ref{secsubs}, \ref{secseries}, \ref{SECPROOF}, \ref{secprop1}, and \ref{secprop2} so that we can construct sufficiently many elements of  \m{\ind_{\INZA}^G \TILDE\SIGMA_{k-2}}.

\begin{lemma}\label{CONSTANCY1}
{If theorem~\ref{MAINTHM1} is true for   \m{k>p^{100}} then it is true for all \m{k}.}
\end{lemma}
\begin{proof}{Proof. }
The representations \m{V_{k,a,\varepsilon}} belong to a certain family \m{V(\alpha,\beta)} of trianguline representations which is defined by Colmez in subsection~4.5 of~\cite{Col08}. Here \m{(\alpha,\beta)} belongs to a parameter space of pairs of characters \m{\QQ_p^\times \to \QQptimes}. If \m{x^{k-1}} denotes the character \m{\QQ_p^\times \to \QQptimes} which maps \m{z} to \m{z^{k-1}}, then, due to proposition~4.13 in~\cite{Col08}, \algn{\textstyle V_{k,a,\varepsilon} = V(\VAREPSILONp\mu_{ap^{1-k}}x^{k-1},\mu_{a^{-1}}).} Proposition~3.9 in~\cite{Che13} implies that  %
 \m{\overline V(\alpha,\beta)} is locally constant in %
 \m{(\alpha,\beta)}. %
For \m{m\in\ZZ_{\GEQSLANT1}}, let \m{k_m=k+(p-1)p^m}. Then 
 \m{\textstyle \mu_{ap^{1-k_m}}x^{{k_m}-1} \to \mu_{ap^{1-k}}x^{k-1}} as \m{m\to \infty}.
So there is an integer \m{m_0} such that \m{\overline{V}_{k,a,\varepsilon} \cong \overline{V}_{k_m,a,\varepsilon}} for all \m{m\GEQSLANT m_0}.
 In particular, \m{\overline{V}_{k,a,\varepsilon} \cong \overline{V}_{k_m,a,\varepsilon}} for some \m{m\GEQSLANT100}, and that completes the proof. %
\end{proof}

In the remainder of this article we make the assumption \m{k>p^{100}}.

\section{Notation}\label{SECnot}

 For \m{h\in\ZZ}, let {\m{I_h}} denote the %
 \m{\FFp [K Z]}-module of degree \m{h} homogeneous functions \m{{\FF_p^2} \to \FFp } that vanish at the origin, with \m{\begin{Smatrix}
    p & 0 \\ 0 & p
    \end{Smatrix}} acting trivially  and \m{K} acting by  \algn{\textstyle{}\left(\begin{Lmatrix}
g_1 & g_2 \\ g_3 & g_4
\end{Lmatrix}f\right)(x,y)=f(g_1x+g_3y,g_2x+g_4y)}  for \m{\begin{Smatrix}
g_1 & g_2 \\ g_3 & g_4
\end{Smatrix} \in K} and \m{f\in I_h}.
   Then \m{\sigma_{\MORE{h}} \subset I_h}, since any element of  \m{\sigma_{\MORE{h}}} is also a function \m{{\FF_p^2 \to \FFp}} which is homogeneous  of degree \m{h} and vanishes at the origin, and the actions of \m{KZ} on \m{\sigma_{\MORE{h}}} and \m{I_h} match. %
  Due to lemma~3.2 in~\cite{AS86}, there is a map
    \Algn{\label{EQUATION1}
    \textstyle f \in I_h \xmapsto{ \ \ } \sum_{u,v} f(u,v) (v  X-u Y)^{\LESS{-h}} \in \sigma_{\LESS{-h}}(h),%
    } which gives an isomorphism  \algn{\textstyle{}I_h/\sigma_{\MORE{h}} \cong \sigma_{\LESS{-h}}(h),
    }and therefore the only two factors of \m{I_h} are \m{\sigma_{\MORE{h}}} (``the submodule'') and \m{\sigma_{\LESS{-h}}(h)} (``the quotient'').
    By comparing the actions of \m{\begin{Smatrix}
    \lambda & 0 \\ 0 & \lambda
    \end{Smatrix}} on \m{\sigma_{\LESS{-h}}(h)} and \m{I_h}, we can conclude that \m{\sigma_{\LESS{-h}}(h)} is a submodule of \m{I_h}
    if and only if  \m{\MORE{h} = p-1}.  For \m{h\in\ZZ}, %
     let {\m{\chi_{h}}} denote the one-dimensional \m{\FFp[I(1)Z]}-module where  \m{{\begin{Smatrix} p & 0 \\ 0 & p \end{Smatrix}}} acts trivially and \m{\begin{Smatrix}
g_1 & g_2 \\ g_3 & g_4
\end{Smatrix} \in I(1)} acts as multiplication by \m{g_4^{h}}. Then \m{\chi_h} is also an \m{\FFp[I(m)Z]}-module for all \m{m\GEQSLANT1}. For \m{\alpha \in \ZZp}, let \m{\BIGO{\alpha}} denote the  sub-\m{\ZZp}-module
    \algn{\alpha \ind_{\INZA}^G(\uline{\Sym}^{k-2}(\overline{\ZZ}_p^2) \otimes \tau) \subseteq \ind_{\INZA}^G(\uline{\Sym}^{k-2}(\overline{\ZZ}_p^2)\otimes \tau).}
We abuse this notation and write \m{f=\BIGO{\alpha}} to represent a term \m{f \in \BIGO{\alpha}}, and we also write \m{h=\BIGO{\alpha}} when \m{h\in\alpha\ZZp}. If \m{f,g \in\ind_{\INZA}^G \TILDE\SIGMA_{k-2}} we write 
\algn{f\EQUIVZEROIDEAL g \Longleftrightarrow f-g\in\ZEROIDEAL.} %
 Let \m{\HECKET\in\End_G(\ind_{\INZA}^{G}\TILDE\SIGMA_{k-2})} be the Hecke operator corresponding to the double coset of \m{\begin{Smatrix}
p & 0 \\ 0 & 1
\end{Smatrix}}. There is the explicit formula  %
    \Algn{\label{eqHecke}\HECKET ({\gamma}\sqbrackets{}{\INZB}{\QQp}{v}) = {\sum_{\MU\in\FF_p} \gamma\begin{Lmatrix}
    p &[\MU] \\ 0 & 1
    \end{Lmatrix} }\sqbrackets{}{\INZB}{\QQp}{\left(\begin{Lmatrix}
    1 & -[\MU] \\ 0 & p
    \end{Lmatrix} \cdot v\right)},}
where \m{[\TEMPORARYA]} denotes the Teichm\"uller lift of \m{\TEMPORARYA \in \FF_p} to \m{\ZZ_p}. This formula can be derived in a similar fashion to the derivation of the formula for \m{T} in section~2.1 of~\cite{Bre03b}. 
For \m{H \in \{KZ\} \cup \{I(m)Z\,|\, m\in\ZZ_{\GEQSLANT1}\}}, \m{s\in \ZZ}, and an  \m{H}-module \m{W}, let \m{W(s)} denote the twist \m{W \otimes \det^s}.  If \m{A\subseteq B\subseteq\ind_{I(n)Z}^{G}\SIGMA_r}, we say that an element \m{x \in \ind_{I(n)Z}^{G}\SIGMA_r} represents an element \m{y \in B/A} if \m{x \in B} and \m{x+A = y \in B/A}. 
    Let \m{r=k-2}.

  \section{Combinatorial definitions}\label{seccombdef}
   
   By a ``set of lifts \m{\LCLASS}'' we mean a set of \m{p} elements of \m{\ZZ_p} such that the set  \m{\{\overline{y} \,|\, y \in \LCLASS\}} of the reductions modulo \m{p} of the elements of \m{\LCLASS} is equal to the whole of \m{\FF_p}.
   Let \m{\LCLASSTEICH} be the set of Teichm\"uller lifts.
    For \m{u,v\in \ZZ}, we use the convention that  \m{\binom{u}{v}} is equal to \m{ \frac1{v!}\prod_{l=0}^{v-1} (u-l)} if \m{v\GEQSLANT0} and equal to \m{0}  %
     if \m{v<0}.
     We also use the convention that when the range of summation is
not specified, it is assumed to be over all of
\m{\ZZ}, so \algn{\sum_{m_1,\ldots,m_k} F(m_1,\ldots,m_k) = \sum_{m_1 \in \ZZ} \cdots \sum_{m_k \in \ZZ}F(m_1,\ldots,m_k).} In all such sums appearing in this article, \m{F} is  supported on a finite subset of \m{\ZZ^k}, so this is well-defined.
 For \m{\tempA\in\ZZ} and \m{m,\tempB\in\ZZ_{\GEQSLANT 0}}, let \m{\MATRIX_{m}(\tempA,\tempB)} denote the \m{(m+1)\times (m+1)} matrix with entries in \m{\ZZ} indexed by \m{0\LEQSLANT i,j\LEQSLANT m} and defined by
    \algn{\textstyle \MATRIX_{m}(\tempA,\tempB) = \left(\binom{\tempA+i}{\tempB+j}\right)_{0\LEQSLANT i,j \LEQSLANT m}.}
For a boolean \m{P} let us define \m{[P]=1} if \m{P} is true and \m{[P]=0} if \m{P} is false. Let \algn{\GEN  & \textstyle{} ={} x^r-2x^{r-p+1}y^{p-1}+x^{r-2p+2}y^{2p-2} \in \ind_{\INZA}^G \TILDE\SIGMA_{k-2}.} For \m{\alpha\in\ZZ_{\GEQSLANT0}}, let \algn{ \GEN_{\alpha} & \textstyle{} ={} \sum_{\MU\in \FF_p} [\MU]^{\alpha} \begin{Lmatrix} 1 & [\MU] \\ 0 & 1
    \end{Lmatrix} \GEN,}
    and for \m{\alpha\in\ZZ_{<0}}, let \m{\GEN_\alpha=0}.
For \m{\alpha \in \{0,\ldots,p-1\}}, \m{s\in \ZZ_{\GEQSLANT1}}, and \m{\TEMPORARYAAAA\in\ZZ_p}, let %
\algn{\textstyle \CONST{\alpha}{s}_{\TEMPORARYAAAA} ={} \sum_{\MU\in\FF_p}[\MU]^{p-1-\alpha} \VAREPSILONp^{-1}\left(1 - \TEMPORARYAAAA[\MU]p^{s}\right).}
 For \m{s\in\ZZ_{\GEQSLANT0}}, let
\m{\mathrm{err}_s = \QUARK^{p^{s+1}} + ap}. 
For \m{\TEMPORARYA\in\FF_p}, let  \algn{\DIFF(\TEMPORARYA)=([\TEMPORARYA]+1-[\TEMPORARYA+1])/p \in \ZZ_p.}
 Let \m{\VECTORSP} denote the \m{\FF_p}-vector space with basis elements \m{\ZZ_{\GEQSLANT0} \times\{1,\ldots,p-1\}}. 
    For \m{(z,t)\in\{0,\ldots,p-1\}\times\{0,\ldots,p-2\}}, let \m{\Xi_{z,t} : \VECTORSP \to \FF_p} denote the linear map which sends the basis element \m{(u,s) \in \ZZ_{\GEQSLANT0} \times\{1,\ldots,p-1\}} to
    \algn{\Xi_{z,t}(u,s)=
     \left\{\begin{array}{rl}\vphantom{(-1)^{u} \sum_{\beta>0 \,\&\, \beta \equiv z \equiv t+u-s+\NEBENPOWEX \bmod p-1} \binom{u}{\beta}}  (-1)^{u} \sum_{\beta>0 \,\&\, \beta \equiv z \equiv t+u-s+\NEBENPOWEX \bmod p-1} \binom{u}{\beta} & \text{if } z>0\text{,} \\[4pt] \vphantom{(-1)^{u} \sum_{\beta>0 \,\&\, \beta \equiv z \equiv t+u-s+\NEBENPOWEX \bmod p-1} \binom{u}{\beta}}(-1)^{u}\VERIFYIF{t\equiv -u+s-\NEBENPOWEY \bmod p-1} & \text{if } z=0\text{.} %
     \end{array}\right.}
     For \m{\DELTAAA \in \ZZ}, let
     \algn{S_\DELTAAA =\{y \,|\, y \in \{0,\ldots,p-1\} \text{ \& } y \equiv \DELTAAA-w \bmod p-1 \text{ for some } w \in \{1,\ldots,\TEMPORARYI-1\}\}.}
     For \m{m\in\ZZ_{\GEQSLANT1}},  \m{\DELTAAA \in \ZZ}, and \m{(u_i)_{1\LEQSLANT i\LEQSLANT m} \in\{0,\ldots,p-1\}^m}, let
 \m{\MMM{m}{(u_i)_{1\LEQSLANT i\LEQSLANT m}}{\DELTAAA}} be the matrix with rows indexed by \m{i \in \{1,\ldots,m\}} and columns indexed by \m{j \in S_\DELTAAA}, with entries in \m{\FF_p} defined by 
\algn{\textstyle \MMM{m}{(u_i)_{1\LEQSLANT i\LEQSLANT m}}{\DELTAAA}= \left(\binom{u_i}{j}\right)_{1\LEQSLANT i\LEQSLANT m,\,j \in S_\DELTAAA}.}

\section{Reference tables}
The following table collates the key notation.

\newcommand{\NEWLINEAA}{\\[4pt]}
\begin{center}
\begin{tabular}{r|l}
    \m{n,k,r} & \m{n\in \ZZ_{\GEQSLANT2}}, \m{k\in\ZZ_{\GEQSLANT2}} is the weight, and \m{r = k-2}; we assume \m{k>p^{100}}. \NEWLINEAA
    \m{\TADICVALUATION} & \m{\TADICVALUATION = p^{n-2}(p-1) \cdot \PADICVALUATION}. \NEWLINEAA
    \m{\QUARK} & \m{\QUARK= p^{1/p^{n-2}(p-1)}}, so that \m{\TADICVALUATION(\QUARK)=1}. \NEWLINEAA
        \m{a} & the eigenvalue; we assume %
 \m{\TADICVALUATION(a) \in (0,\frac{p-1}{2})\backslash\ZZ}. \NEWLINEAA
      \m{\TEMPORARYI} &  \m{\TEMPORARYI = \lceil \TADICVALUATION(a)\rceil \in \{1,\ldots,\frac{p-1}{2}\}}. \NEWLINEAA
 \m{\VAREPSILONp} &  a  ramified character \m{ \ZZ_p^\times\to \QQptimes} which has conductor \m{p^n}.\NEWLINEAA \m{\ZETA_s,\ZETA} &  \m{\ZETA_s=\VAREPSILONp^{-1}(1-p^{n-s})} (\m{s \in \{1,\ldots,n-1\}}),  \m{\ZETA=\ZETA_{n-1}}.\NEWLINEAA
     \m{\NEBENPOWER,\NEBENPOWEX} &  the reduction modulo \m{p} of \m{\VAREPSILONp} is the map \m{\ast \mapsto {\ast^{\NEBENPOWER}}}, and \m{\NEBENPOWEX=r+\NEBENPOWER}.
    \NEWLINEAA
     \m{\MORE{h}} & the number in \m{\{1,\ldots,p-1\}} that is congruent to \m{h} modulo \m{p-1}. \NEWLINEAA
\m{\LESS{h}} & the number in \m{\{0,\ldots,p-2\}} that is congruent to \m{h} modulo \m{p-1}. \NEWLINEAA
\end{tabular}
\end{center}

The following table references the key definitions.

\begin{center}
\begin{tabular}{r|l}
\m{G}, \m{B}, \m{K},  \m{I(m)}, \m{\RHO} &   subsection~\ref{subsecWe}. \NEWLINEAA
 \m{\ind_{H}^G}, \m{\sqbrackets{g}{G,H}{\FF}{w}}, \m{\tau},  \m{\textstyle\TILDE\SIGMA_{k-2}},  \m{\textstyle\SIGMA_{k-2}}, \m{\sigma_h}, \m{\repnn(t,\lambda,\psi)} &   subsection~\ref{loccomp}. \NEWLINEAA
 \m{\Psi}, \m{\ZEROIDEAL}, \m{\Theta_{k,a,\VAREPSILONp{}}}, \m{\INAAA  } &   subsection~\ref{subcompar}. \NEWLINEAA
\m{I_h}, \m{\chi_h}, \m{\BIGO{\alpha}}, \m{\EQUIVZEROIDEAL}, \m{\HECKET} &   section~\ref{SECnot}.  \NEWLINEAA \m{\LCLASS}, \m{\LCLASSTEICH}, \m{\MATRIX_m}, \m{[\ ]}, \m{\GEN}, \m{\GEN_\alpha}, \m{\CONST{\alpha}{s}_{\TEMPORARYAAAA}}, \m{\mathrm{err}_s}, \m{\DIFF}, \m{\VECTORSP}, \m{\Xi_{z,t}}, \m{S_\DELTAAA}, \m{M_m} &   section~\ref{seccombdef}.
\end{tabular}
\end{center}

\section{Combinatorial identities}
\label{seccomb}

In this section we prove five lemmas involving standard combinatorial identities.

\begin{lemma}\label{LEMMAX8}
Suppose that \m{\tempA\in \ZZ} and \m{m,\tempB\in\ZZ_{\GEQSLANT 0}}.\AFIOFJFEWF
\begin{enumerate}[leftmargin=*]
    \item\label{vander1} We have \Algn{\label{eqvandermonde}\textstyle \det\MATRIX_{m}(\tempA,\tempB) =  \prod_{(l,w)\in\{0,\ldots,m\}\times\{0,\ldots,\tempB-1\}} \frac{\tempA+l-w}{\tempB+l-w}.} If \m{\tempB=0} then the product is empty and hence \m{\det\MATRIX_{m}(\tempA,\tempB)=1}. \AFIOFJ
    \item \label{vander2}If \m{v+m<p} and \m{p\nmid \binom{u+m}{v+m}} then \m{\MATRIX_{m}(u,v)}  is invertible over \m{\ZZ_p}.\AFIOFJ
\end{enumerate}
\end{lemma}

\begin{proof}{Proof. } \eqref{vander1}
If \m{\tempB>0}, we have \m{\binom{\tempA+i}{\tempB+j} = \frac{\tempA+i}{\tempB+j}\binom{\tempA+i-1}{\tempB+j-1}} for all \m{i,j\GEQSLANT 0} and therefore %
\algn{\textstyle \det  \MATRIX_{m}(\tempA,\tempB)  =  \det  \MATRIX_{m}(\tempA-1,\tempB-1)\prod_{l\in\{0,\ldots,m\}} \frac{u+l}{v+l}.}
We can use this equation and induction to prove that, for all \m{v \in \ZZ_{\GEQSLANT0}},
\Algn{\label{eq3453}\textstyle \det  \MATRIX_{m}(\tempA,\tempB)  =  \det  \MATRIX_{m}(\tempA-\tempB,0)\prod_{(l,w)\in\{0,\ldots,m\}\times\{0,\ldots,\tempB-1\}} \frac{\tempA+l-w}%
{\tempB+l-w}.}
 Due to Vandermonde's convolution formula, for all \m{i,j\in\{0,\ldots,m\}} we have
\algn{\textstyle \sum_{l} \binom{i}{l}\binom{\tempA-\tempB}{j-l} = \binom{\tempA-\tempB+i}{j},} and therefore
\algn{\textstyle \MATRIX_{m}(\tempA-\tempB,0) = \left(\binom{\tempA-\tempB+i}{j}\right)_{0\LEQSLANT i,j \LEQSLANT m} = \left(\binom{i}{j}\right)_{0\LEQSLANT i,j \LEQSLANT m} \left(\binom{\tempA-\tempB}{j-i}\right)_{0\LEQSLANT i,j \LEQSLANT m}.}
This is the product of a lower triangular matrix with \m{1}'s on the diagonal and an upper triangular matrix with \m{1}'s on the diagonal. So
\m{%
\det \MATRIX_{m}(\tempA-\tempB,0)  = %
1}
and equation~\eqref{eqvandermonde} follows from equation~\eqref{eq3453}.
\FUEWIFFEW
\eqref{vander2} The assumptions \m{v+m<p} and \m{p\nmid %
 \binom{u+m}{v+m}} imply that %
  \m{\frac{\tempA+l-w}%
{\tempB+l-w}} is a unit for all \m{(l,w)\in\{0,\ldots,m\}\times\{0,\ldots,\tempB-1\}}. So equation~\eqref{eqvandermonde}  implies that
\algn{ \det\MATRIX_{m}(u,v)\in\ZZ_p^\times.} %
\end{proof}

\begin{lemma}\label{LEMMAX3.04}
Recall that, for \m{s\in\ZZ_{\GEQSLANT0}}, \m{\mathrm{err}_s = \QUARK^{p^{s+1}} + ap}, and  \m{\ZETA =\varepsilon^{-1}(1-p)}. 
\begin{enumerate}[leftmargin=*]
    \item\label{comb10.1} For \m{s\in \ZZ_{\GEQSLANT0}}, \m{\BIGO{\mathrm{err}_sa^{-1}} = \BIGO{\QUARK^{p^{s+1}}a^{-1}} \cup \BIGO{p} \supseteq \BIGO{p}}. \AFIOFJ
    \item\label{comb10.2} For \m{\TEMPORARYA\in\FF_p}, \m{\DIFF(\TEMPORARYA) = \sum_{j=1}^{p-1}\frac{1}{j}[-\TEMPORARYA]^j+\BIGO{p}}. \AFIOFJ
    \item\label{comb10.3} For \m{\TEMPORARYA,\MU\in\FF_p},
    \algn{\sum_{j=1}^{p-1}\frac1j([\MU]+[\TEMPORARYA])^j[\TEMPORARYA]^{\TEMPORARYI-j}& \textstyle= \sum_{l=1}^{p-1} \frac{(-1)^l}{l} [\MU]^l [\TEMPORARYA]^{\TEMPORARYI-l}+\BIGO{p}.}
    \item \label{comb10.5} For \m{u\in\ZZ} and \m{v,w\in\ZZ_{\GEQSLANT0}}, %
    \algn{\sum_{i} (-1)^{i}\binom{v}{i}\binom{u+i}{w} = (-1)^v \binom{u}{w-v}.} 
    \item\label{comb10.4}  For \m{j\in\{1,\ldots,p-1\}},
    \algn{\sum_{i=1}^{\TEMPORARYI-1} \frac{(-1)^{i}}{i}\binom{\TEMPORARYI-1}{i} \binom{p-i}{j} = \frac{1}{j}\left((-1)^{j+1}-\binom{-\TEMPORARYI}{j-\TEMPORARYI}\right)+ \BIGO{p}.}
    \end{enumerate}
\end{lemma}

\begin{proof}{Proof. }
 Part~\eqref{comb10.1} follows because the assumption \m{\TADICVALUATION(a)\not\in\ZZ} implies that  the two  valuations \m{\TADICVALUATION(\QUARK^{p^{s+1}})} and \m{\TADICVALUATION (ap)} are never the same. For part~\eqref{comb10.2},
 \algn{\textstyle [\TEMPORARYA+1]^p & \textstyle = ([\TEMPORARYA]+1-p\DIFF(\TEMPORARYA))^p \\ & \textstyle = ([\TEMPORARYA]+1)^p + \BIGO{p^2} = [\TEMPORARYA]^p+1+ \sum_{j=1}^{p-1}\binom{p}{j}[\TEMPORARYA]^j+\BIGO{p^2},}
and therefore
\algn{\nonumber\textstyle \DIFF(\TEMPORARYA) & =\textstyle ([\TEMPORARYA]^p+1-[\TEMPORARYA+1]^p)/p \\\nonumber& =-\textstyle \sum_{j=1}^{p-1}\frac{1}{p}\binom{p}{j}[\TEMPORARYA]^j + \BIGO{p} \\
& =-\textstyle \sum_{j=1}^{p-1}\frac{1}{j}\binom{p-1}{j-1}[\TEMPORARYA]^j + \BIGO{p} %
 =
 \textstyle \sum_{j=1}^{p-1}\frac{1}{j}[-\TEMPORARYA]^j+\BIGO{p}.} Part~\eqref{comb10.3} follows from
\algn{\sum_{j=1}^{p-1}\frac1j([\MU]+[\TEMPORARYA])^j[\TEMPORARYA]^{\TEMPORARYI-j}& \textstyle=\vphantom{\sum_{l=0}^{p-1} \left(\sum_{j=1}^{p-1}\frac1j\binom{j}{l}\right) [\MU]^l [\TEMPORARYA]^{\TEMPORARYI-l}}\sum_{j=1}^{p-1}\sum_{l=0}^{j} \frac1j\binom{j}{l} [\MU]^l [\TEMPORARYA]^{\TEMPORARYI-l}\\ & \textstyle\vphantom{\sum_{l=0}^{p-1} \left(\sum_{j=1}^{p-1}\frac1j\binom{j}{l}\right) [\MU]^l [\TEMPORARYA]^{\TEMPORARYI-l}}=\sum_{l=0}^{p-1} \left(\sum_{j=1}^{p-1}\frac1j\binom{j}{l}\right) [\MU]^l [\TEMPORARYA]^{\TEMPORARYI-l}
            \\ & \textstyle=\vphantom{\sum_{l=0}^{p-1} \left(\sum_{j=1}^{p-1}\frac1j\binom{j}{l}\right) [\MU]^l [\TEMPORARYA]^{\TEMPORARYI-l}}\sum_{l=1}^{p-1} \frac{(-1)^l}{l} [\MU]^l [\TEMPORARYA]^{\TEMPORARYI-l}+\BIGO{p},}
            where the third equality follows from the identities \m{\sum_{j=1}^{p-1}\frac1j = \BIGO{p}}, \m{\frac1j \binom{j}{l} = \frac{1}{l}\binom{j-1}{l-1}} for \m{l\in\ZZ_{\GEQSLANT1}}, and \algn{\sum_{j=1}^{p-1} \binom{j-1}{l-1} = \binom{p-1}{l} = (-1)^l+\BIGO{p}}
               for \m{l\in \ZZ_{\GEQSLANT1}}.     Part~\eqref{comb10.5} %
             follows from the identity \algn{\textstyle \sum_{w\GEQSLANT0}\sum_{i} (-1)^{i}\binom{v}{i}\binom{u+i}{w}z^w = (1+z)^u\sum_{i} (-1)^{i}\binom{v}{i}(1+z)^i=(-z)^v (1+z)^u\in\QQ(z).}
 Part~\eqref{comb10.4} follows from the identities
  \newcommand{\tempheightaer}{\vphantom{\frac{1}{j}\left((-1)^{j+1}+\binom{-\TEMPORARYI}{j-\TEMPORARYI}\right)}}\algn{%
    \textstyle\nonumber \sum_{i=1}^{\TEMPORARYI-1} \frac{(-1)^{i}}{i}\binom{\TEMPORARYI-1}{i} \binom{p-i}{j} & \textstyle =\nonumber\tempheightaer \sum_{i=1}^{\TEMPORARYI-1} \frac{(-1)^i}{i}\binom{\TEMPORARYI-1}{i} \binom{-i}{j}+\BIGO{p} 
                \\
                & \textstyle\nonumber =\tempheightaer \frac{1}{j} \sum_{i=1}^{\TEMPORARYI-1}(-1)^{i+1} \binom{\TEMPORARYI-1}{i} \binom{-i-1}{j-1}+\BIGO{p} 
                 \\
                & \textstyle =\nonumber\tempheightaer \frac{(-1)^{j+1}}{j} + \frac{1}{j} \sum_{i=0}^{\TEMPORARYI-1}(-1)^{i+1} \binom{\TEMPORARYI-1}{i} \binom{-i-1}{j-1}+\BIGO{p}
                 \\
                & \textstyle =\tempheightaer \frac{1}{j}\left((-1)^{j+1}-\binom{-\TEMPORARYI}{j-\TEMPORARYI}\right)+ \BIGO{p}
            ,}
          where the fourth equality follows from part~\eqref{comb10.5} with \m{u=w=j-1} and \m{v=\TEMPORARYI-1}.  %
\end{proof}

\begin{lemma}\label{LEMMAX3.01}
Recall that \m{\GEN=x^r-2x^{r-p+1}y^{p-1}+x^{r-2p+2}y^{2p-2} \in \ind_{\INZA}^G \TILDE\SIGMA_{r}}.\AFIOFJFEWF
 \begin{enumerate}[leftmargin=*]
  \item\label{propeta1}  \m{\GEN = (x^{p-1}-y^{p-1})^{2}x^{r-2p+2}}. \AFIOFJ 
 \item\label{propeta2} \m{\GEN(x,p y) = x^{r} + \BIGO{p^2}}.\AFIOFJ 
\item\label{propeta3} If \m{\alpha\in\{1,\ldots,p-1\}} then \algn{\textstyle\GEN_\alpha = (-1)^{\alpha-1} (\alpha-1) x^{r-\alpha}y^\alpha + (-1)^\alpha \alpha x^{r-\alpha-p+1}y^{\alpha+p-1}+\BIGO{p}.}
\item\label{propeta4} \m{\GEN = \GEN_0 - \GEN_{p-1}}.\AFIOFJ 
\item\label{propeta5} \m{\GEN_0 = x^{r} + \BIGO{p}}.\AFIOFJ
    \item\label{propeta6}  For \m{\MU \in \FF_p},
     \algn{\textstyle \begin{Lmatrix}
    1 & [\MU] \\ 0 & 1
    \end{Lmatrix} \GEN = (1-[\MU]^{p-1})\left(\GEN_0-\frac p{p-1}\GEN_{p-1}\right) + \frac{1}{p-1} \sum_{\alpha=1}^{p-1} [\MU]^{p-1-\alpha} \GEN_\alpha.}
\end{enumerate}
\end{lemma}
 
\begin{proof}{Proof. }
Parts~\eqref{propeta1} and~\eqref{propeta2} follow from the definition of \m{\GEN}. If \m{\beta\in\ZZ_{\GEQSLANT1}} then
\algn{\textstyle \sum_{\MU\in\FF_p} [\MU]^\beta = \left\{\begin{array}{rl} p-1 & \text{if }p-1\mid\beta\text{,} \\ 0 & \text{otherwise.} \end{array}\right.} Consequently, if \m{\alpha\in\{1,\ldots,p-1\}} then
\algn{\textstyle
\GEN_\alpha & = \textstyle \sum_{\MU\in\FF_p} [\MU]^\alpha \begin{Lmatrix}
1 & [\MU] \\ 0 & 1
\end{Lmatrix} (x^r-2x^{r-p+1}y^{p-1}+x^{r-2p+2}y^{2p-2})
\\ & = \textstyle \sum_{\MU\in\FF_p} [\MU]^\alpha  (x^r-2x^{r-p+1}([\MU]x+y)^{p-1}+x^{r-2p+2}([\MU]x+y)^{2p-2})
\\ & = \textstyle b_{0} x^{r-\alpha}y^\alpha + b_{1} x^{r-\alpha-p+1}y^{\alpha+p-1}}
with \algn{\textstyle b_0 &\textstyle= -2(p-1)\binom{p-1}{\alpha}+(p-1)\binom{2p-2}{\alpha} = (-1)^\alpha(1-\alpha) + \BIGO{p},\\ b_1 & \textstyle = (p-1)\binom{2p-2}{\alpha+p-1} = -\binom{1}{1}\binom{p-2}{\alpha-1} + \BIGO{p} = (-1)^{\alpha} \alpha + \BIGO{p}.}
This implies part~\eqref{propeta3}. %
Part~\eqref{propeta4} follows from the equation
\algn{\textstyle \GEN_{0}-\GEN_{p-1} = \sum_{\MU\in \FF_p} ([\MU]^{0}-[\MU]^{p-1})\begin{Lmatrix}
    1 & [\MU] \\ 0 & 1
    \end{Lmatrix} \GEN = \begin{Lmatrix}
    1 & 0 \\ 0 & 1
    \end{Lmatrix} \GEN=\GEN.}
Part~\eqref{propeta5} follows from parts~\eqref{propeta3} and~\eqref{propeta4}. If  \m{\MU=0} then part~\eqref{propeta6} is equivalent to part~\eqref{propeta4}, and if \m{\MU\ne 0} then part~\eqref{propeta6} follows from the equations \m{1-[\MU]^{p-1} = 0} and
\algn{ \textstyle
\sum_{\alpha=1}^{p-1} [\MU]^{p-1-\alpha} \GEN_\alpha & \textstyle \vphantom{\sum_{\TEMPORARYAA\in\FF_p} {\left(\sum_{\alpha=1}^{p-1}[\TEMPORARYAA]^\alpha\right)} \begin{Lmatrix}
1 & [\TEMPORARYAA] [\MU] \\ 0 & 1
\end{Lmatrix} \GEN }= \sum_{\alpha=1}^{p-1} [\MU]^{-\alpha} \sum_{\LAMBDA\in\FF_p} [\LAMBDA]^\alpha \begin{Lmatrix}
1 & [\LAMBDA] \\ 0 & 1
\end{Lmatrix} \GEN 
\\ & \textstyle \vphantom{\sum_{\TEMPORARYAA\in\FF_p} {\left(\sum_{\alpha=1}^{p-1}[\TEMPORARYAA]^\alpha\right)} \begin{Lmatrix}
1 & [\TEMPORARYAA] [\MU] \\ 0 & 1
\end{Lmatrix} \GEN }= \sum_{\alpha=1}^{p-1} \sum_{\LAMBDA\in\FF_p} [\LAMBDA/\MU]^\alpha \begin{Lmatrix}
1 & [\LAMBDA/\MU] [\MU] \\ 0 & 1
\end{Lmatrix} \GEN 
\\ & \textstyle \vphantom{\sum_{\TEMPORARYAA\in\FF_p} {\left(\sum_{\alpha=1}^{p-1}[\TEMPORARYAA]^\alpha\right)} \begin{Lmatrix}
1 & [\TEMPORARYAA] [\MU] \\ 0 & 1
\end{Lmatrix} \GEN }= \sum_{\TEMPORARYA\in\FF_p} {\left(\sum_{\alpha=1}^{p-1}[\TEMPORARYA]^\alpha\right)} \begin{Lmatrix}
1 & [\TEMPORARYA] [\MU] \\ 0 & 1
\end{Lmatrix} \GEN 
\\ & \textstyle \vphantom{\sum_{\TEMPORARYAA\in\FF_p} {\left(\sum_{\alpha=1}^{p-1}[\TEMPORARYAA]^\alpha\right)} \begin{Lmatrix}
1 & [\TEMPORARYAA] [\MU] \\ 0 & 1
\end{Lmatrix} \GEN }= (p-1) \begin{Lmatrix}
1 & [\MU] \\ 0 & 1
\end{Lmatrix} \GEN.
}
Here for the third equality we change the variable \m{\LAMBDA} to \m{\TEMPORARYA=\LAMBDA/\MU}, and for the fourth equality we note that
\algn{\textstyle \sum_{\alpha=1}^{p-1}[\TEMPORARYA]^\alpha = \left\{\begin{array}{rl} p-1 & \text{if }\TEMPORARYA=1\text{,} \\ 0 & \text{otherwise.} \end{array}\right.}
\end{proof}

\begin{lemma}\label{LEMMAX3.02}
{Suppose that \m{\TEMPORARYAAAA\in \ZZ_p^\times}. \begin{enumerate}[leftmargin=*]
    \item \label{gammaprop1} If \m{\alpha \in \{0,\ldots,p-1\}} and \m{s\in\ZZ_{\GEQSLANT n}} then \algn{\textstyle \CONST{\alpha}{s}_{\TEMPORARYAAAA} = \sum_{\MU\in\FF_p}[\MU]^{p-1-\alpha} = \left\{\begin{array}{rl} p-1 & \text{if }\alpha= 0\text{,} \\ 0 & \text{if }0 < \alpha < p-1\text{,} \\ p & \text{if }\alpha= p-1\text{.} \end{array}\right.}
    \item \label{gammaprop2}If \m{(\alpha,s)=(p-1,n-1)} then \algn{\CONST{\alpha}{s}_{\TEMPORARYAAAA}=\sum_{\MU\in\FF_p}  \varepsilon^{-1}(1-\TEMPORARYAAAA [\MU]p^{n-1}) = 0.}
    \item \label{gammaprop3}If \m{\alpha\in\{0,\ldots,p-1\}} and \m{s \in \{1,\ldots,n-1\}} and \m{(\alpha,s)\ne(p-1,n-1)} then %
 \algn{\textstyle\CONST{\alpha}{s}_{\TEMPORARYAAAA} = \left\{\begin{array}{rl} \frac{(-1)^{\alpha+1}}{\alpha!}  (1-\ZETA_{\TEMPORARYAAAA}^p)(1-\ZETA_{\TEMPORARYAAAA})^{\alpha-p} (1+\BIGO{\QUARK})  & \text{if } s<n-1\text{,} \\ \frac{(-1)^{\alpha}}{\alpha!}  p(1-\ZETA_{\TEMPORARYAAAA})^{\alpha-p+1} (1+\BIGO{\QUARK}) & \text{if } s=n-1\text{,} \end{array}\right.} %
where \m{\ZETA_{\TEMPORARYAAAA}=\VAREPSILONp^{-1}\left(1-\TEMPORARYAAAA p^{s} \right)}.
\end{enumerate}  }
\end{lemma}
\begin{proof}{Proof. } Parts~\eqref{gammaprop1} and~\eqref{gammaprop2} follow from the definition of \m{\CONST{\alpha}{s}_{\TEMPORARYAAAA}}. %

\eqref{gammaprop3} For  \m{\beta\in\{0,\ldots,p-1\}}, let
\algn{\textstyle \delta_\beta = \delta_\beta(\ZETA_\TEMPORARYAAAA) = \left\{\begin{array}{rl} \frac{(-1)^{\beta+1}}{\beta!}  (1-\ZETA_{\TEMPORARYAAAA}^p)(1-\ZETA_{\TEMPORARYAAAA})^{\beta-p}    & \text{if } s<n-1\text{,} \\ \frac{(-1)^{\beta}}{\beta!}  p(1-\ZETA_{\TEMPORARYAAAA})^{\beta-p+1}  & \text{if } s=n-1\text{,} \end{array}\right.}
and
\algn{f_\beta(z) = \sum_{m=0}^{p-1} m^{p-1-\beta}z^m \in \ZZ_p[z].} We want to prove the equation  \Algn{\label{eq2353t3}\CONST{\alpha}{s}_{\TEMPORARYAAAA} = \delta_\alpha%
(1 + \BIGO{\QUARK}).}
Note that, since \m{(\alpha,s)\ne(p-1,n-1)}, \Algn{\label{eq34t3g3g4g545g5g4}\TADICVALUATION(\delta_\alpha ) = \alpha\TADICVALUATION(1-\ZETA_\TEMPORARYAAAA) \LEQSLANT \min\{\TADICVALUATION(1-\ZETA_\TEMPORARYAAAA^p),\TADICVALUATION(p)\}-1.} Moreover,
\algn{\VAREPSILONp^{-1}\left(1 - \TEMPORARYAAAA[\MU]p^{s}\right) = \ZETA_\TEMPORARYAAAA^{[\MU]}+\BIGO{1-\ZETA_\TEMPORARYAAAA^p},} %
and consequently
\algn{\CONST{\alpha}{s}_{\TEMPORARYAAAA} \textstyle{} = \sum_{\MU\in\FF_p}[\MU]^{p-1-\alpha} \ZETA_\TEMPORARYAAAA^{[\MU]} + \BIGO{1-\ZETA_\TEMPORARYAAAA^p} %
 = f_\alpha( \ZETA_\TEMPORARYAAAA) + \BIGO{1-\ZETA_\TEMPORARYAAAA^p} + \BIGO{p}.}
So in order to prove equation~\eqref{eq2353t3}, it is enough to prove the equivalent equation
\Algn{\label{eq54t4grth}\textstyle{}  f_\alpha( \ZETA_\TEMPORARYAAAA) = \delta_\alpha %
(1+\BIGO{\QUARK}).}
We have
\Algn{\label{eqr43r3g45g4h}
    \textstyle (1-\ZETA_\TEMPORARYAAAA)f_{\alpha}(\ZETA_\TEMPORARYAAAA) &\textstyle \vphantom{\sum_{m=0}^{p-1} (m^{p-1-\alpha}-(m-1)^{p-1-\alpha})\ZETA_\TEMPORARYAAAA^m +  \BIGO{1-\ZETA_\TEMPORARYAAAA^p} + \BIGO{p}}= 
    \sum_{m=0}^{p-1} m^{p-1-\alpha}\ZETA_\TEMPORARYAAAA^m-   \sum_{m=1}^{p}(m-1)^{p-1-\alpha}\ZETA_\TEMPORARYAAAA^{m} 
    \nonumber
    \\ &\textstyle \vphantom{\sum_{m=0}^{p-1} (m^{p-1-\alpha}-(m-1)^{p-1-\alpha})\ZETA_\TEMPORARYAAAA^m +  \BIGO{1-\ZETA_\TEMPORARYAAAA^p} + \BIGO{p}}= 
    \sum_{m=0}^{p-1} (m^{p-1-\alpha}-(m-1)^{p-1-\alpha})\ZETA_\TEMPORARYAAAA^m + (-1)^{\alpha} -(p-1)^{p-1-\alpha}\ZETA_\TEMPORARYAAAA^p %
    \nonumber
    \\ &\textstyle \vphantom{\sum_{m=0}^{p-1} (m^{p-1-\alpha}-(m-1)^{p-1-\alpha})\ZETA_\TEMPORARYAAAA^m +  \BIGO{1-\ZETA_\TEMPORARYAAAA^p} + \BIGO{p}}=\sum_{\beta=\alpha+1}^{p-1} b_\beta f_\beta(\ZETA_\TEMPORARYAAAA) + \BIGO{1-\ZETA_\TEMPORARYAAAA^p} + \BIGO{p},
}
where \m{b_{\alpha+1},\ldots,b_{p-1}\in\ZZ_p}
are the coefficients of the polynomial
\algn{b_{\alpha+1}z^{p-2-\alpha} + \cdots + b_{p-1} = z^{p-1-\alpha}-(z-1)^{p-1-\alpha} \in \ZZ_p[z],}
so that \m{b_{\alpha+1}=p-1-\alpha}.
   Suppose that \m{s<n-1}. Then equation~\eqref{eq54t4grth}  follows for \m{\alpha=p-1} from   \algn{f_{p-1}(z) \textstyle{}=\frac{1-z^{p}}{1-z}.}
If \m{\alpha\in\{0,\ldots,p-2\}}, and if \m{f_\beta(\ZETA_\TEMPORARYAAAA)=\delta_\beta(1+\BIGO{\QUARK})} for all \m{\beta\in\{\alpha+1,\ldots,p-1\}}, then 
equation~\eqref{eq54t4grth} follows from equations~\eqref{eq34t3g3g4g545g5g4} and~\eqref{eqr43r3g45g4h}. Therefore we can prove by induction that equation~\eqref{eq54t4grth} is true for all \m{\alpha\in\{0,\ldots,p-1\}}. %
 Suppose that \m{s=n-1}. Then equation~\eqref{eq54t4grth}  follows for \m{\alpha=p-2} from  \algn{f_{p-2}(z)  \textstyle{}= z\left(\frac{\partial}{\partial z}f_{p-1}(z)\right) = -\frac{pz^p}{1-z} - \frac{z(1-z^p)}{(1-z)^2},} 
 and we can similarly use equation~\eqref{eqr43r3g45g4h} to prove by induction that equation~\eqref{eq54t4grth} is true for all \m{\alpha\in\{0,\ldots,p-2\}}.
\end{proof}

\begin{lemma}\label{LEMMAX3.03}
Recall that \m{\ZETA =\varepsilon^{-1}(1-p)}. If \m{\alpha\in\{1,\ldots,p-1\}} then
\algn{ \textstyle  \sum_{l=1}^{\alpha} (-1)^{l}\CONST{l}{1}_1\CONST{{\alpha-l}}{1}_1 = - \frac{(\ZETA-1)^\alpha}{\alpha!}(1+\BIGO{\QUARK}).}
\end{lemma}

\begin{proof}{Proof. }
If \m{n>2} then part~\eqref{gammaprop3} of lemma~\ref{LEMMAX3.02} implies that \algn{\textstyle
    \sum_{l=1}^{\alpha} (-1)^{l}\CONST{l}{1}_1\CONST{{\alpha-l}}{1}_1   & %
    \textstyle = \sum_{l=1}^{\alpha}  \frac{(-1)^{\alpha+l}(1-\ZETA ^p)^2}{l!(\alpha-l)!(1-\ZETA )^{2p-\alpha}} (1+\BIGO{\QUARK})
    \\ & \textstyle = \frac{(\ZETA-1)^\alpha}{\alpha!}\sum_{l=1}^{\alpha} (-1)^{l} \binom{\alpha}{l} \frac{(1-\ZETA^p)^2}{(1-\ZETA)^{2p}} (1+\BIGO{\QUARK})
        \\ & \textstyle = -\frac{(\ZETA-1)^\alpha}{\alpha!} (1+\BIGO{\QUARK}).
} For the third equality we use the equations \algn{\sum_{l=1}^{\alpha} (-1)^{l} \binom{\alpha}{l} = -1, \text{ and }  (1-\ZETA^p)(1-\ZETA)^{-p} = 1 + \BIGO{\QUARK}.}
Similarly, if \m{n=2} then part~\eqref{gammaprop3} of lemma~\ref{LEMMAX3.02} implies that
\algn{\textstyle
    \sum_{l=1}^{\alpha} (-1)^{l}\CONST{l}{1}_1\CONST{{\alpha-l}}{1}_1   & %
    \textstyle = \sum_{l=1}^{\alpha}  \frac{(-1)^{\alpha+l}p^2}{l!(\alpha-l)!(1-\ZETA )^{2p-\alpha-2}} (1+\BIGO{\QUARK})
    \\ & \textstyle = \frac{(\ZETA-1)^\alpha}{\alpha!}\sum_{l=1}^{\alpha} (-1)^{l} \binom{\alpha}{l} \frac{p^2}{(1-\ZETA)^{2p-2}} (1+\BIGO{\QUARK})
        \\ & \textstyle = -\frac{(\ZETA-1)^\alpha}{\alpha!} (1+\BIGO{\QUARK}).
} For the third equality we use the equation \algn{p^2(1-\ZETA)^{2-2p} = \CONST{0}{1}_1^2(1+\BIGO{\QUARK}) = 1+\BIGO{\QUARK},}
which follows from part~\eqref{gammaprop3} of lemma~\ref{LEMMAX3.02}.
\end{proof}

\section{%
Concrete elements of the ideal \texorpdfstring{\m{{\ZEROIDEAL}}}{I}} %

In this section we prove that the image of  \m{\HECKET-a} is contained in \m{\ZEROIDEAL}, which implies that integral elements in the image of  \m{\HECKET-a}  reduce modulo \m{\maxideal} to elements that belong to \m{\INAAA}. Finding such elements is relevant to computing \m{\overline{\Theta}_{k,a,\VAREPSILONp}} because equation~\eqref{eq34f3g3g3g45g445g45gg45} implies that \Algn{\label{eqf34g549gi409i} \overline{\Theta}_{k,a,\VAREPSILONp} \cong \ind_{\INZA}^G \SIGMA_{r}/\INAAA.}
  
\begin{lemma}\label{DECOMP1}
The group \m{G} can be decomposed as the disjoint union of double cosets
\m{\cup_{0\LEQSLANT j\LEQSLANT n} Bx_jI(n)}, where \m{x_i = \begin{Smatrix}
    1 & 0 \\ p^i & 1
    \end{Smatrix}} for \m{i\in\{0,\ldots,n-1\}} and \m{x_n=\begin{Smatrix}
    1 & 0 \\ 0 & 1
    \end{Smatrix}}.
\end{lemma}

\begin{proof}{Proof. }
  We can write \m{K} as the  (not disjoint) union of cosets  \m{\cup_{y\in S} \hspace{1pt}yI(n)} with
  \algn{S = \{\begin{Lmatrix}
     1 & 0 \\ l & 1
    \end{Lmatrix} \,|\, l\in\ZZ_p \} \cup \{\begin{Lmatrix}
     lp & -1 \\ 1 & 0
    \end{Lmatrix} \,|\, l\in\ZZ_p\}.}
   Indeed, if \m{\begin{Smatrix} \aaa & \bbb \\ \ccc & \ddd
  \end{Smatrix} \in K} then
  \algn{
 \begin{Lmatrix} \aaa & \bbb \\ \ccc & \ddd
  \end{Lmatrix}  = 
  \left\{\begin{array}{rl} \begin{Lmatrix} 1 & 0 \\ \ccc/\aaa & 1
  \end{Lmatrix} \begin{Lmatrix} \aaa & \bbb \\ 0 & (\aaa\ddd-\bbb\ccc)/\aaa
  \end{Lmatrix} & \text{if } \aaa \in \ZZ_p^\times\text{,}  \\[2pt] \begin{Lmatrix} \aaa/\ccc & -1 \\ 1 & 0
  \end{Lmatrix} \begin{Lmatrix} \ccc & \ddd \\ 0 & (\aaa\ddd-\bbb\ccc)/\ccc
  \end{Lmatrix}  & \text{if } \aaa \in p\ZZ_p \text{ (and so }\ccc\in\ZZ_p^\times \text{)}\text{.} \end{array}\right.
  }
 Together with the Iwasawa decomposition \m{G=BK}, this results in a decomposition of \m{G} as the (not disjoint) union   \m{G = \cup_{y\in S} ByI(n)}. Since
\algn{ 
\textstyle\begin{Lmatrix}
1 & 0 \\ up^i & 1
\end{Lmatrix} & \textstyle{}\vphantom{\begin{Lmatrix}
1 & -1/(1-tp) \\ 0 & 1/(1-tp)
\end{Lmatrix}x_0 \begin{Lmatrix}
1-tp&-1\\0&1
\end{Lmatrix} \in I(n)x_0B}\in \hspace{2pt} x_nI(n) \text{ for } u\in \ZZ_p^\times \text{ and }{i\in\ZZ_{\GEQSLANT n}},\\
\textstyle\begin{Lmatrix}
1 & 0 \\ up^i & 1
\end{Lmatrix} & \textstyle{}\vphantom{\begin{Lmatrix}
1 & -1/(1-tp) \\ 0 & 1/(1-tp)
\end{Lmatrix}x_0 \begin{Lmatrix}
1-tp&-1\\0&1
\end{Lmatrix} \in I(n)x_0B}= \begin{Lmatrix}
1 & 0 \\ 0 & u
\end{Lmatrix}x_i \begin{Lmatrix}
1&0\\0&1/u
\end{Lmatrix} \in Bx_iI(n) \text{ for } u\in \ZZ_p^\times \text{ and }{i\in\{0,\ldots,n-1\}},
\\ \textstyle\begin{Lmatrix}
lp & -1 \\ 1 & 0
\end{Lmatrix} & \vphantom{\begin{Lmatrix}
1 & -1 \\ 0 & 1
\end{Lmatrix}x_0 \begin{Lmatrix}
1&1+tp\\0&-1
\end{Lmatrix} \in Bx_0I(n)}= \begin{Lmatrix}
1 & lp-1 \\ 0 & 1
\end{Lmatrix}x_0 \begin{Lmatrix}
1&-1\\0&1
\end{Lmatrix} \in Bx_0I(n) \text{ for } l\in \ZZ_p,
} the double cosets of \m{x_0,\ldots,x_n} cover \m{G}. It is easy to verify that they are disjoint because if
\algn{ \begin{Lmatrix}
1 & 0 \\ q_1 & 1
\end{Lmatrix} = \begin{Lmatrix}
x & y \\ 0 & z
\end{Lmatrix}\begin{Lmatrix}
1 & 0 \\ q_2 & 1
\end{Lmatrix} \begin{Lmatrix}
\aaa&\bbb\\\ccc p^n&\ddd
\end{Lmatrix}    \text{ for }\begin{Lmatrix}
x & y \\ 0 & z
\end{Lmatrix} \in B \text{ and }\begin{Lmatrix}
\aaa&\bbb\\\ccc p^n&\ddd
\end{Lmatrix} \in I(n),} %
then we can solve this equation for \m{q_1} to get that \m{\textstyle q_1=\frac{\aaa q_2+ \ccc p^n}{\ddd+\bbb q_2}\in\ZZ_p},
 so if \m{q_1} and \m{q_2} are in \m{\{p^i \,|\, i\in\{0,\ldots,n-1\}\}\cup\{0\}} then \m{q_1=q_2}. %
\end{proof}

\begin{lemma}\label{DECOMP2}
Let \m{\delta} be %
 as in equation~\eqref{eq3goigj45fjiugifweko3ig}, and let \m{x_0,\ldots,x_n} be %
  as in lemma~\ref{DECOMP1}.
\begin{enumerate}[leftmargin=*]
    \item\label{decomp1.0} For \m{i\in\{0,\ldots,n\}},
\m{\sum_{\MU\in\FF_p}\delta\left(x_i\begin{Lmatrix}
    p & [\MU] \\ 0 & 1
    \end{Lmatrix}\right) = \VERIFYIF{i=0}ap^{2-k}}.\AFIOFJ
    \item\label{decomp1.1}
    For \m{g \in G},
  \m{\sum_{\MU\in\FF_p}\delta\left(g\begin{Lmatrix}
    p & [\MU] \\ 0 & 1
    \end{Lmatrix}\right) = ap^{2-k}\delta(g)}.\AFIOFJ
    \end{enumerate}
\end{lemma}

\begin{proof}{Proof. }
 \eqref{decomp1.0}   If \m{i=n} then \algn{\textstyle \sum_{\MU\in\FF_p}\delta\left(x_i\begin{Lmatrix}
    p & [\MU] \\ 0 & 1
    \end{Lmatrix}\right) = \sum_{\MU\in\FF_p}\delta\left(\begin{Lmatrix}
    p & [\MU] \\ 0 & 1
    \end{Lmatrix}\right) = a p^{2-k}.}
    If \m{i=n-1} then, for \m{\MU\in\FF_p}, \algn{ \textstyle 
    x_i\begin{Lmatrix}
    p & [\MU] \\ 0 & 1
    \end{Lmatrix} & \textstyle = \begin{Lmatrix}
    p & [\MU]/(1+[\MU]p^{n-1}) \\ 0 & 1
    \end{Lmatrix}\begin{Lmatrix}
    1/(1+[\MU]p^{n-1}) & 0 \\ p^{n} & 1+[\MU]p^{n-1}
    \end{Lmatrix} \in BI(n),}
 and consequently \algn{\textstyle\sum_{\MU\in\FF_p}\delta\left(x_i\begin{Lmatrix}
    p & [\MU] \\ 0 & 1
    \end{Lmatrix}\right) = a p^{1-k} \sum_{\MU\in\FF_p} \VAREPSILONp^{-1}(1+[\MU]p^{n-1}) = 0.} %
    If \m{i\in \{0,\ldots,n-2\}} then  the equation  \algn{\textstyle x_i\begin{Lmatrix}
    p & [\MU] \\ 0 & 1
    \end{Lmatrix} = \begin{Lmatrix}
    p & [\MU] \\ p^{i+1} & 1+[\MU]p^i
    \end{Lmatrix} = \begin{Lmatrix}
x & y \\ 0 & z
\end{Lmatrix} \begin{Lmatrix}
\aaa&\bbb\\\ccc p^n&\ddd
\end{Lmatrix}}
    has no solutions in \m{\MU\in\FF_p}, \m{\begin{Smatrix}
    x & y \\ 0 & z
    \end{Smatrix}\in B}, and \m{\begin{Smatrix}
    \aaa & \bbb \\ \ccc p^n & \ddd
    \end{Smatrix}\in I(n)} (a direct calculation shows that the existence of such a solution would imply that \m{p^{i+1}=\BIGO{p^n}}), so 
    \m{x_i\begin{Smatrix}
    p & [\MU] \\ 0 & 1
    \end{Smatrix}\not\in BI(n)} for all \m{\MU\in\FF_p} and consequently \algn{\textstyle \sum_{\MU\in\FF_p}\delta\left(x_i\begin{Lmatrix}
    p & [\MU] \\ 0 & 1
    \end{Lmatrix}\right) = 0.}%
    \eqref{decomp1.1} Let \m{g = bx_iu} for some \m{b \in B}, \m{i \in \{0,\ldots,n\}}, and \m{u = \begin{Smatrix} \aaa & \bbb \\ \ccc p^n & \ddd
  \end{Smatrix} \in I(n)}. Then \algn{\begin{Lmatrix}
    \aaa & \bbb \\ \ccc p^n & \ddd
    \end{Lmatrix}\begin{Lmatrix}
  p & [\MU] \\ 0 & 1
    \end{Lmatrix} = \begin{Lmatrix}
    p & [(\aaa \MU+\bbb)/\ddd] \\ 0 & 1
    \end{Lmatrix}\begin{Lmatrix}
   \aaa - \ccc [(\aaa \MU+\bbb)/\ddd] p^n & 
   (\aaa[\MU]+\bbb -\ddd[(\aaa \MU+\bbb)/\ddd])/p - \ccc [(\aaa \MU+\bbb)\MU/\ddd]p^{n-1}
   \\ \ccc p^{n+1} & \ddd+\ccc [\MU]p^n
    \end{Lmatrix},}
so, for \m{\MU\in\FF_p},
    \algn{\delta\left(g\begin{Lmatrix}
    p & [\MU] \\ 0 & 1
    \end{Lmatrix}\right) = \delta(b)\, \delta\left(x_iu\begin{Lmatrix}
    p & [\MU] \\ 0 & 1
    \end{Lmatrix}\right) = \delta(b)\, \delta\left(x_i\begin{Lmatrix}
    p & [(\aaa\MU+\bbb)/\ddd] \\ 0 & 1
    \end{Lmatrix}\right)\delta(u).}
    Therefore, by part~\eqref{decomp1.0},
    \algn{\sum_{\MU\in\FF_p}\delta\left(g\begin{Lmatrix}
    p & [\MU] \\ 0 & 1
    \end{Lmatrix}\right) & \textstyle\vphantom{\delta(bu) \sum_{\LAMBDA\in\FF_p} \delta\left(x_i\begin{Lmatrix}
    p & [\LAMBDA] \\ 0 & 1
    \end{Lmatrix}\right) } = \delta(bu) \sum_{\MU\in\FF_p} \delta\left(x_i\begin{Lmatrix}
    p & [(\aaa\MU+\bbb)/\ddd] \\ 0 & 1
    \end{Lmatrix}\right) \\ & \textstyle\vphantom{\delta(bu) \sum_{\LAMBDA\in\FF_p} \delta\left(x_i^{-1}\begin{Lmatrix}
    p & [\LAMBDA] \\ 0 & 1
    \end{Lmatrix}\right) } =  \delta(bu) \sum_{\LAMBDA\in\FF_p} \delta\left(x_i\begin{Lmatrix}
    p & [\LAMBDA] \\ 0 & 1
    \end{Lmatrix}\right)  \\ & \textstyle\vphantom{\delta(bu) \sum_{\LAMBDA\in\FF_p} \delta\left(x_i^{-1}\begin{Lmatrix}
    p & [\LAMBDA] \\ 0 & 1
    \end{Lmatrix}\right) } = \VERIFYIF{i=0}ap^{2-k}\delta(bu) = ap^{2-k}\delta(g).}
\end{proof}

\begin{lemma}\label{LEMMASETUP1}
{
The image \m{\Img(\HECKET-a)\subseteq \ind_{\INZA}^{G}\TILDE\SIGMA_{k-2}} is contained in \m{\ZEROIDEAL}.
}
\end{lemma}

\begin{proof}{Proof. }
  For all \m{\gamma\in G} %
   and \m{v\in\TILDE \SIGMA_r} we have
\Algn{\label{eqg98j5349e3ju3h4ytujrew}\nonumber
    & \textstyle %
   \vphantom{\tau(u_0) \COMPARR\left( {\sum_{\MU\in\FF_p} (\gamma u_0^{-1})\begin{Lmatrix}
    p &[\MU] \\ 0 & 1
    \end{Lmatrix} }\sqbrackets{}{\INZB}{\QQp}{\left(\begin{Lmatrix}
    1 & -[\MU] \\ 0 & p
    \end{Lmatrix} \cdot (u_0\cdot v)\right)}\right)} \COMPARR\left(\HECKET({\gamma}\sqbrackets{}{\INZB}{\QQp}{v})\right) 
        \\&\textstyle\nonumber\null\qquad\vphantom{\tau(u_0) \COMPARR\left( {\sum_{\MU\in\FF_p} (\gamma u_0^{-1})\begin{Lmatrix}
    p &[\MU] \\ 0 & 1
    \end{Lmatrix} }\sqbrackets{}{\INZB}{\QQp}{\left(\begin{Lmatrix}
    1 & -[\MU] \\ 0 & p
    \end{Lmatrix} \cdot (u_0\cdot v)\right)}\right)} =  \COMPARR\left( {\sum_{\MU\in\FF_p} \gamma \begin{Lmatrix}
    p &[\MU] \\ 0 & 1
    \end{Lmatrix} }\sqbrackets{}{\INZB}{\QQp}{\left(\begin{Lmatrix}
    1 & -[\MU] \\ 0 & p
    \end{Lmatrix} \cdot v\right)}\right)
       \\&\textstyle\nonumber\null\qquad\vphantom{\tau(u_0) \COMPARR\left( {\sum_{\MU\in\FF_p} (\gamma u_0^{-1})\begin{Lmatrix}
    p &[\MU] \\ 0 & 1
    \end{Lmatrix} }\sqbrackets{}{\INZB}{\QQp}{\left(\begin{Lmatrix}
    1 & -[\MU] \\ 0 & p
    \end{Lmatrix} \cdot (u_0\cdot v)\right)}\right)} = \left(\left(\gamma \begin{Lmatrix}
    p & [\MU] \\ 0 & 1
    \end{Lmatrix}\begin{Lmatrix}
    1 & -[\MU] \\ 0 & p
    \end{Lmatrix}\right)\cdot v\right)\otimes \left(g\mapsto{}\sum_{\MU\in\FF_p}\delta\left(g\gamma \begin{Lmatrix}
    p & [\MU] \\ 0 & 1
    \end{Lmatrix}\right)L(1)\right)
         \\&\textstyle\null\qquad\vphantom{\tau(u_0) \COMPARR\left( {\sum_{\MU\in\FF_p} (\gamma u_0^{-1})\begin{Lmatrix}
    p &[\MU] \\ 0 & 1
    \end{Lmatrix} }\sqbrackets{}{\INZB}{\QQp}{\left(\begin{Lmatrix}
    1 & -[\MU] \\ 0 & p
    \end{Lmatrix} \cdot (u_0\cdot v)\right)}\right)} = a \left(\gamma\cdot v\right)\otimes \left(g\mapsto{}\delta(g\gamma)L(1)\right) = \COMPARR\left(a{\gamma}\sqbrackets{}{\INZB}{\QQp}{v}\right) .
} 
 The first equality follows from equation~\eqref{eqHecke}, and the second equality follows from equation~\eqref{eqf34ufh348gh348hf384h}.
The third equality follows from part~\eqref{decomp1.1} of lemma~\ref{DECOMP2} and the equation
\algn{\textstyle \left(\gamma \begin{Lmatrix}
    p & [\MU] \\ 0 & 1
    \end{Lmatrix}\begin{Lmatrix}
    1 & -[\MU] \\ 0 & p
    \end{Lmatrix}\right) \cdot v = \left(\gamma \begin{Lmatrix}
    p & 0 \\ 0 & p
    \end{Lmatrix}\right) \cdot v = p^{k-2} \left(\gamma \cdot v\right).}
Equation~\eqref{eqg98j5349e3ju3h4ytujrew} implies that
  \m{\Img(\HECKET-a)\subseteq\Ker(\COMPARR)}.
\end{proof}

\section{%
Subquotients of \texorpdfstring{\m{{\ind_{\INZA}^{G}\SIGMA_r}}}{ind(\textSigma)}} \label{secsubs}

In this section we   analyze the factors of \m{\ind_{\INZA}^{KZ}\SIGMA_r}. The goal is to use equation~\eqref{eqf34g549gi409i} to deduce information about \m{\overline{\Theta}_{k,a,\VAREPSILONp}} by finding elements in the factors of \m{\ind_{\INZA}^{KZ}\SIGMA_r} that are represented by elements of \m{\INAAA}.

\begin{lemma}\label{LEMMAX1.0}
{
Let \m{h\in\ZZ}, and let \m{v} be a generator of \m{\chi_h}.
 There is an isomorphism \algn{
 \Lambda_h : I_{h} \xrightarrow{\hspace{\SIZEEE} {\cong} \hspace{\SIZEEE}} \ind_{I(1)Z}^{KZ} \chi_h}
 given by the map that sends
\m{f \in I_{h}} to
\algn{
        \textstyle   & \textstyle  \sum_{\TEMPORARYA \in \FF_p} 
    f(-\TEMPORARYA,1)\begin{Lmatrix}
    1 & 0 \\ [\TEMPORARYA]  & 1
    \end{Lmatrix}  \sqbrackets{}{KZ,I(1)Z}{\FFp}{v} %
    + \textstyle 
    f(-1,0)\begin{Lmatrix}
    0 & -1 \\ 1 & 0
    \end{Lmatrix}  \sqbrackets{}{KZ,I(1)Z}{\FFp}{v}.
    }
}
\end{lemma}
 
 \begin{proof}{Proof. }
   We can decompose \m{KZ} as the  union of cosets  \m{\cup_{y\in S} \hspace{1pt}yI(1)Z} with
  \algn{S = \{\begin{Lmatrix}
     1 & 0 \\ [\TEMPORARYA] & 1
    \end{Lmatrix} \,|\, \TEMPORARYA\in\FF_p \} \cup \{\begin{Lmatrix}
     0 & -1 \\ 1 & 0
    \end{Lmatrix} \}.}
   Indeed, if \m{\begin{Smatrix} \aaa & \bbb \\ \ccc & \ddd
  \end{Smatrix} \in K} then
  \algn{
 \begin{Lmatrix} \aaa & \bbb \\ \ccc & \ddd
  \end{Lmatrix}  = 
  \left\{\begin{array}{rl} \begin{Lmatrix} 1 & 0 \\ \ccc/\aaa -sp & 1
  \end{Lmatrix} \begin{Lmatrix} \aaa & \bbb \\ \aaa sp & \bbb sp+(\aaa\ddd-\bbb\ccc)/\aaa
  \end{Lmatrix} & \text{if } \aaa \in \ZZ_p^\times\text{,}  \\[2pt] \begin{Lmatrix} \aaa/\ccc+tp & -1 \\ 1 & 0
  \end{Lmatrix} \begin{Lmatrix} \ccc & \ddd \\ \ccc tp & \ddd tp+(\aaa\ddd-\bbb\ccc)/\ccc
  \end{Lmatrix}  & \text{if } \aaa \in p\ZZ_p \text{ (and so }\ccc\in\ZZ_p^\times \text{)}\text{,} \end{array}\right.
  }
  and after choosing suitable \m{s,t \in \ZZ_p} we can conclude that \m{\begin{Smatrix} \aaa & \bbb \\ \ccc & \ddd
  \end{Smatrix} \in yI(1)} for some \m{y\in S}. It is easy to show by a direct calculation that this is a disjoint union of cosets. Therefore \m{\Lambda_h} %
   is bijective. It is also evidently linear. For \m{\TEMPORARYA\in\FF_p}, let
  \algn{ \DOUBBR{-\TEMPORARYA,1} = \begin{Lmatrix}
    1 & 0 \\ [\TEMPORARYA]  & 1
    \end{Lmatrix} \sqbrackets{}{KZ,I(1)Z}{\FFp}{v}, \text{ and } \DOUBBR{-1,0} = \begin{Lmatrix} 0 & -1 \\ 1 & 0 \end{Lmatrix}  \sqbrackets{}{KZ,I(1)Z}{\FFp}{v}.} 
    If \m{\begin{Smatrix} \aaa & \bbb \\ \ccc & \ddd
  \end{Smatrix} \in K} then
  \algn{\begin{Lmatrix} \aaa & \bbb \\ \ccc & \ddd
  \end{Lmatrix} \begin{Lmatrix}
    1 & 0 \\ [\TEMPORARYA]  & 1
    \end{Lmatrix}  & \textstyle= \left\{\begin{array}{rl}   \vphantom{\begin{Lmatrix} 1 & 0 \\ (\ccc+\ddd[\TEMPORARYA])/(\aaa+\bbb[\TEMPORARYA]) & 1
  \end{Lmatrix} \begin{Lmatrix} \aaa+\bbb[\TEMPORARYA] & \bbb \\ 0 & (\aaa\ddd-\bbb\ccc)/(\aaa+\bbb[\TEMPORARYA])
  \end{Lmatrix}}\begin{Lmatrix} 1 & 0 \\ (\ccc+\ddd[\TEMPORARYA])/(\aaa+\bbb[\TEMPORARYA]) & 1
  \end{Lmatrix} \begin{Lmatrix} \aaa+\bbb[\TEMPORARYA] & \bbb \\ 0 & (\aaa\ddd-\bbb\ccc)/(\aaa+\bbb[\TEMPORARYA])
  \end{Lmatrix} & \text{if } \aaa+\bbb[\TEMPORARYA]\in\ZZ_p^\times\text{,}  \\[2pt]  \vphantom{\begin{Lmatrix} 1 & 0 \\ (\ccc+\ddd[\TEMPORARYA])/(\aaa+\bbb[\TEMPORARYA]) & 1
  \end{Lmatrix} \begin{Lmatrix} \aaa+\bbb[\TEMPORARYA] & \bbb \\ 0 & (\aaa\ddd-\bbb\ccc)/(\aaa+\bbb[\TEMPORARYA])
  \end{Lmatrix}}\begin{Lmatrix} 0 & -1 \\ 1 & 0 \end{Lmatrix} \begin{Lmatrix} \ccc+\ddd[\TEMPORARYA] & \ddd \\ -\aaa-\bbb[\TEMPORARYA] & -\bbb
  \end{Lmatrix}  & \text{if } \aaa+\bbb[\TEMPORARYA]\in p\ZZ_p\text{,} \end{array}\right.
  \\[3pt]
  \begin{Lmatrix} \aaa & \bbb \\ \ccc & \ddd
  \end{Lmatrix} \begin{Lmatrix} 0 & -1 \\ 1 & 0 \end{Lmatrix} 
  & \textstyle= 
  \left\{\begin{array}{rl}  
  \vphantom{\begin{Lmatrix} 1 & 0 \\ (\ccc+\ddd[\TEMPORARYA])/(\aaa+\bbb[\TEMPORARYA]) & 1
  \end{Lmatrix} \begin{Lmatrix} \aaa+\bbb[\TEMPORARYA] & \bbb \\ 0 & (\aaa\ddd-\bbb\ccc)/(\aaa+\bbb[\TEMPORARYA])
  \end{Lmatrix}}\begin{Lmatrix} 1 & 0 \\ \ddd/\bbb & 1
  \end{Lmatrix} \begin{Lmatrix} \bbb & -\aaa \\ 0 & (\aaa\ddd-\bbb\ccc)/\bbb
  \end{Lmatrix} & \text{if } \bbb\in\ZZ_p^\times\text{,} 
  \\[2pt]  \vphantom{\begin{Lmatrix} 1 & 0 \\ (\ccc+\ddd[\TEMPORARYA])/(\aaa+\bbb[\TEMPORARYA]) & 1
  \end{Lmatrix} \begin{Lmatrix} \aaa+\bbb[\TEMPORARYA] & b \\ 0 & (\aaa\ddd-\bbb\ccc)/(\aaa+\bbb[\TEMPORARYA])
  \end{Lmatrix}} \begin{Lmatrix} 0 & -1 \\ 1 & 0 \end{Lmatrix} \begin{Lmatrix} \ddd & -\ccc \\ -\bbb & \aaa
  \end{Lmatrix}  & \text{if } \bbb\in p\ZZ_p\text{,} \end{array}\right.}
  and consequently
  \Algn{\label{eqfi3j498gh5y8turigturh64j4oe}\nonumber\begin{Lmatrix} \aaa & \bbb \\ \ccc & \ddd
  \end{Lmatrix}\DOUBBR{-\TEMPORARYA,1}  & \textstyle= \left\{\begin{array}{rl}  ((\aaa\ddd-\bbb\ccc)/(\aaa+\bbb\TEMPORARYA))^h \DOUBBR{-(\ccc+\ddd\TEMPORARYA)/(\aaa+\bbb\TEMPORARYA),1} & \text{if } \aaa+\bbb[\TEMPORARYA]\in\ZZ_p^\times\text{,}  \\[2pt] (-\bbb)^h\DOUBBR{-1,0}   & \text{if } \aaa+\bbb[\TEMPORARYA]\in p\ZZ_p\text{,} \end{array}\right.
  \\[3pt]
  \begin{Lmatrix} \aaa & \bbb \\ \ccc & \ddd
  \end{Lmatrix} \DOUBBR{-1,0}
  & \textstyle= 
  \left\{\begin{array}{rl}  
((\aaa\ddd-\bbb\ccc)/\bbb)^h\DOUBBR{-\ddd/\bbb,1} 
  & \text{if } \bbb\in\ZZ_p^\times\text{,}
  \\[2pt]  \aaa^h\DOUBBR{-1,0} & \text{if } \bbb\in p\ZZ_p\text{.} \end{array}\right.}
 For \m{t = (\LAMBDA,\MU) \in \FF_p^2\backslash\{(0,0)\}}, let us define
 \algn{\DOUBBR{t} = \left\{\begin{array}{rl}\MU^{-h}\DOUBBR{\LAMBDA/\MU,1} & \text{if } \MU \in \FF_p^\times\text{,}  \\[2pt] (-\LAMBDA)^{-h} \DOUBBR{-1,0}  & \text{if } \MU=0\text{.} \end{array}\right.}
Then we can summarize equation~\eqref{eqfi3j498gh5y8turigturh64j4oe} as
 \m{\begin{Smatrix} \aaa & \bbb \\ \ccc & \ddd
  \end{Smatrix} \DOUBBR{t} =  \DOUBBR{t\begin{Smatrix} \aaa & \bbb \\ \ccc & \ddd
  \end{Smatrix}^{-1}}} for all \m{\begin{Smatrix} \aaa & \bbb \\ \ccc & \ddd
  \end{Smatrix} \in K} and  \m{t \in \FF_p^2\backslash\{(0,0)\}}.  In particular, for all \m{\LAMBDA\in\FF_p} and  \m{t \in \FF_p^2\backslash\{(0,0)\}}, \m{\DOUBBR{\LAMBDA t} = \LAMBDA^{-h}\DOUBBR{t}}, so \algn{s \in \bb P^1(\FF_p) \xmapsto{ \ \ } f(s)\DOUBBR{s}\in \ind_{I(1)Z}^{KZ}\chi_h} is a well-defined map
 and
\algn{\textstyle \begin{Lmatrix} \aaa & \bbb \\ \ccc & \ddd
  \end{Lmatrix} \Lambda_h(f) &\textstyle \vphantom{\begin{Lmatrix} \aaa & \bbb \\ \ccc & \ddd
  \end{Lmatrix}\left(\sum_{s\in\bb P^1(\FF_p)} 
    f(s)\DOUBBR{s}\right)}= \begin{Lmatrix} \aaa & \bbb \\ \ccc & \ddd
  \end{Lmatrix}\left(\sum_{\TEMPORARYA \in \FF_p} 
    f(-\TEMPORARYA,1)\DOUBBR{-\TEMPORARYA,1} %
    + \textstyle 
    f(-1,0)\DOUBBR{-1,0}\right)
    \\ &\textstyle\vphantom{\begin{Lmatrix} \aaa & \bbb \\ \ccc & \ddd
  \end{Lmatrix}\left(\sum_{s\in\bb P^1(\FF_p)} 
    f(s)\DOUBBR{s}\right)} = 
 \textstyle\begin{Lmatrix} \aaa & \bbb \\ \ccc & \ddd
  \end{Lmatrix}\left(\sum_{s\in\bb P^1(\FF_p)} 
    f(s)\DOUBBR{s}\right)
        \\ &\textstyle\vphantom{\begin{Lmatrix} \aaa & \bbb \\ \ccc & \ddd
  \end{Lmatrix}\left(\sum_{s\in\bb P^1(\FF_p)} 
    f(s)\DOUBBR{s}\right)} = 
 \textstyle\sum_{s\in\bb P^1(\FF_p)} 
    f(s)\DOUBBR{s\begin{Lmatrix} \aaa & \bbb \\ \ccc & \ddd
  \end{Lmatrix}^{-1}}
          \\ &\textstyle\vphantom{\begin{Lmatrix} \aaa & \bbb \\ \ccc & \ddd
  \end{Lmatrix}\left(\sum_{s\in\bb P^1(\FF_p)} 
    f(s)\DOUBBR{s}\right)} = 
 \textstyle\sum_{w\in\bb P^1(\FF_p)} 
    f(w\begin{Lmatrix} \aaa & \bbb \\ \ccc & \ddd
  \end{Lmatrix})\DOUBBR{w}
 = \Lambda_h(\begin{Lmatrix} \aaa & \bbb \\ \ccc & \ddd
  \end{Lmatrix} f),
    } implying that \m{\Lambda_h} is a homomorphism.
 \end{proof}

\begin{lemma}\label{LEMMAX1.02}
{
Let \m{h\in\ZZ} and \m{m\in \ZZ_{\GEQSLANT1}}, and let \m{v} be a generator of \m{\chi_h}. %
 Let
\algn{
        \textstyle   & \textstyle e_i =   \sum_{ \TEMPORARYA \in \FF_p} \TEMPORARYA^{i}
    \begin{Lmatrix}
    1 & 0 \\ [\TEMPORARYA]p^m & 1
    \end{Lmatrix}  \sqbrackets{}{I(m),I(m+1)}{\FFp}{v}\in\ind_{I(m+1)Z}^{I(m)Z}\chi_{h} \text{ for } { i \in \{0,\ldots, p-1\}}.  %
    }
 There is a series
    \algn{\ind_{I(m+1)Z}^{I(m)Z} \chi_{h}=N_{p}\supset N_{p-1}\supset \cdots \supset N_{0} = \{0\}}
with the following properties. \begin{enumerate}[leftmargin=*]
    \item\label{dec0.01} For \m{i\in\{0,\ldots,p-1\}}, \m{N_{i+1}} is the linear span of all \m{e_{j}} for all \m{j \in \{0,\ldots,i\}}. %
    \AFIOFJ
    \item\label{dec0.02} %
    For \m{i\in\{0,\ldots,p-1\}},  there is an isomorphism \algn{{\chi_{h-2i}(i)} \xrightarrow{\hspace{\SIZEEE} {\cong} \hspace{\SIZEEE}} N_{i+1}/N_{i}}
 given by the map that sends a generator \m{v_i} of \m{\chi_{h-2i}(i)} to  
    \m{e_{i} + N_{i}\in N_{i+1}/N_{i}}.\AFIOFJ
\end{enumerate}
}
\end{lemma}

\begin{proof}{Proof. }
 It is easy to show just like in the proof of lemma~\ref{LEMMAX1.0} that we can decompose \m{I(m)Z} as the  disjoint union of cosets  \m{\cup_{y\in S} \hspace{1pt}yI(m+1)Z} with
  \algn{S = \{\begin{Lmatrix}
     1 & 0 \\ [\TEMPORARYA]p^m & 1
    \end{Lmatrix} \,|\, \TEMPORARYA\in\FF_p \}.} For \m{\TEMPORARYA\in\FF_p}, let
  \algn{ \DOUBBR{\TEMPORARYA} = \begin{Lmatrix}
    1 & 0 \\ [\TEMPORARYA]p^m  & 1
    \end{Lmatrix} \sqbrackets{}{I(m),I(m+1)Z}{\FFp}{v}.} Then \m{\ind_{I(m+1)Z}^{I(m)Z} \chi_{h}} is  spanned by \m{\{\DOUBBR{\TEMPORARYA}\,|\,\TEMPORARYA\in\FF_p\}}. 
    For \m{\begin{Smatrix} \aaa & \bbb \\ \ccc p^m & \ddd
  \end{Smatrix} \in I(m)} and \m{\TEMPORARYA\in\FF_p},
  \algn{\begin{Lmatrix} \aaa & \bbb \\ \ccc p^m & \ddd
  \end{Lmatrix}\begin{Lmatrix} 1 & 0 \\ [\TEMPORARYA]p^m & 1
  \end{Lmatrix} = \begin{Lmatrix} 1 & 0 \\ (\ccc+\ddd [\TEMPORARYA])p^m/(\aaa+\bbb[\TEMPORARYA] p^m) & 1
  \end{Lmatrix}\begin{Lmatrix} \aaa+\bbb [\TEMPORARYA] p^m & \bbb \\ 0 & (\aaa\ddd-\bbb\ccc p^m)/(\aaa+\bbb[\TEMPORARYA] p^m)
  \end{Lmatrix},}
  so
  \algn{\begin{Lmatrix} \aaa & \bbb \\ \ccc p^m & \ddd
  \end{Lmatrix} \DOUBBR{\TEMPORARYA} = \ddd^h \DOUBBR{(\ccc+\ddd\TEMPORARYA)/\aaa}.} Therefore, for \m{i\in \{0,\ldots,p-1\}} and \m{j\in \{0,\ldots,i\}},
  \Algn{\label{eq3oigj4oigjg5}\nonumber\begin{Lmatrix} \aaa & \bbb \\ \ccc p^m & \ddd
  \end{Lmatrix} e_j & \textstyle= \begin{Lmatrix} \aaa & \bbb \\ \ccc p^m & \ddd
  \end{Lmatrix}\sum_{\TEMPORARYA\in\FF_p} \TEMPORARYA^j \DOUBBR{\TEMPORARYA} \\\nonumber&\textstyle{}= \ddd^h  \sum_{\TEMPORARYA\in\FF_p}  \TEMPORARYA^j \DOUBBR{(\ccc+\ddd\TEMPORARYA)/\aaa}
  \\\nonumber&\textstyle{}= \ddd^h  \sum_{\TEMPORARYAA\in\FF_p}  ((\aaa\TEMPORARYAA-\ccc)/\ddd)^j \DOUBBR{\TEMPORARYAA}
  \\&\textstyle{}=  \ddd^{h-j} \sum_{l=0}^j \binom{j}{l}\aaa^l (-\ccc)^{j-l} e_l,} implying that \m{N_{i+1}} is closed under the action of \m{I(m)Z}. Moreover, for \m{\LAMBDA\in \FF_p}, if \m{f_\LAMBDA\in \FFp[z]} is the polynomial \m{\prod_{\MU\in\FF_p\backslash\{\LAMBDA\}} (z-\MU)} of degree \m{p-1}, then
  \algn{\DOUBBR{\LAMBDA} = -\sum_{\TEMPORARYA\in\FF_p} f_\LAMBDA(\TEMPORARYA)\DOUBBR{\TEMPORARYA} \in N_p,} so \m{N_p = \ind_{I(m+1)Z}^{I(m)Z} \chi_{h}}. By definition, for \m{i\in\{0,\ldots,p-1\}}, \m{ N_{i+1}/N_{i}} is generated by \m{{e_{i} + N_{i}}} so it is at most one-dimensional. Since \m{N_p} has dimension \m{p}, it follows that each \m{N_{i+1}/N_i} must be exactly one-dimensional. Finally, equation~\eqref{eq3oigj4oigjg5} implies that the actions of \m{I(m)Z} on \m{\chi_{h-2i}(i) } and \m{N_{i+1}/N_i} match, so the two representations are isomorphic.
\end{proof}

\begin{lemma}\label{COROLEMMAX1}
{
Let \m{h\in\ZZ}, and let \m{v} be a generator of \m{\chi_h}. For \m{ i \in \{0,\ldots, p^{n-1}-1\}} and \m{f \in I_{h-2i}(i)}, let
\m{e_{i,f}\in \ind_{\INZA}^{KZ}\chi_{h}}
denote the element
\algn{
        \textstyle   & \textstyle  \sum_{\TEMPORARYA, \TEMPORARYA_1,\ldots,\TEMPORARYA_{n-1} \in \FF_p} \TEMPORARYA_1^{i_1}\cdots \TEMPORARYA_{n-1}^{i_{n-1}}
    f(-\TEMPORARYA,1)\begin{Lmatrix}
    1 & 0 \\ [\TEMPORARYA] + [\TEMPORARYA_1]p+\cdots+[\TEMPORARYA_{n-1}]p^{n-1} & 1
    \end{Lmatrix}  \sqbrackets{}{KZ,I(n)Z}{\FFp}{v} \\
    & \null\qquad + \textstyle \sum_{\TEMPORARYA_1,\ldots,\TEMPORARYA_{n-1} \in \FF_p} \TEMPORARYA_1^{i_1}\cdots \TEMPORARYA_{n-1}^{i_{n-1}}
    f(-1,0)\begin{Lmatrix}
     -[\TEMPORARYA_1]p-\cdots-[\TEMPORARYA_{n-1}]p^{n-1} & -1 \\ 1 & 0
    \end{Lmatrix}  \sqbrackets{}{KZ,I(n)Z}{\FFp}{v},
    }
where \m{i_{1},\ldots,i_{n-1} \in \{0,\ldots,p-1\}} %
are such that \m{i=i_1+\cdots+i_{n-1}p^{n-2}}. There is a filtration
    \algn{\ind_{\INZA}^{KZ} \chi_{h}=M_{p^{n-1}}\supset M_{p^{n-1}-1}\supset \cdots \supset M_{0} = \{0\}}
with the following properties. \begin{enumerate}[leftmargin=*]
    \item\label{dec1} For \m{i\in\{0,\ldots,p^n-1\}}, \m{M_{i+1}} is the linear span of all \m{e_{j,f}} for all \m{j \in \{0,\ldots,i\}} and all \m{f \in I_{h-2j}(j)}. %
    \AFIOFJ
    \item\label{dec2} %
    For \m{i\in\{0,\ldots,p^n-1\}}, there is an isomorphism \algn{I_{h-2i}(i) \xrightarrow{\hspace{\SIZEEE} {\cong} \hspace{\SIZEEE}} M_{i+1}/M_{i}}
 given by the map that sends 
\m{f \in I_{h-2i}(i)} to
    \m{e_{i,f} + M_{i}\in M_{i+1}/M_{i}}.\AFIOFJ
\end{enumerate}
}
\end{lemma}

\begin{proof}{Proof. }
Lemma~\ref{LEMMAX1.02} applied \m{n-1} times implies that there is a series
    \algn{\ind_{I(2)Z}^{I(1)Z}\cdots\ind_{\INZA}^{I(n-1)Z} \chi_{h} & \textstyle \cong \ind_{\INZA}^{I(1)Z} \chi_{h}\\ &  \textstyle=N_{p^{n-1}}^\prime\supset N_{p^{n-1}-1}^\prime\supset \cdots \supset N_{0}^\prime = \{0\}}
such that, for \m{i\in\{0,\ldots,p^{n-1}-1\}}, \algn{N_{i+1}^\prime/N_i^\prime\cong \chi_{h-2(i_1+\cdots+i_{n-1})}(i_1+\cdots+i_{n-1}),} where \m{i_{1},\ldots,i_{n-1} \in \{0,\ldots,p-1\}} %
are such that \m{i=i_1+\cdots+i_{n-1}p^{n-2}}. Roughly speaking, this is because (counting the factors from the bottom and starting from zero) the \m{i}th factor \m{N_{i+1}^\prime/N_i^\prime} of \m{\ind_{\INZA}^{I(1)Z} \chi_{h}} is the \m{i_1}th factor of the induction to \m{I(1)Z} of the \m{i_2}th factor of the induction to \m{I(2)Z} of \m{\cdots} of the \m{i_{n-1}}th factor of the induction to \m{I(n-1)Z} of \m{\chi_h}.  Since \algn{i \equiv \sum_{j=1}^{n-1} i_jp^{j-1}\equiv \sum_{j=1}^{n-1} i_j \bmod {p-1},} we also have \m{N_{i+1}^\prime/N_i^\prime\cong \chi_{h-2i}(i)}. If, for \m{ i \in \{0,\ldots, p^{n-1}-1\}}, 
\algn{
        \textstyle   & \textstyle e_i =  \sum_{ \TEMPORARYA_1,\ldots,\TEMPORARYA_{n-1} \in \FF_p} \TEMPORARYA_1^{i_1}\cdots \TEMPORARYA_{n-1}^{i_{n-1}}
    \begin{Lmatrix}
    1 & 0 \\ [\TEMPORARYA_1]p+\cdots+[\TEMPORARYA_{n-1}]p^{n-1} & 1
    \end{Lmatrix}  \sqbrackets{}{I(1),I(n)}{\FFp}{v} \in \ind_{\INZA}^{I(1)Z}\chi_{h}, %
    }
where \m{i_{1},\ldots,i_{n-1} \in \{0,\ldots,p-1\}} %
are such that \m{i=i_1+\cdots+i_{n-1}p^{n-2}}, then the description of the factors in lemma~\ref{LEMMAX1.02} implies that \m{N_{i+1}^\prime} is the linear span of all \m{e_{j}} for all \m{j \in \{0,\ldots,i\}}, and \m{N_{i+1}^\prime/N_{i}^\prime} is generated by
    \m{e_{i} + N_{i}^\prime\in N_{i+1}^\prime/N_{i}^\prime}. Lemma~\ref{LEMMAX1.0} then completes the proof.
\end{proof}

\begin{lemma}\label{LEMMAX2}
{
Recall that \m{\ast \mapsto {\ast^{\NEBENPOWEX-r}}} is the reduction modulo \m{p} of \m{\VAREPSILONp}. There is a series
 \Algn{%
 \label{eq3i0ijg3iogfi43gj4}\textstyle \SIGMA_r = \Span\{y^{r},xy^{r-1},\ldots\} \supset \Span\{xy^{r-1},\ldots\} \supset \cdots \supset \Span\{x^r\} \supset  \{0\}.%
 }
 The factors of this series are
 the one-dimensional modules \algn{\textstyle\chi_{2r-\NEBENPOWEY}(\NEBENPOWEX-r),\chi_{2r-\NEBENPOWEY-2}(\NEBENPOWEX-r+1),\ldots,\chi_{2-\NEBENPOWEY}(\NEBENPOWEX-1),\chi_{-\NEBENPOWEY}(\NEBENPOWEX).}
}
\end{lemma}

\begin{proof}{Proof. } It is clear that the factors of the series in equation~\eqref{eq3i0ijg3iogfi43gj4} are one-dimensional. The description of the  factors follows from the fact that \m{\begin{Smatrix}
\aaa&\bbb\\\ccc p^n&\ddd
\end{Smatrix}\in I(n)} %
 maps \algn{x^jy^{r-j}\xmapsto{ \ \ } \VAREPSILONp(\aaa)\aaa^j\ddd^{r-j} x^jy^{r-j}+\sum_{l=j+1}^{r} b_{\aaa,\bbb,\ddd,j,l} x^{l}y^{r-l} \text{ for some  } b_{\aaa,\bbb,\ddd,j,l},}
 and \m{\VAREPSILONp(u)u^jv^{r-j} = u^{\NEBENPOWEX-r+j}v^{r-j} = \det (\begin{Smatrix}
\aaa&\bbb\\\ccc p^n&\ddd
\end{Smatrix})^{\NEBENPOWEX-r+j} v^{2r-\NEBENPOWEX-2j}}.
\end{proof}

\begin{lemma}\label{LEMMAX3}
{ Let \m{U} be the submodule  \algn{\langle \{x^r\} \cup \{(x y^p-x^p y)x^{r-p-1-i}y^{i} \,|\, 0 \LEQSLANT i \LEQSLANT r-p-1\}\rangle \subseteq \SIGMA_r,}
 and let
 \m{\QUOTIENT =\SIGMA_r/U}. 
If %
 \m{h \in U} then \m{\textstyle{1}\sqbrackets{}{\INZB}{\FFp}{h} \in \INAAA  }. Hence there are the inclusions \algn{\ind_{\INZA}^G U \subseteq \INAAA   \subseteq \ind_{\INZA}^G \SIGMA_r,} and there is a surjective homomorphism 
\algn{\textstyle \QUOTIENTMAP : \ind_{\INZA}^G \QUOTIENT \xtwoheadrightarrow{ \ \ } \overline\Theta_{k,a,\varepsilon}.}}
\end{lemma}
\begin{proof}{Proof. }
 The proof is nearly identical to the proofs of lemmas~4.1 and~4.3 in~\cite{BG09}: we use the explicit formula for \m{\HECKET}  to show that
\Algn{\label{eq938jg93jg943jg4i9ug}
    \nonumber & \nonumber\textstyle (\HECKET-a)\left(-a^{-1}\sqbrackets{}{\INZB}{\QQp}{(x y^p-x^p y)x^{r-p-1-i}y^{i}}\right)\nonumber \\ & \textstyle \qquad = {1}\sqbrackets{}{\INZB}{\QQp}{(x y^p-x^p y)x^{r-p-1-i}y^{i}} + \BIGO{pa^{-1}}
}
for all \m{i\in\{0,\ldots,r-p+1\}}, and 
\Algn{\label{eq3ijg39gj93gj39}
     & \textstyle \textstyle (\HECKET-a)\left({-\begin{Lmatrix}
    1 & 0 \\ 0 & p
    \end{Lmatrix}}\sqbrackets{}{\INZB}{\QQp}{(x y^p-x^p y)x^{r-p}y^{-1}}\right)={1}\sqbrackets{}{\INZB}{\QQp}{x^{r}} + \BIGO{a}.
}
 The reductions modulo~\m{\maxideal} of the left sides of equations~\eqref{eq938jg93jg943jg4i9ug} and~\eqref{eq3ijg39gj93gj39} are in \m{ \INAAA }. The assumption \m{0 < \TADICVALUATION(a) < \frac{p-1}{2}} implies  \m{0<\PADICVALUATION(a)<1}, so the reductions
 modulo~\m{\maxideal} of the right sides of equations~\eqref{eq938jg93jg943jg4i9ug} and~\eqref{eq3ijg39gj93gj39} are
 \m{{1}\sqbrackets{}{\INZB}{\FFp}{(x  y^p-x^py)x^{r-p-1-i}y^{i}}} and 
 \m{{1}\sqbrackets{}{\INZB}{\FFp}{x^{r}}}, respectively.
\end{proof}

\begin{lemma}\label{LEMMAX3.043}
{Let \m{\QUOTIENT} and \m{U} be as in lemma~\ref{LEMMAX3}. Then \m{\QUOTIENT}  has a series
\algn{\QUOTIENT = \Span\{\GENN_1+U,\ldots,\GENN_{p-1}+U\} \supset \Span\{\GENN_1+U,\ldots,\GENN_{p-2}+U\}\supset\cdots%
\supset U} whose factors are the one-dimensional modules %
  \algn{\textstyle\chi_{-\NEBENPOWEY}(\NEBENPOWEX),\ldots,\chi_{2-\NEBENPOWEY}(\NEBENPOWEX-1),} so that, for \m{i\in\{1,\ldots,p-1\}}, \m{\overline{\GEN}_i} represents a generator of \m{\chi_{2i-\NEBENPOWEY}(\NEBENPOWEX-i)}.}  Therefore there is a filtration of \m{\ind_{I(n)Z}^{KZ} \QUOTIENT} %
  whose subquotients are
 \algn{\left\{I_{i,j}\cong I_{2r-\NEBENPOWEX-2(i+j)}(\NEBENPOWEX-r+i+j)\,|\,i\in\{0,\ldots,p^{n-1}-1\},\, j\in\{1,\ldots,p-1\}\right\},} %
  with \m{f\in I_{i,j}} represented by
    \Algn{\label{eq3t094u43u984398}\nonumber
       & \textstyle  \sum_{\TEMPORARYA ,\TEMPORARYA_1,\ldots,\TEMPORARYA_{n-1} \in \FF_p} \TEMPORARYA_1^{i_1}\cdots \TEMPORARYA_{n-1}^{i_{n-1}}
    f(-\TEMPORARYA,1)\begin{Lmatrix}
    1 & 0 \\ [\TEMPORARYA] + [\TEMPORARYA_1]p+\cdots+[\TEMPORARYA_{n-1}]p^{n-1} & 1
    \end{Lmatrix}  \sqbrackets{}{KZ,I(n)Z}{\FFp}{\GENN_j} \\ & \null\qquad + \textstyle \sum_{\TEMPORARYA_1,\ldots,\TEMPORARYA_{n-1} \in \FF_p} \TEMPORARYA_1^{i_1}\cdots \TEMPORARYA_{n-1}^{i_{n-1}}
    f(-1,0)\begin{Lmatrix}
 -[\TEMPORARYA_1]p-\cdots-[\TEMPORARYA_{n-1}]p^{n-1} & -1 \\ 1 & 0
\end{Lmatrix}  \sqbrackets{}{KZ,I(n)Z}{\FFp}{\GENN_j},
    }
where \m{i_{1},\ldots,i_{n-1} \in \{0,\ldots,p-1\}} %
are such that \m{i=i_1+\cdots+i_{n-1}p^{n-2}}.
\end{lemma}

\begin{proof}{Proof. }
Let %
 \m{h\in\SIGMA_r}. 
If \m{s\GEQSLANT p} then %
\algn{\textstyle x^{r-s}y^s = x^{r-s+p-1}y^{s-p+1} + (xy^p-x^py)x^{r-s-1}y^{s-p} \in x^{r-s+p-1}y^{s-p+1} + U.}
By repeatedly using this equation we can write \m{h} as \m{h_0+h^{\prime}}, where \m{h_0} is in the linear span of \m{\{x^{r-i}y^i \,|\, i\in\{0,\ldots,p-1\}\}} and \m{h^{\prime}\in U}. 
Since \m{x^r \in U}, it follows that \m{\SIGMA_r} is the linear span of \m{\textstyle U \cup \{x^{r-i}y^i \,|\, i \in \{1,\ldots,p-1\}\}}. %
The linear subspace \algn{\Span \{x^{r-i}y^i \,|\, i \in \{0,\ldots,p-1\}\} \subseteq \SIGMA_r} is a submodule of \m{\SIGMA_r}, %
 so
if, for \m{i \in \{1,\ldots,p-1\}}, \m{g_i=x^{r-i}y^i + U \in \QUOTIENT}, then there is a series
\algn{\QUOTIENT = \Span\{g_1,\ldots,g_{p-1}\}  \supseteq \Span\{g_1,\ldots,g_{p-2}\}\supseteq\cdots\supseteq \Span\{g_1\} \supseteq \{0\}.}
Since each of the factors of this series is generated by a single element and since \m{\dim\QUOTIENT = p-1}, it follows that all of these inclusions are proper. Therefore this is a series of
 \m{\QUOTIENT }  whose factors are 
one-dimensional, and one can find that they are  the representations \m{\textstyle\chi_{-\NEBENPOWEY}(\NEBENPOWEX),\ldots,\chi_{2-\NEBENPOWEY}(\NEBENPOWEX-1)}
 just as in the proof of lemma~\ref{LEMMAX2}. Part~\eqref{propeta3} of lemma~\ref{LEMMAX3.01} implies that
\m{g_i = x^{r-i}y^i + U = (-1)^i\GENN_i + U } for \m{i \in \{1,\ldots,p-1\}}.  The rest follows immediately from lemma~\ref{COROLEMMAX1}.
\end{proof}

\section{%
 Series expansions of elements of \texorpdfstring{\m{{\ind_{\INZA}^{G}\TILDE\SIGMA_r}}}{ind(\textSigma)}} \label{secseries}
 In this section we prove eight lemmas about the first few terms of a certain series expansion of
\algn{\textstyle %
\sum_{\TEMPORARYA \in \LCLASS}  {c_\TEMPORARYA \begin{Lmatrix}
    1 & 0 \\ \LIFTA{\TEMPORARYA} p^s & 1
    \end{Lmatrix}}\sqbrackets{}{\INZB}{\QQp}{\GEN}.}  
    Here \m{s\in\ZZ_{\GEQSLANT0}}, \m{c_\TEMPORARYA\in\ZZp[G]}, and \m{\LCLASS} is an arbitrary set of lifts (as defined in section~\ref{seccombdef}). %

\begin{lemma}\label{LEMMAX4}
Let \m{s \in\{0,\ldots,n-1\}} and let \m{\LCLASS} be an arbitrary set of lifts.
Recall that \m{\mathrm{err}_s = \QUARK^{p^{s+1}} + ap}.
For each \m{\TEMPORARYA \in\LCLASS} let \m{c_\TEMPORARYA\in\ZZp[G]}, and suppose that \algn{\textstyle \sum_{\TEMPORARYA \in\LCLASS} c_\TEMPORARYA = \BIGO{p}} in the sense that the coefficient of each \m{g\in G \subset \ZZp[G]} in \m{\textstyle \sum_{\TEMPORARYA \in\LCLASS} c_\TEMPORARYA} is \m{\BIGO{p}}. 
Then   \algn{%
    & \textstyle \sum_{\TEMPORARYA \in \LCLASS}  {c_\TEMPORARYA \begin{Lmatrix}
    1 & 0 \\ \LIFTA{\TEMPORARYA} p^s & 1
    \end{Lmatrix}}\sqbrackets{}{\INZB}{\QQp}{\GEN} 
    \\&\nonumber\textstyle\null\qquad \EQUIVZEROIDEAL a^{-1} \begin{Lmatrix}
    p & 0 \\ 0 & 1
    \end{Lmatrix} \sum_{\TEMPORARYA \in \LCLASS}  c_\TEMPORARYA       \begin{Lmatrix}
    1 & 0 \\ \LIFTA{\TEMPORARYA} p^{s+1} & 1
    \end{Lmatrix}{}\sqbrackets{}{\INZB}{\QQp}{ \sum_{\alpha=1}^{p-1}\CONST{\alpha}{s+1}_{\LIFTA{\TEMPORARYA}} \GEN_\alpha} + \BIGO{\mathrm{err}_sa^{-2}}.
}
\end{lemma}

\begin{proof}{Proof. }
To begin with, it is enough to prove the lemma when each \m{c_\TEMPORARYA} is in \m{\ZZp}, since then one can conclude the lemma for general \m{c_\TEMPORARYA \in \ZZp[G]} by 
writing
\algn{\textstyle c_\TEMPORARYA = \sum_{g\in G} c_{\TEMPORARYA}(g) g \in \ZZp[G],}
and by applying the lemma with the constants \m{c_{\TEMPORARYA}(g) \in\ZZp} 
for each \m{g\in G \subset \ZZp[G]}.
So let us assume that \m{c_\TEMPORARYA \in \ZZp}. Moreover, it is enough to prove the lemma when \m{\textstyle \sum_{\TEMPORARYA \in\LCLASS} c_\TEMPORARYA = 0}. This is because part~\eqref{comb10.1} of lemma~\ref{LEMMAX3.04} implies that \m{\BIGO{\mathrm{err}_sa^{-1}}  \supseteq \BIGO{p}}, so
 we can tweak \m{c_0} by adding a term \m{\BIGO{p}} to it to make \m{\sum_{\TEMPORARYA \in\LCLASS} c_\TEMPORARYA = 0} and in the process change the left and the right side of the main equation by terms that %
 belong to \m{\BIGO{\mathrm{err}_sa^{-2}}}. So let us  assume that \m{\sum_{\TEMPORARYA \in\LCLASS} c_\TEMPORARYA = 0}.
The proof proceeds by induction on \m{s}, starting from \m{s=n-1} and going down to \m{0}. Directly from the definition of \m{\HECKET} we have
\algn{
    \nonumber & \nonumber\textstyle \HECKET \left(1\sqbrackets{}{\INZB}{\QQp}{\GEN}\right)\nonumber %
     =  {\sum_{\MU\in\FF_p} \begin{Lmatrix}
    p &[\MU] \\ 0 & 1
    \end{Lmatrix} }\sqbrackets{}{\INZB}{\QQp}{\left(\begin{Lmatrix}
    1 & -[\MU] \\ 0 & p
    \end{Lmatrix} \cdot \GEN\right)}.
}
If \m{\MU\ne0} then
\algn{\textstyle \begin{Lmatrix}
    1 & -[\MU] \\ 0 & p
    \end{Lmatrix} \cdot \GEN  & \textstyle = (x^{p-1}-(py-[\MU]x)^{p-1})^2x^{r-2(p-1)}\\ & \textstyle = \left(p\binom{p-1}{1}[\MU]^{p-2}x^{p-2}y+\BIGO{p^2}\right)^2 x^{r-2(p-1)} = \BIGO{p^2},}
so that
\Algn{\label{eq3498ut4398jg394j3498}
     & \textstyle \HECKET \left(1\sqbrackets{}{\INZB}{\QQp}{\GEN}\right) %
     =  { \begin{Lmatrix}
    p & 0 \\ 0 & 1
    \end{Lmatrix} }\sqbrackets{}{\INZB}{\QQp}{\GEN(x,py)} + \BIGO{p^2}.
}
Therefore we have
   \algn{
    & \textstyle \HECKET \left(\sum_{\TEMPORARYA \in \LCLASS}  {c_\TEMPORARYA \begin{Lmatrix}
    1 & 0 \\ \TEMPORARYA p^{n-1} & 1
    \end{Lmatrix}}\sqbrackets{}{\INZB}{\QQp}{\GEN}\right)
    \\&\textstyle\null\qquad =  \sum_{\TEMPORARYA \in \LCLASS}  {c_\TEMPORARYA \begin{Lmatrix}
    p & 0 \\ \TEMPORARYA p^{n} & 1
    \end{Lmatrix}}\sqbrackets{}{\INZB}{\QQp}{ \GEN(x,p y)} + \BIGO{p^2}
    \\&\textstyle\null\qquad =  
     {\begin{Lmatrix}
    p & 0 \\ 0 & 1
    \end{Lmatrix}}\sqbrackets{}{\INZB}{\QQp}{\sum_{\TEMPORARYA \in \LCLASS}  c_\TEMPORARYA \begin{Lmatrix}
    1 & 0 \\ \TEMPORARYA p^{n} & 1
    \end{Lmatrix} \GEN(x,p y)} + \BIGO{p^2}     \\&\textstyle\null\qquad =  
     {\begin{Lmatrix}
    p & 0 \\ 0 & 1
    \end{Lmatrix}}\sqbrackets{}{\INZB}{\QQp}{\sum_{\TEMPORARYA \in \LCLASS}  c_\TEMPORARYA \GEN(x,p y)} + \BIGO{p^2} = \BIGO{p^2}.
  }
  The third equality is true because %
  \m{n\GEQSLANT 2}. %
  So
  \algn{
    & \textstyle \textstyle \sum_{\TEMPORARYA \in \LCLASS} { c_\TEMPORARYA \begin{Lmatrix}
    1 & 0 \\ \TEMPORARYA p^{n-1} & 1
    \end{Lmatrix}}\sqbrackets{}{\INZB}{\QQp}{\GEN}
    \\ & \null\qquad\textstyle = \BIGO{p^2a^{-1}} - (\HECKET-a)\left(a^{-1}\sum_{\TEMPORARYA \in \LCLASS} { c_\TEMPORARYA \begin{Lmatrix}
    1 & 0 \\ \TEMPORARYA p^{n-1} & 1
    \end{Lmatrix}}\sqbrackets{}{\INZB}{\QQp}{\GEN}\right),
  }
and therefore
 \Algn{\label{eqn33}%
    & \textstyle \sum_{\TEMPORARYA \in \LCLASS} { c_\TEMPORARYA \begin{Lmatrix}
    1 & 0 \\ \TEMPORARYA p^{n-1} & 1
    \end{Lmatrix}}\sqbrackets{}{\INZB}{\QQp}{\GEN} \EQUIVZEROIDEAL \BIGO{p^2a^{-1}}.
  }
This implies the base case \m{s=n-1}, and in fact it is a slightly stronger statement since \m{\BIGO{p^2a^{-1}} \subset \BIGO{\mathrm{err}_{n-1}a^{-2}}}.
 Now let \m{s\in\{0,\ldots,n-2\}} and suppose that the theorem is true for \m{s+1}. 
If \m{s<n-2} then \m{\textstyle\CONST{\alpha}{s+2}_\TEMPORARYA = \BIGO{\mathrm{err}_s}}  and therefore the induction hypothesis implies that
\Algn{\label{eq1423r23t34t}\textstyle  \sum_{\TEMPORARYA \in \LCLASS^\prime}  {c_{\TEMPORARYA} \begin{Lmatrix}
    1 & 0 \\ \TEMPORARYA p^{s+1} & 1
    \end{Lmatrix}}\sqbrackets{}{\INZB}{\QQp}{\GEN} \EQUIVZEROIDEAL \BIGO{\mathrm{err}_sa^{-1}}}
for any set of lifts \m{\LCLASS^\prime}. 
And, if \m{s=n-2} then \m{\textstyle\BIGO{p^2a^{-1}} \subset \BIGO{\mathrm{err}_sa^{-1}}} and therefore equation~\eqref{eq1423r23t34t} is implied by the stronger statement~\eqref{eqn33}. So equation~\eqref{eq1423r23t34t} is true regardless of whether \m{s<n-2} or \m{s=n-2}. %
 For \m{\MU\in\FF_p}, let \m{\LCLASS_\MU} be the set of lifts for which there is a bijection
 \algn{\TEMPORARYA\in\LCLASS \xmapsto{ \ \ } \TEMPORARYA_\MU=\TEMPORARYA/(1-\TEMPORARYA[\MU]p^{s+1}) \in \LCLASS_\MU.}
 If we apply equation~\eqref{eq1423r23t34t} with the set of lifts  \m{\LCLASS_\MU} and the constants \m{c_{\TEMPORARYA_\MU} = c_{\TEMPORARYA}},
 multiply it by \m{\begin{Smatrix}
    1 & [\MU] \\ 0 & 1
    \end{Smatrix}}, and then sum over all \m{\MU\in\FF_p}, we get that
    \algn{\textstyle \sum_{\MU\in\FF_p}\begin{Lmatrix}
    1 & [\MU] \\ 0 & 1
    \end{Lmatrix}\sum_{\TEMPORARYA_\MU \in \LCLASS_\MU}  {c_{\TEMPORARYA_\MU} \begin{Lmatrix}
    1 & 0 \\ \TEMPORARYA_\MU p^{s+1} & 1
    \end{Lmatrix}}\sqbrackets{}{\INZB}{\QQp}{\GEN} \EQUIVZEROIDEAL \BIGO{\mathrm{err}_sa^{-1}}.}
For \m{q=p^{s+1}} we have  \m{1+\TEMPORARYA_\MU[\MU]q = 1/(1-\TEMPORARYA [\MU]q)}. Then the equations
 \algn{\textstyle \begin{Lmatrix}
    1 & [\MU] \\ 0 & 1
    \end{Lmatrix}\begin{Lmatrix}
    1 & 0 \\ \TEMPORARYA_\MU q & 1
    \end{Lmatrix} = \begin{Lmatrix}
    1 & 0 \\  \TEMPORARYA_\MU q/(1+\TEMPORARYA_\MU [\MU]q) & 1
    \end{Lmatrix}\begin{Lmatrix}
    1 + \TEMPORARYA_\MU [\MU]q & [\MU] \\ 0 & 1/(1+\TEMPORARYA_\MU[\MU]q)
    \end{Lmatrix} = \begin{Lmatrix}
    1 & 0 \\  \TEMPORARYA  q & 1
    \end{Lmatrix}\begin{Lmatrix}
    1/(1-\TEMPORARYA [\MU]q) & [\MU] \\ 0 & 1-\TEMPORARYA[\MU]q
    \end{Lmatrix} }
 and \algn{\textstyle\begin{Lmatrix}
    1/(1-\TEMPORARYA [\MU]q) & [\MU] \\ 0 & 1-\TEMPORARYA[\MU]q
    \end{Lmatrix}\GEN
    = \VAREPSILONp^{-1}(1 - \TEMPORARYA[\MU]q)  \begin{Lmatrix}
    1 & [\MU] \\ 0 & 1
    \end{Lmatrix}\GEN + \BIGO{p}}
imply that
\Algn{\label{eqfj3948j3948j3}\textstyle  \sum_{\TEMPORARYA \in \LCLASS}  {c_{\TEMPORARYA} \begin{Lmatrix}
    1 & 0 \\ \TEMPORARYA p^{s+1} & 1
    \end{Lmatrix}}\sqbrackets{}{\INZB}{\QQp}{\left(\sum_{\MU\in\FF_p} \VAREPSILONp^{-1}(1 - \TEMPORARYA[\MU]q)  \begin{Lmatrix}
    1 & [\MU] \\ 0 & 1
    \end{Lmatrix}\GEN\right)} \EQUIVZEROIDEAL \BIGO{\mathrm{err}_sa^{-1}},}
    since \m{\BIGO{p} \subseteq \BIGO{\mathrm{err}_sa^{-1}}}.
Due to part~\eqref{propeta6} of lemma~\ref{LEMMAX3.01},
\algn{\textstyle \begin{Lmatrix}
    1 & [\MU] \\ 0 & 1
    \end{Lmatrix} \GEN = (1-[\MU]^{p-1})\GEN_0 - \sum_{\alpha=1}^{p-1} [\MU]^{p-1-\alpha} \GEN_\alpha + \BIGO{p}.}
    Since  \m{\sum_{\MU\in\FF_p} \VAREPSILONp^{-1}(1 - \TEMPORARYA[\MU]q)(1-[\MU]^{p-1})=1}, it follows that
\algn{\textstyle & \textstyle \sum_{\MU\in\FF_p} \VAREPSILONp^{-1}(1 - \TEMPORARYA[\MU]q)  \begin{Lmatrix}
    1 & [\MU] \\ 0 & 1
    \end{Lmatrix}\GEN
    \\ & \textstyle \null\qquad =
    \GEN_0- \sum_{\MU\in\FF_p} \VAREPSILONp^{-1}(1 - \TEMPORARYA[\MU]q)\sum_{\alpha=1}^{p-1}  [\MU]^{p-1-\alpha} \GEN_\alpha + \BIGO{p}
    \\ & \textstyle \null\qquad =
    \GEN_0 - \sum_{\alpha=1}^{p-1}   \CONST{\alpha}{s+1}_{\TEMPORARYA}  \GEN_\alpha + \BIGO{p}.}
Therefore, since \m{\BIGO{p} \subseteq \BIGO{\mathrm{err}_sa^{-1}}}, equation~\eqref{eqfj3948j3948j3} implies that
\Algn{\label{eqn44}%
 \textstyle  \sum_{\TEMPORARYA \in \LCLASS}  {c_{\TEMPORARYA} \begin{Lmatrix}
    1 & 0 \\ \TEMPORARYA p^{s+1} & 1
    \end{Lmatrix}}\sqbrackets{}{\INZB}{\QQp}{\left(\GEN_0 - \sum_{\alpha=1}^{p-1}   \CONST{\alpha}{s+1}_{\TEMPORARYA}  \GEN_\alpha\right)} \EQUIVZEROIDEAL \BIGO{\mathrm{err}_sa^{-1}}.}
Due to parts~\eqref{propeta2} and~\eqref{propeta5} of lemma~\ref{LEMMAX3.01} we have \m{\GEN(x,p y) = \GEN_0 + \BIGO{p}}, so
  \algn{
    & \textstyle \HECKET \left(\sum_{\TEMPORARYA \in \LCLASS}  {c_\TEMPORARYA \begin{Lmatrix}
    1 & 0 \\ \TEMPORARYA p^{s} & 1
    \end{Lmatrix}}\sqbrackets{}{\INZB}{\QQp}{\GEN}\right)
    \\&\textstyle\null\qquad =  \sum_{\TEMPORARYA \in \LCLASS}  {c_\TEMPORARYA \begin{Lmatrix}
    p & 0 \\ \TEMPORARYA p^{s+1} & 1
    \end{Lmatrix}}\sqbrackets{}{\INZB}{\QQp}{ \GEN(x,p y)} + \BIGO{p^2}
    \\&\textstyle\null\qquad =  
     \begin{Lmatrix}
    p & 0 \\ 0 & 1
    \end{Lmatrix} \sum_{\TEMPORARYA \in \LCLASS}  {c_\TEMPORARYA \begin{Lmatrix}
    1 & 0 \\ \TEMPORARYA p^{s+1} & 1
    \end{Lmatrix} }\sqbrackets{}{\INZB}{\QQp}{\GEN(x,p y)} + \BIGO{p^2}
    \\&\textstyle\null\qquad =  
     \begin{Lmatrix}
    p & 0 \\ 0 & 1
    \end{Lmatrix} \sum_{\TEMPORARYA \in \LCLASS}  {c_\TEMPORARYA \begin{Lmatrix}
    1 & 0 \\ \TEMPORARYA p^{s+1} & 1
    \end{Lmatrix} }\sqbrackets{}{\INZB}{\QQp}{ \GEN_0} + \BIGO{p}.
  }
Therefore
  \algn{
    & \textstyle  a^{-1}  \begin{Lmatrix}
    p & 0 \\ 0 & 1
    \end{Lmatrix} \sum_{\TEMPORARYA \in \LCLASS}  {c_\TEMPORARYA    \begin{Lmatrix}
    1 & 0 \\ \TEMPORARYA p^{s+1} & 1
    \end{Lmatrix}
}\sqbrackets{}{\INZB}{\QQp}{\GEN_0} %
    \EQUIVZEROIDEAL   
      \sum_{\TEMPORARYA \in \LCLASS}  {c_\TEMPORARYA \begin{Lmatrix}
    1 & 0 \\ \TEMPORARYA p^{s} & 1
    \end{Lmatrix}}\sqbrackets{}{\INZB}{\QQp}{\GEN} + \BIGO{pa^{-1}},
  }
which implies that equation~\eqref{eqn44} is equivalent to
\algn{
    & \textstyle  \sum_{\TEMPORARYA \in \LCLASS}  {c_\TEMPORARYA      \begin{Lmatrix}
    1 & 0 \\ \TEMPORARYA p^{s} & 1
    \end{Lmatrix}
}\sqbrackets{}{\INZB}{\QQp}{\GEN} 
    \\&\textstyle\null\qquad \EQUIVZEROIDEAL a^{-1}\begin{Lmatrix}
    p & 0 \\ 0 & 1
    \end{Lmatrix} \sum_{\TEMPORARYA \in \LCLASS}  c_\TEMPORARYA       \begin{Lmatrix}
    1 & 0 \\ \TEMPORARYA p^{s+1} & 1
    \end{Lmatrix}{}\sqbrackets{}{\INZB}{\QQp}{ \sum_{\alpha=1}^{p-1}\CONST{\alpha}{s+1}_{\TEMPORARYA} \GEN_\alpha} + \BIGO{\mathrm{err}_sa^{-2}}.
  }
 This concludes the induction step and completes the proof. 
\end{proof}

\begin{lemma}\label{LEMMAX6prime}
{Recall that, for \m{s\in\ZZ_{\GEQSLANT0}}, \m{\mathrm{err}_s = \QUARK^{p^{s+1}} + ap}. If \m{\beta\in\ZZ_p} then
\algn{{\begin{Lmatrix} 1 & 0 \\ \beta p & 1
\end{Lmatrix}}\sqbrackets{}{\INZB}{\QQp}{\GEN} \EQUIVZEROIDEAL 1 \sqbrackets{}{\INZB}{\QQp}{\GEN} + \BIGO{\mathrm{err}_{\PADICVALUATION(\beta)}a^{-1}}.}
}
\end{lemma}

\begin{proof}{Proof. } If \m{\PADICVALUATION(\beta)\GEQSLANT n-1} then the claim follows from part~\eqref{comb10.1} of lemma~\ref{LEMMAX3.04} and \algn{{\begin{Lmatrix} 1 & 0 \\ \delta p^n & 1
\end{Lmatrix}} \sqbrackets{}{\INZB}{\QQp}{\GEN} = 1 \sqbrackets{}{\INZB}{\QQp}{\GEN}+\BIGO{p}} for \m{\delta = \beta p^{1-n}\in\ZZ_p}. Let us suppose that \m{\PADICVALUATION(\beta)=s-1} for \m{s \in \{1,\ldots,n-1\}}, and let us write \m{\delta=\beta p^{1-s}\in\ZZ_p^\times}.  
Let us apply lemma~\ref{LEMMAX4} with a set of lifts \m{\LCLASS_\delta} that contains \m{0} and \m{\delta}, and
\algn{c_\TEMPORARYA=\left\{\begin{array}{rl} 1 & \text{if } \TEMPORARYA = \delta\text{,} \\ -1 & \text{if } \TEMPORARYA = 0\text{,} \\ 0 & \text{otherwise.} \end{array}\right.}
We get that 
\algn{{\begin{Lmatrix} 1 & 0 \\ \beta p & 1
\end{Lmatrix}}\sqbrackets{}{\INZB}{\QQp}{\GEN} - 1 \sqbrackets{}{\INZB}{\QQp}{\GEN} \EQUIVZEROIDEAL  \ERRR} for some %
 error term \m{\ERRR} that is explicitly described by lemma~\ref{LEMMAX4}. If \m{s<n-1} then \m{\ERRR=\BIGO{\mathrm{err}_{s-1}a^{-1}}} because all of the coefficients \m{a^{-1}\CONST{\alpha}{s+1}_\TEMPORARYA} appearing in \m{\ERRR} belong to \m{\textstyle %
\BIGO{\mathrm{err}_{s-1}a^{-1}}}. If \m{s=n-1} then equation~\eqref{eqn33} implies that \m{\ERRR} belongs to \m{\BIGO{p^2a^{-1}}\subset%
\BIGO{\mathrm{err}_{s-1}a^{-1}}}. %
\end{proof}

\begin{lemma}\label{LEMMAX6primeq}
{Recall that, for \m{s\in\ZZ_{\GEQSLANT0}}, \m{\mathrm{err}_s = \QUARK^{p^{s+1}} + ap}. If \m{\beta,\TEMPORARYA\in\ZZ_p}  then
\algn{{\begin{Lmatrix} 1 & 0 \\ \beta p & 1
\end{Lmatrix}}{\begin{Lmatrix} 1 & \TEMPORARYA \\ 0 & 1
\end{Lmatrix}}\sqbrackets{}{\INZB}{\QQp}{\GEN} \EQUIVZEROIDEAL \ZETA^{-\beta\TEMPORARYA} \begin{Lmatrix}
            1 & \TEMPORARYA  \\ 0 & 1
            \end{Lmatrix}  \sqbrackets{}{\INZB}{\QQp}{\GEN} + \BIGO{\mathrm{err}_{0}a^{-1}}.}
}
\end{lemma}

\begin{proof}{Proof. }
    Since
        \m{\textstyle \begin{Lmatrix}
            1 & 0 \\ \beta p & 1
   \end{Lmatrix}\begin{Lmatrix}
            1 & \TEMPORARYA  \\ 0 & 1
            \end{Lmatrix} = \begin{Lmatrix}
     1/(1+\beta\TEMPORARYA p) & \TEMPORARYA \\ 0 & 1+\beta\TEMPORARYA p
     \end{Lmatrix}\begin{Lmatrix}
            1 & 0  \\ \beta p/(1+\beta\TEMPORARYA p) & 1
            \end{Lmatrix},}
     \algn{\textstyle  \begin{Lmatrix}
            1 & 0 \\ \beta p & 1
   \end{Lmatrix}\begin{Lmatrix}
            1 & \TEMPORARYA  \\ 0 & 1
            \end{Lmatrix} \sqbrackets{}{\INZB}{\QQp}{\GEN} & \textstyle \EQUALZEROIDEAL \begin{Lmatrix}
     1/(1+\beta\TEMPORARYA p) & \TEMPORARYA \\ 0 & 1+\beta\TEMPORARYA p
     \end{Lmatrix}\begin{Lmatrix}
            1 & 0  \\ \beta p/(1+\beta\TEMPORARYA p) & 1
            \end{Lmatrix} \sqbrackets{}{\INZB}{\QQp}{\GEN}  
          \\& \textstyle \EQUIVZEROIDEAL \begin{Lmatrix}
     1/(1+\beta\TEMPORARYA p) & \TEMPORARYA \\ 0 & 1+\beta\TEMPORARYA p
          \end{Lmatrix} \sqbrackets{}{\INZB}{\QQp}{\GEN} + \BIGO{\mathrm{err}_{\PADICVALUATION(\beta)}a^{-1}}
          \\ & \textstyle \textstyle \EQUALZEROIDEAL \ZETA^{-\beta\TEMPORARYA} \begin{Lmatrix}
            1 & \TEMPORARYA  \\ 0 & 1
            \end{Lmatrix}  \sqbrackets{}{\INZB}{\QQp}{\GEN} + \BIGO{\mathrm{err}_{0}a^{-1}}  .}
          The congruence follows from lemma~\ref{LEMMAX6prime}, and the second equality follows from %
\Algn{\label{eqf34j49j5g98jg}\textstyle \VAREPSILONp^{-1}(1+\beta \TEMPORARYA p) = \ZETA^{-\beta\TEMPORARYA} + \BIGO{1-\ZETA^p}=\ZETA^{-\beta\TEMPORARYA} + \BIGO{\mathrm{err}_{0}a^{-1}}.}
\end{proof}

\begin{lemma}\label{LEMMAX6primep}
Let \m{w \in \{1,\ldots,p-1\}}. Recall that \m{\mathrm{err}_0 = \QUARK^{p} + ap}. If \m{\TEMPORARYA\in\ZZ_p} then
\algn{& \textstyle{}\begin{Lmatrix}
1 &  0 \\ \TEMPORARYA p & 1
\end{Lmatrix} \sqbrackets{}{\INZB}{\QQp}{\GEN_w} %
 \EQUIVZEROIDEAL -{1}\sqbrackets{}{\INZB}{\QQp}{\sum_{l=0}^{p-2} %
\CONST{l}{1}_{-\TEMPORARYA}\GEN_{\MORE{w+l}}} + \BIGO{\mathrm{err}_0a^{-1}}.}
\end{lemma}

\begin{proof}{Proof. }
For \m{\MU \in \FF_p}, let us apply lemma~\ref{LEMMAX6prime} with \m{\beta  = -\TEMPORARYA/(1+\TEMPORARYA[\MU]p)}, multiply it by \m{\textstyle \begin{Lmatrix}
1 &  0 \\ \TEMPORARYA p & 1
\end{Lmatrix}  [\MU]^w \begin{Lmatrix}
1 & [\MU] \\ 0 & 1
\end{Lmatrix} \in \ZZp[G]},
 and then sum over all \m{ {\MU\in\FF_p}}.
 On the right side we get
\Algn{\label{eqf3948jf39juf398}
 &    \textstyle \begin{Lmatrix}
1 &  0 \\ \TEMPORARYA p & 1
\end{Lmatrix} \sum_{\MU\in\FF_p}  [\MU]^w \begin{Lmatrix}
1 & [\MU] \\ 0 & 1
\end{Lmatrix} \sqbrackets{}{\INZB}{\QQp}{\GEN}
 = \textstyle   \begin{Lmatrix}
1 &  0 \\ \TEMPORARYA p & 1
\end{Lmatrix} \sqbrackets{}{\INZB}{\QQp}{\GEN_w}.
}
On the left side we  get
\Algn{\label{eqf34ojg3ogij}
 &    \textstyle \begin{Lmatrix}
1 &  0 \\ \TEMPORARYA p & 1
\end{Lmatrix} \sum_{\MU\in\FF_p}  [\MU]^w \begin{Lmatrix}
1 & [\MU] \\ 0 & 1
\end{Lmatrix} {\begin{Lmatrix} 1 & 0 \\ -\TEMPORARYA p/(1+\TEMPORARYA[\MU]p) & 1
\end{Lmatrix}}\sqbrackets{}{\INZB}{\QQp}{\GEN}\nonumber
\\ & \null\qquad = \textstyle  \begin{Lmatrix}
1 &  0 \\ \TEMPORARYA p & 1
\end{Lmatrix} {\begin{Lmatrix} 1 & 0 \\ -\TEMPORARYA p & 1
\end{Lmatrix}} \sum_{\MU\in\FF_p}  [\MU]^w \begin{Lmatrix}
1 & [\MU] \\ 0 & 1
\end{Lmatrix}\VAREPSILONp^{-1}(1+\TEMPORARYA[\MU]p) \sqbrackets{}{\INZB}{\QQp}{\GEN} + \BIGO{p}\nonumber
\\ & \null\qquad = \textstyle  1 \sqbrackets{}{\INZB}{\QQp}{\sum_{\MU\in\FF_p}  [\MU]^w \begin{Lmatrix}
1 & [\MU] \\ 0 & 1
\end{Lmatrix} \VAREPSILONp^{-1}(1+\TEMPORARYA[\MU]p)\GEN} + \BIGO{p}\nonumber
\\ & \null\qquad =\textstyle -{1}\sqbrackets{}{\INZB}{\QQp}{\sum_{l=0}^{p-2} \sum_{\MU\in\FF_p} [\MU]^{p-1-l}\VAREPSILONp^{-1}(1+\TEMPORARYA[\MU]p)\GEN_{\MORE{w+l}}} + \BIGO{p}.
}
The first equality follows from the equation
\algn{\textstyle \begin{Lmatrix}
    1 & [\MU] \\ 0 & 1
    \end{Lmatrix}\begin{Lmatrix}
    1 & 0 \\ -\TEMPORARYA p/(1+\TEMPORARYA[\MU]p) & 1
    \end{Lmatrix} = \begin{Lmatrix}
    1 & 0 \\  -\TEMPORARYA p & 1
    \end{Lmatrix}\begin{Lmatrix}
    1/(1+\TEMPORARYA [\MU]p) & [\MU] \\ 0 & 1+\TEMPORARYA[\MU]p
    \end{Lmatrix}.} %
 The third equality follows from part~\eqref{propeta6} of lemma~\ref{LEMMAX3.02}, which implies that
    \algn{\textstyle \begin{Lmatrix}
    1 & [\MU] \\ 0 & 1
    \end{Lmatrix} \GEN = (1-[\MU]^{p-1})\GEN_0 - \sum_{u=1}^{p-1} [\MU]^{p-1-u} \GEN_u + \BIGO{p}}
    and, since \m{w>0} and therefore
 \m{([\MU]^w-[\MU]^{p-1+w}) \GEN_0 =    0},  that \algn{
        & \textstyle \sum_{\MU\in\FF_p}  [\MU]^w \begin{Lmatrix}
1 & [\MU] \\ 0 & 1
\end{Lmatrix} \VAREPSILONp^{-1}(1+\TEMPORARYA[\MU]p)\GEN\nonumber
\\ & \textstyle \null\qquad = 
-\sum_{\MU\in\FF_p}  \sum_{u=1}^{p-1} [\MU]^{p-1-u+w} \VAREPSILONp^{-1}(1+\TEMPORARYA[\MU]p)\GEN_u + \BIGO{p}\nonumber
\\ & \textstyle \null\qquad = 
-\sum_{\MU\in\FF_p} \sum_{u=w}^{p-2+w} [\MU]^{p-1-u+w} \VAREPSILONp^{-1}(1+\TEMPORARYA[\MU]p)\GEN_{\MORE{u}} + \BIGO{p}\nonumber
\\ & \textstyle \null\qquad = 
-\sum_{\MU\in\FF_p} \sum_{l=0}^{p-2} [\MU]^{p-1-l} \VAREPSILONp^{-1}(1+\TEMPORARYA[\MU]p)\GEN_{\MORE{w+l}} + \BIGO{p}.
    }
Here for the second equality we %
  change the range of summation from \m{1\LEQSLANT u\LEQSLANT p-1} to \m{w \LEQSLANT u\LEQSLANT p-2+w} (which we can do since \m{[\MU]^{p-1-u+w} \VAREPSILONp^{-1}(1+\TEMPORARYA[\MU]p)\GEN_{\MORE{u}}}   depends only on the congruence class of \m{u} modulo \m{p-1}), and for the third equality we change the variable \m{u} to \m{l=u-w}. Equations~\eqref{eqf3948jf39juf398} and~\eqref{eqf34ojg3ogij} complete the proof.
\end{proof}

\begin{lemma}\label{LEMMAX6}
{Let \m{\alpha\in\{1,\ldots,p-1\}} and let \m{\LCLASS} be an arbitrary set of lifts.  Recall that  %
 \m{\mathrm{err}_0 = \QUARK^{p} + ap}. Let \m{\CONSTANT{\alpha}=-\frac{(1-\ZETA)^{\alpha}}{\alpha!}}. Then
\algn{ & \textstyle \sum_{\TEMPORARYA \in \LCLASS} {\TEMPORARYA^{p-1-\alpha} \begin{Lmatrix}
    1 & 0 \\ \TEMPORARYA & 1
    \end{Lmatrix}}\sqbrackets{}{\INZB}{\QQp}{\GEN}%
    \EQUIVZEROIDEAL \textstyle %
    (1+\BIGO{\QUARK})\COVERA{\alpha} \begin{Lmatrix} p & 0 \\  0 & 1
\end{Lmatrix} \sqbrackets{}{\INZB}{\QQp}{\GEN_\alpha}+\BIGO{\mathrm{err}_0a^{-2}}.}
}
\end{lemma}

\begin{proof}{Proof. }
 Due to lemma~\ref{LEMMAX4} with \m{s=0}, the set of lifts \m{\LCLASS}, and \m{c_\TEMPORARYA = \TEMPORARYA^{p-1-\alpha}} for \m{\TEMPORARYA\in\LCLASS},
\algn{
    & \textstyle \sum_{\TEMPORARYA \in \LCLASS} {\TEMPORARYA^{p-1-\alpha} \begin{Lmatrix}
    1 & 0 \\ \TEMPORARYA & 1
    \end{Lmatrix}}\sqbrackets{}{\INZB}{\QQp}{\GEN} 
    \\&\textstyle\null\qquad \EQUIVZEROIDEAL {a^{-1}\begin{Lmatrix}
    p & 0 \\ 0 & 1
    \end{Lmatrix}\sum_{w=1}^{p-1} \sum_{\TEMPORARYA\in\LCLASS} \TEMPORARYA^{p-1-\alpha}\begin{Lmatrix}
    1 & 0 \\ \TEMPORARYA p & 1
    \end{Lmatrix} \left({\sum_{\LAMBDA\in\FF_p} [\LAMBDA]^{p-1-w}\ZETA^{\TEMPORARYA[\LAMBDA]} }\right)}\sqbrackets{}{\INZB}{\QQp}{\GEN_{w}} \\ & \textstyle{}\null\qquad\qquad\qquad\qquad\qquad\qquad\qquad\qquad\qquad\qquad\qquad\qquad\qquad\quad\qquad \ + \BIGO{\mathrm{err}_0a^{-2}}
    \\&\textstyle\null\qquad \EQUIVZEROIDEAL
    {a^{-1}\begin{Lmatrix}
    p & 0 \\ 0 & 1
    \end{Lmatrix}\sum_{j=1}^{p-1} C_{\alpha,j} }\sqbrackets{}{\INZB}{\QQp}{\GEN_{j}} + \BIGO{\mathrm{err}_0a^{-2}}.
  }
For the first congruence we use equation~\eqref{eqf34j49j5g98jg} to replace \m{\CONST{w}{1}_\TEMPORARYA} with \algn{\sum_{\LAMBDA\in\FF_p} [\LAMBDA]^{p-1-w}\ZETA^{\TEMPORARYA[\LAMBDA]}+\BIGO{\mathrm{err}_0a^{-1}}.} For the
second congruence we use lemma~\ref{LEMMAX6primep} to replace  \m{\begin{Smatrix}
    1 & 0 \\ \TEMPORARYA p & 1
    \end{Smatrix} \sqbrackets{}{\INZB}{\QQp}{\GEN_{w}} } with 
    \algn{\textstyle -{1}\sqbrackets{}{\INZB}{\QQp}{\sum_{l=0}^{p-2} \sum_{\MU\in\FF_p} [\MU]^{p-1-l}\ZETA^{-\TEMPORARYA[\MU]}\GEN_{\MORE{w+l}}} + \BIGO{\mathrm{err}_0a^{-1}}.}
    Since for each  \m{j\in\{1,\ldots,p-1\}} and for each \m{w \in\{1,\ldots,p-1\}} there is precisely one  \m{l\in\{0,\ldots,p-2\}} such that 
    \m{\MORE{w+l}=j \in\{1,\ldots,p-1\}}, we can write %
    \m{C_{\alpha,j}} as 
\algn{\textstyle
C_{\alpha,j} & \textstyle = -\sum_{w=1}^{p-1}\sum_{\TEMPORARYA\in\LCLASS}\sum_{\LAMBDA,\MU\in\FF_p} \TEMPORARYA^{p-1-\alpha} { [\LAMBDA]^{p-1-w} }  [\MU]^{\MORE{w-j}}\ZETA^{\TEMPORARYA([\LAMBDA]-[\MU])}+ \BIGO{\mathrm{err}_0a^{-1}}.
}
We can change the set of lifts \m{\LCLASS} with the set \m{\LCLASSTEICH} of Teichm\"uller lifts in this equation, as that only changes the right side by a term that belongs to \m{\BIGO{\mathrm{err}_0a^{-1}}}.
Then the term in ``\m{\sum_{\TEMPORARYA\in\LCLASSTEICH}}'' such that \m{\TEMPORARYA=0} is \algn{\textstyle -\VERIFYIF{\alpha=p-1}\sum_{w=1}^{p-1} {\sum_{\LAMBDA,\MU\in\FF_p} [\LAMBDA]^{p-1-w} }  [\MU]^{\MORE{w-j}} & \textstyle = \VERIFYIF{\alpha=p-1}(1-p) \sum_{\LAMBDA \in\FF_p} [\LAMBDA]^{p-1-j} \\ & \textstyle= \BIGO{\mathrm{err}_0a^{-1}}.} So if \m{\alpha,j \in\{1,\ldots, p-1\}}  then \m{C_{\alpha,j}} is equal to
\algn{\textstyle
& \textstyle 
-\sum_{w=1}^{p-1}\sum_{\TEMPORARYA\in \FF_p^\times} \sum_{\LAMBDA,\MU\in\FF_p} [\TEMPORARYA]^{p-1-\alpha} { [\LAMBDA]^{p-1-w} }  [\MU]^{\MORE{w-j}}\ZETA^{[\TEMPORARYA\LAMBDA]-[\TEMPORARYA\MU]}+ \BIGO{\mathrm{err}_0a^{-1}}
\\ & \textstyle \null\qquad =
-\sum_{w=1}^{p-1}\sum_{\TEMPORARYA\in \FF_p^\times} \sum_{\LAMBDA,\MU\in\FF_p} [\TEMPORARYA]^{{j-\alpha}} { [\TEMPORARYA\LAMBDA]^{p-1-w} }  [\TEMPORARYA\MU]^{\MORE{w-j}}\ZETA^{[\TEMPORARYA\LAMBDA]-[\TEMPORARYA\MU]}+ \BIGO{\mathrm{err}_0a^{-1}}
\\ & \textstyle \null\qquad = -\sum_{w=1}^{p-1}\sum_{\TEMPORARYA\in \FF_p^\times} \sum_{\TEMPORARYA_\LAMBDA,\TEMPORARYA_\MU\in\FF_p} [\TEMPORARYA]^{{j-\alpha}} { [\TEMPORARYA_\LAMBDA]^{p-1-w} }  [-\TEMPORARYA_\MU]^{\MORE{w-j}}\ZETA^{[\TEMPORARYA_\LAMBDA]+[\TEMPORARYA_\MU]}+ \BIGO{\mathrm{err}_0a^{-1}}
\\ & \textstyle \null\qquad = \VERIFYIF{j=\alpha}\sum_{w=1}^{p-1}\sum_{\TEMPORARYA_\LAMBDA,\TEMPORARYA_\MU\in\FF_p} (-1)^{j-w}{ [\TEMPORARYA_\LAMBDA]^{p-1-w} }  [\TEMPORARYA_\MU]^{\MORE{w-j}}\ZETA^{[\TEMPORARYA_\LAMBDA]+[\TEMPORARYA_\MU]}+ \BIGO{\mathrm{err}_0a^{-1}}.
}
For the second equality we change the variable \m{\LAMBDA} to \m{\TEMPORARYA_\LAMBDA=\TEMPORARYA\LAMBDA} and the variable \m{\MU} to \m{\TEMPORARYA_\MU = -\TEMPORARYA\MU}. %
 For the third equality we use the equation \algn{\textstyle \sum_{\TEMPORARYA\in\FF_p^\times} [\TEMPORARYA]^{{j-\alpha}} = \left\{\begin{array}{rl} p-1 & \text{if } p-1\mid j-\alpha\text{,} \\ 0 & \text{otherwise.} \end{array}\right.}
Therefore \m{C_{\alpha,j}=\BIGO{\mathrm{err}_0a^{-1}}} when \m{\alpha\ne j}, and
     \algn{C_{\alpha,\alpha} & \textstyle{} = \sum_{w=1}^{p-1} (-1)^{\alpha-w}\CONST{w}{1}_1\CONST{\LESS{\alpha-w}}{1}_1 + \BIGO{\mathrm{err}_0a^{-1}} \\ & = \textstyle{} (-1)^\alpha \sum_{w=1}^{\alpha} (-1)^{w}\CONST{w}{1}_1\CONST{{\alpha-w}}{1}_1 + \BIGO{\mathrm{err}_0a^{-1}}
    \\ & = \textstyle{}
    C_\alpha(1+\BIGO{\QUARK})+\BIGO{\mathrm{err}_0a^{-1}}.} For the second equality we can change the range of summation  because part~\eqref{gammaprop3} of lemma~\ref{LEMMAX3.02} implies that  \m{\TADICVALUATION(\CONST{\beta}{1}_1)\GEQSLANT\beta} for \m{\beta\in\{0,\ldots,p-1\}}. The third equality follows from lemma~\ref{LEMMAX3.03}. %
\end{proof}

\begin{lemma}\label{LEMMAX5}
{Let \m{u \in\ZZ_{\GEQSLANT 0}} and \m{s \in \{1,\ldots,p-1\}}. For \m{\TEMPORARYA\in\FF_p}, let
\algn{c_\TEMPORARYA\textstyle = (-1)^{u+\NEBENPOWEX} \sum_{\MU \in \FF_p} \sum_{j=0}^{u} \binom{u}{j}[\MU]^{u-j} [\TEMPORARYA]^{\LESS{-j+s-\NEBENPOWEY}} \begin{Lmatrix}
        0 & -1 \\ 1 & [\MU]
        \end{Lmatrix} \in \ZZ_p[G].}
Then
\Algn{\label{eqf34398jf983j498fj934}\nonumber
& \textstyle\sum_{\TEMPORARYA \in \FF_p} {c_\TEMPORARYA \begin{Lmatrix}
    1 & 0 \\ [\TEMPORARYA] & 1
    \end{Lmatrix}}\sqbrackets{}{\INZB}{\QQp}{\GEN} \\ \textstyle & \textstyle\qquad  =\nonumber
\sum_{\TEMPORARYA \in \FF_p} {[\TEMPORARYA]^{u} \begin{Lmatrix}
    1 & 0 \\ [\TEMPORARYA] & 1
    \end{Lmatrix}}\sqbrackets{}{\INZB}{\QQp}{\GEN_s} \\& \null\qquad\qquad\textstyle + (-1)^{u+\NEBENPOWEX} \sum_{\beta \equiv s-\NEBENPOWEY \bmod{p-1}}{\binom{u}{\beta}\begin{Lmatrix}
0 & -1 \\ 1 & 0
\end{Lmatrix}}\sqbrackets{}{\INZB}{\QQp}{\GEN_{u-\beta}} + \ERRR,}
where, for any sets of lifts \m{(\LCLASS_\TEMPORARYAA)_{\TEMPORARYAA\in\FF_p^\times}},
\Algn{\label{eq0f3j394jf9834jf39}\nonumber\ERRR & \textstyle\vphantom{(\ZETA-1) \sum_{\TEMPORARYA \in \LCLASS,\,\TEMPORARYAA\in\FF_p^\times} \sum_{j=1}^{p-1} 
        \frac{(-1)^j}{j} \TEMPORARYA^{u+j} [\TEMPORARYAA]^{\NEBENPOWEX-s-j} \begin{Lmatrix}
        1 & 0 \\  \TEMPORARYA & 1
        \end{Lmatrix}\begin{Lmatrix}
        1/[\TEMPORARYAA] & 1 \\ 0 & [\TEMPORARYAA]
        \end{Lmatrix}\sqbrackets{}{\INZB}{\QQp}{\GEN} + \BIGO{\QUARK^2}} \EQUALZEROIDEAL \sum_{\TEMPORARYAA\in\FF_p^\times,\,\LAMBDA \in \FF_p} { [\LAMBDA]^u [\TEMPORARYAA]^{{\NEBENPOWEX-s}} \left(\begin{Lmatrix}
        1 & 0 \\  [\LAMBDA+\TEMPORARYAA]-[\TEMPORARYAA] & 1
        \end{Lmatrix}-\begin{Lmatrix}
        1 & 0 \\  [\LAMBDA] & 1
        \end{Lmatrix}\right)\begin{Lmatrix}
        1/[\TEMPORARYAA] & 1 \\ 0 & [\TEMPORARYAA]
        \end{Lmatrix}}\sqbrackets{}{\INZB}{\QQp}{\GEN} + \BIGO{p} \\ &\nonumber \textstyle\vphantom{(\ZETA-1) \sum_{\TEMPORARYA \in \LCLASS,\,\TEMPORARYAA\in\FF_p^\times} \sum_{j=1}^{p-1} 
        \frac{(-1)^j}{j} \TEMPORARYA^{u+j} [\TEMPORARYAA]^{\NEBENPOWEX-s-j} \begin{Lmatrix}
        1 & 0 \\  \TEMPORARYA & 1
        \end{Lmatrix}\begin{Lmatrix}
        1/[\TEMPORARYAA] & 1 \\ 0 & [\TEMPORARYAA]
        \end{Lmatrix}\sqbrackets{}{\INZB}{\QQp}{\GEN} + \BIGO{\QUARK^2}} \EQUIVZEROIDEAL (\ZETA-1) \sum_{\TEMPORARYAA\in\FF_p^\times,\,\TEMPORARYA \in \LCLASS_\TEMPORARYAA} \sum_{j=1}^{p-1} 
        \frac{(-1)^j}{j} \TEMPORARYA^{u+j} [\TEMPORARYAA]^{\NEBENPOWEX-s-j} \begin{Lmatrix}
        1 & 0 \\  \TEMPORARYA & 1
        \end{Lmatrix}\begin{Lmatrix}
        1/[\TEMPORARYAA] & 1 \\ 0 & [\TEMPORARYAA]
        \end{Lmatrix}\sqbrackets{}{\INZB}{\QQp}{\GEN} + \BIGO{\QUARK^2} \\ & \textstyle\vphantom{(\ZETA-1) \sum_{\TEMPORARYA \in \LCLASS,\,\TEMPORARYAA\in\FF_p^\times} \sum_{j=1}^{p-1} 
        \frac{(-1)^j}{j} \TEMPORARYA^{u+j} [\TEMPORARYAA]^{\NEBENPOWEX-s-j} \begin{Lmatrix}
        1 & 0 \\  \TEMPORARYA & 1
        \end{Lmatrix}\begin{Lmatrix}
        1/[\TEMPORARYAA] & 1 \\ 0 & [\TEMPORARYAA]
        \end{Lmatrix}\sqbrackets{}{\INZB}{\QQp}{\GEN} + \BIGO{\QUARK^2}} \EQUALZEROIDEAL \BIGO{\QUARK}. }
}\end{lemma}

\begin{proof}{Proof. }
For \m{\TEMPORARYA\in\FF_p}, let us write \m{B(\TEMPORARYA) = c_\TEMPORARYA \begin{Smatrix}
1 & 0 \\ [\TEMPORARYA] & 1
\end{Smatrix}.} %
 If \m{\beta\in\{1,\ldots,p-2\}} then \m{[0]^\beta = 0}, and \m{[0]^0 = 1}.
 Therefore
\Algn{\label{eqt43gj49jg4j59}
& \textstyle B(0)\sqbrackets{}{\INZB}{\QQp}{\GEN} \nonumber\\ & \textstyle\null\qquad = \textstyle \sum_{\MU \in \FF_p} {(-1)^{u+\NEBENPOWEX} \sum_{j} \binom{u}{j}[\MU]^{u-j} [0]^{\LESS{-j+s-\NEBENPOWEY}} \begin{Lmatrix}
        0 & -1 \\ 1 & [\MU]
        \end{Lmatrix}}\sqbrackets{}{\INZB}{\QQp}{\GEN}\nonumber
        \\ & \textstyle\null\qquad= \textstyle  
        \textstyle \sum_{\MU \in \FF_p} {(-1)^{u+\NEBENPOWEX} \sum_{\beta \equiv s-\NEBENPOWEY \bmod{p-1}} \binom{u}{\beta}[\MU]^{u-\beta} \begin{Lmatrix}
        0 & -1 \\ 1 & [\MU]
        \end{Lmatrix}}\sqbrackets{}{\INZB}{\QQp}{\GEN}\nonumber
  \\ & \textstyle\null\qquad =(-1)^{u+\NEBENPOWEX} \sum_{\beta \equiv s-\NEBENPOWEY \bmod{p-1}}{\binom{u}{\beta}\begin{Lmatrix}
0 & -1 \\ 1 & 0
\end{Lmatrix}}\sqbrackets{}{\INZB}{\QQp}{\sum_{\MU \in \FF_p}[\MU]^{u-\beta}\begin{Lmatrix}
        1 & [\MU] \\ 0 & 1
        \end{Lmatrix} \GEN}   \nonumber\\ & \textstyle\null\qquad =(-1)^{u+\NEBENPOWEX} \sum_{\beta \equiv s-\NEBENPOWEY \bmod{p-1}}{\binom{u}{\beta}\begin{Lmatrix}
0 & -1 \\ 1 & 0
\end{Lmatrix}}\sqbrackets{}{\INZB}{\QQp}{\GEN_{u-\beta}}.
}
 We also have
            \Algn{\label{eq12rf0jwf98j3498fj3984jf3498j}\nonumber
\textstyle & \textstyle  \sum_{\TEMPORARYA\in\FF_p^\times}B(\TEMPORARYA)\sqbrackets{}{\INZB}{\QQp}{\GEN}
\\ &\nonumber \textstyle \null\qquad =
\sum_{\TEMPORARYA\in\FF_p^\times,\,\MU \in \FF_p} {(-1)^{u+\NEBENPOWEX} \sum_{j} \binom{u}{j}[\MU]^{u-j} [\TEMPORARYA]^{{-j+s-\NEBENPOWEY}} \begin{Lmatrix}
        0 & -1 \\ 1 & [\MU]
        \end{Lmatrix} \begin{Lmatrix}
    1 & 0 \\ [\TEMPORARYA] & 1
    \end{Lmatrix}}\sqbrackets{}{\INZB}{\QQp}{\GEN}
        \\ &\nonumber \textstyle \null\qquad =
\sum_{\TEMPORARYA\in\FF_p^\times,\,\MU \in \FF_p} {(-1)^{u} \sum_{j} \binom{u}{j}[\MU]^{u-j} [\TEMPORARYA]^{{-j+s-\NEBENPOWEY}} \begin{Lmatrix}
        1 & 0 \\  -[\MU]-1/[\TEMPORARYA] & 1
        \end{Lmatrix}\begin{Lmatrix}
        [\TEMPORARYA] & 1 \\ 0 & 1/[\TEMPORARYA]
        \end{Lmatrix}}\sqbrackets{}{\INZB}{\QQp}{\GEN}
               \\ &\nonumber \textstyle \null\qquad =
\sum_{\TEMPORARYAA\in\FF_p^\times,\,\MU \in \FF_p} {(-1)^{u} \sum_{j} \binom{u}{j}[\MU]^{u-j} [\TEMPORARYAA]^{{j-s+\NEBENPOWEX}} \begin{Lmatrix}
        1 & 0 \\  -[\MU]-[\TEMPORARYAA] & 1
        \end{Lmatrix}\begin{Lmatrix}
        1/[\TEMPORARYAA] & 1 \\ 0 & [\TEMPORARYAA]
        \end{Lmatrix}}\sqbrackets{}{\INZB}{\QQp}{\GEN}
                       \\ &\nonumber \textstyle \null\qquad =
\sum_{\TEMPORARYAA\in\FF_p^\times,\,\LAMBDA \in \FF_p} {(-1)^{u} \sum_{j} \binom{u}{j}[-\LAMBDA-\TEMPORARYAA]^{u-j} [\TEMPORARYAA]^{{j-s+\NEBENPOWEX}} \begin{Lmatrix}
        1 & 0 \\  [\LAMBDA+\TEMPORARYAA]-[\TEMPORARYAA] & 1
        \end{Lmatrix}\begin{Lmatrix}
        1/[\TEMPORARYAA] & 1 \\ 0 & [\TEMPORARYAA]
        \end{Lmatrix}}\sqbrackets{}{\INZB}{\QQp}{\GEN}
 \\ &\nonumber \textstyle \null\qquad =
\sum_{\TEMPORARYAA\in\FF_p^\times,\,\LAMBDA \in \FF_p} { [\LAMBDA]^u [\TEMPORARYAA]^{{\NEBENPOWEX-s}} \begin{Lmatrix}
        1 & 0 \\  [\LAMBDA+\TEMPORARYAA]-[\TEMPORARYAA] & 1
        \end{Lmatrix}\begin{Lmatrix}
        1/[\TEMPORARYAA] & 1 \\ 0 & [\TEMPORARYAA]
        \end{Lmatrix}}\sqbrackets{}{\INZB}{\QQp}{\GEN}+ \BIGO{p}
        \\ &\nonumber \textstyle \null\qquad =
\sum_{\TEMPORARYAA\in\FF_p^\times,\,\LAMBDA \in \FF_p} { [\LAMBDA]^u [\TEMPORARYAA]^{{\NEBENPOWEX-s}} \begin{Lmatrix}
        1 & 0 \\ [\LAMBDA] & 1
        \end{Lmatrix}\begin{Lmatrix}
        1/[\TEMPORARYAA] & 1 \\ 0 & [\TEMPORARYAA]
        \end{Lmatrix}}\sqbrackets{}{\INZB}{\QQp}{\GEN} + \ERRR+\BIGO{p} \nonumber
                \\ &\nonumber \textstyle \null\qquad =
\sum_{\LAMBDA \in \FF_p} { [\LAMBDA]^u \begin{Lmatrix}
        1 & 0 \\ [\LAMBDA] & 1
        \end{Lmatrix} }\sqbrackets{}{\INZB}{\QQp}{\sum_{\TEMPORARYAA\in\FF_p^\times} [\TEMPORARYAA]^{{\NEBENPOWEX-s}} \begin{Lmatrix}
        1/[\TEMPORARYAA] & 1 \\ 0 & [\TEMPORARYAA]
        \end{Lmatrix} \GEN} +\ERRR+\BIGO{p}  \nonumber  \\ &\nonumber \textstyle \null\qquad =
\sum_{\LAMBDA \in \FF_p} { [\LAMBDA]^u \begin{Lmatrix}
        1 & 0 \\ [\LAMBDA] & 1
        \end{Lmatrix} }\sqbrackets{}{\INZB}{\QQp}{\sum_{\TEMPORARYA\in\FF_p^\times} [\TEMPORARYA]^{s} \begin{Lmatrix}
        1 & [\TEMPORARYA] \\ 0 & 1
        \end{Lmatrix}\GEN} + \ERRR+\BIGO{p}\nonumber
          \\ & \textstyle \null\qquad =
\sum_{\LAMBDA \in \FF_p} { [\LAMBDA]^u \begin{Lmatrix}
        1 & 0 \\ [\LAMBDA] & 1
        \end{Lmatrix} }\sqbrackets{}{\INZB}{\QQp}{\GEN_s} + \ERRR+\BIGO{p}.} The first equality follows directly from the definition of \m{B(\TEMPORARYA)}. The second equality follows from  the fact that for \m{\TEMPORARYA\in\FF_p^\times} we have
\algn{\begin{Lmatrix}
        0 & -1 \\ 1 & [\MU]
        \end{Lmatrix} \begin{Lmatrix}
    1 & 0 \\ [\TEMPORARYA] & 1
    \end{Lmatrix} = \begin{Lmatrix}
    -1 & 0 \\ 0 & -1
    \end{Lmatrix} \begin{Lmatrix}
        1 & 0 \\  -[\MU]-1/[\TEMPORARYA] & 1
        \end{Lmatrix}\begin{Lmatrix}
        [\TEMPORARYA] & 1 \\ 0 & 1/[\TEMPORARYA]
        \end{Lmatrix},} and the fact that \m{\begin{Smatrix}
        -1 & 0 \\ 0 & -1
        \end{Smatrix}} acts on \m{\ind_{\INZA}^G \TILDE\SIGMA_{r}} as multiplication by \m{(-1)^{\NEBENPOWEX}}. For the third equality we change the variable \m{\TEMPORARYA} to \m{\TEMPORARYAA=1/\TEMPORARYA}. For the fourth equality we change the variable \m{\MU} to \m{\LAMBDA = -\MU-\TEMPORARYAA} (we can do this if we treat  \m{\sum_{\TEMPORARYAA\in\FF_p^\times}} as the outer sum and  \m{\sum_{\mu\in\FF_p}} as the inner sum). %
        The fifth equality follows from the equations \algn{[-\LAMBDA-\TEMPORARYAA] = -[\LAMBDA]-[\TEMPORARYAA] +\BIGO{p}, \text{ and } \textstyle \sum_{j} \binom{u}{j}(-[\LAMBDA]-[\TEMPORARYAA])^{u-j} [\TEMPORARYAA]^{j} = (-[\LAMBDA])^u.}
        The sixth equality follows from the definition of \m{\ERRR}. %
 The eighth equality follows from the equation
  \Algn{\label{eq935gj4948ug985u} \textstyle [\TEMPORARYAA]^{{\NEBENPOWEX-s}} \begin{Lmatrix}
        1/[\TEMPORARYAA] & 1 \\ 0 & [\TEMPORARYAA]
        \end{Lmatrix}\GEN  & \textstyle = [\TEMPORARYAA]^{{-s}} \begin{Lmatrix}
        1 & 1 \\ 0 & [\TEMPORARYAA]
        \end{Lmatrix}\GEN  \textstyle = [\TEMPORARYA]^{s} \begin{Lmatrix}
        1 & [\TEMPORARYA] \\ 0 & 1
        \end{Lmatrix}\GEN. }
 Equations~\eqref{eqt43gj49jg4j59} and~\eqref{eq12rf0jwf98j3498fj3984jf3498j}  imply equation~\eqref{eqf34398jf983j498fj934}. Then
 \Algn{\label{eqf304jfk39fk3049k3f}\nonumber\ERRR & \textstyle\vphantom{(\ZETA-1) \sum_{\TEMPORARYA \in \LCLASS,\,\TEMPORARYAA\in\FF_p^\times} \sum_{j=1}^{p-1} 
        \frac{(-1)^j}{j} \TEMPORARYA^{u+j} [\TEMPORARYAA]^{\NEBENPOWEX-s-j} \begin{Lmatrix}
        1 & 0 \\  \TEMPORARYA & 1
        \end{Lmatrix}\begin{Lmatrix}
        1/[\TEMPORARYAA] & 1 \\ 0 & [\TEMPORARYAA]
        \end{Lmatrix}\sqbrackets{}{\INZB}{\QQp}{\GEN} + \BIGO{\QUARK^2}} \EQUALZEROIDEAL \sum_{\TEMPORARYAA\in\FF_p^\times,\,\LAMBDA \in \FF_p} { [\LAMBDA]^u [\TEMPORARYAA]^{{\NEBENPOWEX-s}} \begin{Lmatrix}
        1 & 0 \\  [\LAMBDA] & 1
        \end{Lmatrix}\left(\begin{Lmatrix}
        1 & 0 \\  [\LAMBDA+\TEMPORARYAA]-[\TEMPORARYAA]-[\LAMBDA] & 1
        \end{Lmatrix}-1\right)\begin{Lmatrix}
        1/[\TEMPORARYAA] & 1 \\ 0 & [\TEMPORARYAA]
        \end{Lmatrix}}\sqbrackets{}{\INZB}{\QQp}{\GEN} + \BIGO{p} \\ &\nonumber  \textstyle\vphantom{(\ZETA-1) \sum_{\TEMPORARYA \in \LCLASS,\,\TEMPORARYAA\in\FF_p^\times} \sum_{j=1}^{p-1} 
        \frac{(-1)^j}{j} \TEMPORARYA^{u+j} [\TEMPORARYAA]^{\NEBENPOWEX-s-j} \begin{Lmatrix}
        1 & 0 \\  \TEMPORARYA & 1
        \end{Lmatrix}\begin{Lmatrix}
        1/[\TEMPORARYAA] & 1 \\ 0 & [\TEMPORARYAA]
        \end{Lmatrix}\sqbrackets{}{\INZB}{\QQp}{\GEN} + \BIGO{\QUARK^2}} \EQUIVZEROIDEAL \sum_{\TEMPORARYAA\in\FF_p^\times,\,\LAMBDA \in \FF_p} { [\LAMBDA]^u [\TEMPORARYAA]^{{\NEBENPOWEX-s}} \begin{Lmatrix}
        1 & 0 \\  [\LAMBDA] & 1
        \end{Lmatrix}\left(\ZETA^{\DIFF(\LAMBDA/\TEMPORARYAA)}-1\right)\begin{Lmatrix}
        1/[\TEMPORARYAA] & 1 \\ 0 & [\TEMPORARYAA]
        \end{Lmatrix}}\sqbrackets{}{\INZB}{\QQp}{\GEN} + \BIGO{\mathrm{err}_0a^{-1}}
        \\ &\nonumber  \textstyle\vphantom{(\ZETA-1) \sum_{\TEMPORARYA \in \LCLASS,\,\TEMPORARYAA\in\FF_p^\times} \sum_{j=1}^{p-1} 
        \frac{(-1)^j}{j} \TEMPORARYA^{u+j} [\TEMPORARYAA]^{\NEBENPOWEX-s-j} \begin{Lmatrix}
        1 & 0 \\  \TEMPORARYA & 1
        \end{Lmatrix}\begin{Lmatrix}
        1/[\TEMPORARYAA] & 1 \\ 0 & [\TEMPORARYAA]
        \end{Lmatrix}\sqbrackets{}{\INZB}{\QQp}{\GEN} + \BIGO{\QUARK^2}} \EQUALZEROIDEAL (\ZETA-1) \sum_{\TEMPORARYAA\in\FF_p^\times,\,\LAMBDA \in \FF_p} { \DIFF(\LAMBDA/\TEMPORARYAA)  [\LAMBDA]^u [\TEMPORARYAA]^{{\NEBENPOWEX-s}} \begin{Lmatrix}
        1 & 0 \\  [\LAMBDA] & 1
        \end{Lmatrix}\begin{Lmatrix}
        1/[\TEMPORARYAA] & 1 \\ 0 & [\TEMPORARYAA]
        \end{Lmatrix}}\sqbrackets{}{\INZB}{\QQp}{\GEN} + \BIGO{\QUARK^2}
        \\ &  \textstyle\vphantom{(\ZETA-1) \sum_{\TEMPORARYA \in \LCLASS,\,\TEMPORARYAA\in\FF_p^\times} \sum_{j=1}^{p-1} 
        \frac{(-1)^j}{j} \TEMPORARYA^{u+j} [\TEMPORARYAA]^{\NEBENPOWEX-s-j} \begin{Lmatrix}
        1 & 0 \\  \TEMPORARYA & 1
        \end{Lmatrix}\begin{Lmatrix}
        1/[\TEMPORARYAA] & 1 \\ 0 & [\TEMPORARYAA]
        \end{Lmatrix}\sqbrackets{}{\INZB}{\QQp}{\GEN} + \BIGO{\QUARK^2}} \EQUALZEROIDEAL
        (\ZETA-1)\sum_{\TEMPORARYAA\in\FF_p^\times,\, \TEMPORARYA\in\LCLASSTEICH} \sum_{j=1}^{p-1} 
        \frac{(-1)^j}{j} \TEMPORARYA^{u+j} [\TEMPORARYAA]^{\NEBENPOWEX-s-j} \begin{Lmatrix}
        1 & 0 \\  \TEMPORARYA & 1
        \end{Lmatrix}\begin{Lmatrix}
        1/[\TEMPORARYAA] & 1 \\ 0 & [\TEMPORARYAA]
        \end{Lmatrix}\sqbrackets{}{\INZB}{\QQp}{\GEN} + \BIGO{\QUARK^2}.}
        The congruence follows from equation~\eqref{eq935gj4948ug985u} and lemma~\ref{LEMMAX6primeq}, and the second equality follows from \m{\BIGO{\mathrm{err}_0a^{-1}}\subseteq \BIGO{\QUARK^{2}}} and the fact that, for \m{t\in \ZZ_p},  \algn{\ZETA^t-1=t(\ZETA-1)+\BIGO{\QUARK^2}.}  The third equality follows from part~\eqref{comb10.2} of lemma~\ref{LEMMAX3.04}. Equation~\eqref{eq935gj4948ug985u} and lemma~\ref{LEMMAX6primeq} imply that if, for each \m{\TEMPORARYAA\in\FF_p^\times},   we replace the  inner sum ``\m{\sum_{\TEMPORARYA\in\LCLASSTEICH}}''  in the fourth line of equation~\eqref{eqf304jfk39fk3049k3f} with ``\m{\sum_{\TEMPORARYA\in\LCLASS_\TEMPORARYAA}}'' then that only changes the expression by a term that is in \m{\BIGO{\QUARK^2}} and therefore results in the desired formula for \m{\ERRR}. %
\end{proof}

\begin{lemma}\label{LEMMAX7}
{
Let \m{m\in\ZZ_{\GEQSLANT1}}. For \m{i\in\{1,\ldots,m\}}, let \m{(u_i,s_i)\in\ZZ_{\GEQSLANT0}\times\{1,\ldots,p-1\}} and \m{\lambda_i \in \FF_p}. Suppose  that
\algn{\sum_{i=1}^m\lambda_i (u_i,s_i) \in \VECTORSP} is in the kernel of \m{\Xi_{z,t}} for all \m{(z,t) \in \{0,\ldots,p-1\}\times\{p-\TEMPORARYI,\ldots,p-2\}}. Then
 \algn{
 & \textstyle
\sum_{i=1}^m [\lambda_i] \sum_{\TEMPORARYA \in \FF_p} {[\TEMPORARYA]^{u_i} \begin{Lmatrix}
    1 & 0 \\ [\TEMPORARYA] & 1
    \end{Lmatrix}}\sqbrackets{}{\INZB}{\QQp}{\GEN_{s_i}} \\& \null\qquad\textstyle + \sum_{i=1}^m (-1)^{u_i+\NEBENPOWEX}[\lambda_i]  \sum_{\beta \equiv s_i-\NEBENPOWEY \bmod{p-1}}{\binom{u_i}{\beta}\begin{Lmatrix}
0 & -1 \\ 1 & 0
\end{Lmatrix}}\sqbrackets{}{\INZB}{\QQp}{\GEN_{u_i-\beta}} %
 \\ & \null\qquad\qquad \textstyle \EQUIVZEROIDEAL \sum_{\TEMPORARYA \in \FF_p} {A_\TEMPORARYA \begin{Lmatrix}
    1 & 0 \\ [\TEMPORARYA] & 1
    \end{Lmatrix}}\sqbrackets{}{\INZB}{\QQp}{\GEN} + \BIGO{\QUARK},}
    where, for \m{\TEMPORARYA\in\FF_p},
    \algn{A_\TEMPORARYA & \textstyle \vphantom{\sum_{i=1}^m (-1)^{u_i+\NEBENPOWEX} [\lambda_i] \sum_{\MU \in \FF_p} \sum_{j\equiv s_i+\TEMPORARYI-\NEBENPOWEX\bmod{p-1}} \binom{u_i}{j}[\MU]^{u_i-j} [\TEMPORARYA]^{p-1-\TEMPORARYI} \begin{Lmatrix}
        0 & -1 \\ 1 & [\MU]
 \end{Lmatrix}}\EQUALZEROIDEAL \sum_{i=1}^m (-1)^{u_i+\NEBENPOWEX} [\lambda_i] \sum_{\MU \in \FF_p} \sum_{j\equiv s_i+\TEMPORARYI-\NEBENPOWEX\bmod{p-1}} \binom{u_i}{j}[\MU]^{u_i-j} [\TEMPORARYA]^{p-1-\TEMPORARYI} \begin{Lmatrix}
        0 & -1 \\ 1 & [\MU]
 \end{Lmatrix}
 \\&\textstyle{}\vphantom{\sum_{i=1}^m (-1)^{u_i+\NEBENPOWEX} [\lambda_i] \sum_{\MU \in \FF_p} \sum_{j\equiv s_i+\TEMPORARYI-\NEBENPOWEX\bmod{p-1}} \binom{u_i}{j}[\MU]^{u_i-j} [\TEMPORARYA]^{p-1-\TEMPORARYI} \begin{Lmatrix}
        0 & -1 \\ 1 & [\MU]
 \end{Lmatrix}}\EQUIVZEROIDEAL\BIGO{\QUARK^\TEMPORARYI a^{-1}}.}
} 
\end{lemma}

\begin{proof}{Proof. }
    For \m{i\in\{1,\ldots,m\}}, let us apply lemma~\ref{LEMMAX5} with \m{(u,s)=(u_i,s_i)}, and let \m{c_\TEMPORARYA(u_i,s_i)} be the corresponding constant
    \algn{c_\TEMPORARYA(u_i,s_i)\textstyle = (-1)^{u_i+\NEBENPOWEX} \sum_{\MU \in \FF_p} \sum_{j=0}^{u_i} \binom{u_i}{j}[\MU]^{u_i-j} [\TEMPORARYA]^{\LESS{-j+s_i-\NEBENPOWEY}} \begin{Lmatrix}
        0 & -1 \\ 1 & [\MU]
        \end{Lmatrix} \in \ZZ_p[G].}
These \m{m} applications of lemma~\ref{LEMMAX5} result in \m{m} equations. %
 The equation
        \Algn{\label{eq34oi3jiofj3oi4fj43f4gkjt5gnt}
 & \textstyle 
\sum_{i=1}^m [\lambda_i]\sum_{\TEMPORARYA \in \FF_p} {[\TEMPORARYA]^{u_i} \begin{Lmatrix}
    1 & 0 \\ [\TEMPORARYA] & 1
    \end{Lmatrix}}\sqbrackets{}{\INZB}{\QQp}{\GEN_{s_i}} \nonumber\\& \null\qquad\textstyle +\sum_{i=1}^m [\lambda_i] (-1)^{u_i+\NEBENPOWEX} \sum_{\beta \equiv s_i-\NEBENPOWEY \bmod{p-1}}{\binom{u_i}{\beta}\begin{Lmatrix}
0 & -1 \\ 1 & 0
\end{Lmatrix}}\sqbrackets{}{\INZB}{\QQp}{\GEN_{u_i-\beta}} \nonumber
 \\& \null\qquad\qquad\textstyle \EQUIVZEROIDEAL  \textstyle\sum_{\TEMPORARYA \in \FF_p} {\sum_{i=1}^m [\lambda_i] c_\TEMPORARYA(u_i,s_i) \begin{Lmatrix}
    1 & 0 \\ [\TEMPORARYA] & 1
    \end{Lmatrix}}\sqbrackets{}{\INZB}{\QQp}{\GEN} + \BIGO{\QUARK}}
       follows as a linear combination of these \m{m} equations. 
    The linear maps    \m{\Xi_{z,t}} satisfy 
    \algn{ c_\TEMPORARYA(u_i,s_i) = (-1)^\NEBENPOWEX \sum_{z=0}^{p-1}\sum_{t=0}^{p-2} [\Xi_{z,t}(u_i,s_i)] \sum_{\MU\in\FF_p}[\MU]^z [\TEMPORARYA]^t \begin{Smatrix}
    0 & -1 \\ 1 & [\MU]
    \end{Smatrix}+ \BIGO{p}.}
   Since \m{\sum_{i=1}^m\lambda_i (u_i,s_i) \in \Ker \Xi_{z,t}} for all \m{(z,t) \in \{0,\ldots,p-1\}\times\{p-\TEMPORARYI,\ldots,p-2\}},
    \algn{&\textstyle{} \sum_{i=1}^m [\lambda_i] c_\TEMPORARYA(u_i,s_i) \\ & \null\qquad\textstyle= (-1)^\NEBENPOWEX \sum_{z=0}^{p-1} \sum_{t=0}^{p-1-\TEMPORARYI} \sum_{i=1}^m[\lambda_i\Xi_{z,t}(u_i,s_i)] \sum_{\MU\in\FF_p}[\MU]^z [\TEMPORARYA]^t \begin{Smatrix}
    0 & -1 \\ 1 & [\MU]
    \end{Smatrix}+ \BIGO{p}.}
    For \m{t\in\{0,\ldots,p-2\}}, lemma~\ref{LEMMAX6} applied to \m{\alpha=p-1-t} and the set of lifts \m{\LCLASSTEICH} implies that
\algn{\textstyle \sum_{\TEMPORARYA \in \FF_p} {[\TEMPORARYA]^t \begin{Lmatrix}
    1 & 0 \\ [\TEMPORARYA] & 1
    \end{Lmatrix}}\sqbrackets{}{\INZB}{\QQp}{\GEN} \EQUIVZEROIDEAL \BIGO{\QUARK^{p-1-t} a^{-1}} + \BIGO{\mathrm{err}_0a^{-2}}.} In particular, for \m{t \in \{0,\ldots,p-2-\TEMPORARYI\}}, the right side is in \m{\BIGO{\QUARK^{\TEMPORARYI+1} a^{-1}}}. Therefore
     \Algn{\label{eq985jug4jr398r9834jr983uj}\nonumber & \textstyle\sum_{\TEMPORARYA \in \FF_p} {\sum_{i=1}^m [\lambda_i] c_\TEMPORARYA(u_i,s_i) \begin{Lmatrix}
    1 & 0 \\ [\TEMPORARYA] & 1
    \end{Lmatrix}}\sqbrackets{}{\INZB}{\QQp}{\GEN} \\ &\nonumber \null\qquad\textstyle{}
    \EQUALZEROIDEAL  {  \sum_{t=0}^{p-1-\TEMPORARYI} B_t \sum_{\TEMPORARYA \in \FF_p} [\TEMPORARYA]^t   \begin{Lmatrix}
    1 & 0 \\ [\TEMPORARYA] & 1
    \end{Lmatrix}}\sqbrackets{}{\INZB}{\QQp}{\GEN} + \BIGO{p}
    \\ & \null\qquad\textstyle{}
    \EQUIVZEROIDEAL  {   B_{p-1-\TEMPORARYI} \sum_{\TEMPORARYA \in \FF_p} [\TEMPORARYA]^{p-1-\TEMPORARYI}  \begin{Lmatrix}
    1 & 0 \\ [\TEMPORARYA] & 1
    \end{Lmatrix}}\sqbrackets{}{\INZB}{\QQp}{\GEN} + \BIGO{\QUARK^{\TEMPORARYI+1} a^{-1}},}
    where, for \m{t\in\{0,\ldots,p-1-\TEMPORARYI\}}, \algn{B_{t} = (-1)^\NEBENPOWEX \sum_{z=0}^{p-1} \sum_{i=1}^m [\lambda_i\Xi_{z,t}(u_i,s_i)] \sum_{\MU\in\FF_p}[\MU]^z  \begin{Smatrix}
    0 & -1 \\ 1 & [\MU]
    \end{Smatrix}. }
    Since
    \algn{\sum_{\TEMPORARYA \in \FF_p} {A_\TEMPORARYA \begin{Lmatrix}
    1 & 0 \\ [\TEMPORARYA] & 1
    \end{Lmatrix}}\sqbrackets{}{\INZB}{\QQp}{\GEN} = {   B_{p-1-\TEMPORARYI} \sum_{\TEMPORARYA \in \FF_p} [\TEMPORARYA]^{p-1-\TEMPORARYI}  \begin{Lmatrix}
    1 & 0 \\ [\TEMPORARYA] & 1
    \end{Lmatrix}}\sqbrackets{}{\INZB}{\QQp}{\GEN} \EQUIVZEROIDEAL\BIGO{\QUARK^\TEMPORARYI a^{-1}},}
    equations~\eqref{eq34oi3jiofj3oi4fj43f4gkjt5gnt} and~\eqref{eq985jug4jr398r9834jr983uj} complete the proof.
\end{proof}    

\begin{lemma}\label{LEMMAX9}
{
Let us consider the quadruple
\Algn{\label{eq34f3oijgi3ojg4igj4iughi4ugh4g}\left(m,\DELTAAA,(u_i)_{1\LEQSLANT i \LEQSLANT m},(\LAMBDA_i)_{1\LEQSLANT i\LEQSLANT m}\right)}
with
\m{m\in\ZZ_{\GEQSLANT1}}, \m{\DELTAAA \in \ZZ}, and, for \m{i\in\{1,\ldots,m\}}, \m{u_i\in\{0,\ldots,p-1\}} and \m{\LAMBDA_i \in \FF_p}. %
If
\Algn{\label{eqg3oigj98j45g9j58ghufieofji3u5}\textstyle ((-1)^{u_1}\lambda_1,\ldots,(-1)^{u_m}\lambda_m) \, \MMM{m}{(u_i)_{1\LEQSLANT i\LEQSLANT m}}{\DELTAAA} = 0}
then
 \algn{
 & \textstyle
\sum_{i=1}^m [\lambda_i] \sum_{\TEMPORARYA \in \FF_p} {[\TEMPORARYA]^{u_i} \begin{Lmatrix}
    1 & 0 \\ [\TEMPORARYA] & 1
    \end{Lmatrix}}\sqbrackets{}{\INZB}{\QQp}{\GEN_{\MORE{u_i+\NEBENPOWEX-\DELTAAA}}} \\& \null\qquad\textstyle + \sum_{i=1}^m (-1)^{u_i+\NEBENPOWEX}[\lambda_i]  \sum_{w \equiv \DELTAAA \bmod{p-1}}{\binom{u_i}{w}\begin{Lmatrix}
0 & -1 \\ 1 & 0
\end{Lmatrix}}\sqbrackets{}{\INZB}{\QQp}{\GEN_{w}} %
  \\ & \null\qquad\qquad \textstyle \EQUIVZEROIDEAL \sum_{\TEMPORARYA \in \FF_p} {A_\TEMPORARYA \begin{Lmatrix}
    1 & 0 \\ [\TEMPORARYA] & 1
    \end{Lmatrix}}\sqbrackets{}{\INZB}{\QQp}{\GEN} + \BIGO{\QUARK},}
    where, for \m{\TEMPORARYA\in\FF_p},
    \algn{A_\TEMPORARYA  & \textstyle\vphantom{\sum_{i=1}^m (-1)^{u_i+\NEBENPOWEX} [\lambda_i] \sum_{\MU \in \FF_p} \sum_{j\equiv \delta-\TEMPORARYI\bmod{p-1}} \binom{u_i}{j}[\MU]^{j} [\TEMPORARYA]^{p-1-\TEMPORARYI} \begin{Lmatrix}
        0 & -1 \\ 1 & [\MU]
 \end{Lmatrix} }\EQUALZEROIDEAL \sum_{i=1}^m (-1)^{u_i+\NEBENPOWEX} [\lambda_i] \sum_{\MU \in \FF_p} \sum_{j\equiv \delta-\TEMPORARYI\bmod{p-1}} \binom{u_i}{j}[\MU]^{j} [\TEMPORARYA]^{p-1-\TEMPORARYI} \begin{Lmatrix}
        0 & -1 \\ 1 & [\MU]
 \end{Lmatrix}  \\ & \textstyle
 \vphantom{\sum_{i=1}^m (-1)^{u_i+\NEBENPOWEX} [\lambda_i] \sum_{\MU \in \FF_p} \sum_{j\equiv \delta-\TEMPORARYI\bmod{p-1}} \binom{u_i}{j}[\MU]^{j} [\TEMPORARYA]^{p-1-\TEMPORARYI} \begin{Lmatrix}
        0 & -1 \\ 1 & [\MU]
 \end{Lmatrix} }\EQUIVZEROIDEAL\BIGO{\QUARK^\TEMPORARYI a^{-1}}.
}} 
\end{lemma}

\begin{proof}{Proof. } Let us define \m{s_i=\MORE{u_i+\NEBENPOWEX-\DELTAAA}}.
 If   \m{(z,t) \in \{0,\ldots,p-1\}\times\{p-\TEMPORARYI,\ldots,p-2\}} and \m{i\in\{1,\ldots,m\}}, and if \m{\Xi_{z,t}(u_i,s_i)\ne 0}, %
  then from the definition of \m{\Xi_{z,t}} we have
\algn{\textstyle z-t \equiv u_i-s_i+\NEBENPOWEX \equiv \DELTAAA \bmod {p-1}, \text{ and }  \textstyle \Xi_{z,t}(u_i,s_i) = (-1)^{u_i}\binom{u_i}{z}.}
Since \m{t=p-1-w} for some \m{w\in\{1,\ldots,\TEMPORARYI-1\}}, we have \m{z\equiv \DELTAAA-w\bmod{p-1}}, so \m{z\in S_\DELTAAA}. %
 As equation~\eqref{eqg3oigj98j45g9j58ghufieofji3u5} implies that  \algn{\sum_{i=1}^m (-1)^{u_i}{\LAMBDA}_i \binom{u_i}{j}=0} for all \m{j\in S_\DELTAAA}, it follows that \m{{\Xi_{z,t}(u_i,s_i)= 0}}, which is a contradiction. This means that the parameters   \m{(u_i,s_i,\LAMBDA_i)_{1\LEQSLANT i\LEQSLANT m}} satisfy the conditions of lemma~\ref{LEMMAX7}, so the desired result follows  from lemma~\ref{LEMMAX7}.
\end{proof}

\NEWPAGE\section[Proof of theorem~\ref{MAINTHM2}]{Proof of theorem~\ref{MAINTHM1}}\label{SECPROOF}

 Lemma~\ref{LEMMAX3.043}
gives a list of subquotients \algn{\left\{\ind_{KZ}^G I_{i,j} \,|\, i\in \{0,\ldots,p^{n-1}-1\}, \,j\in\{1,\ldots,p-1\}\right\}} of \m{\ind_{\INZA}^{G} \QUOTIENT}, and explicitly describes their elements. %
 If \m{i\in\{0,\ldots,p^{n-1}-2\}} and \m{i_{1},\ldots,i_{n-1} \in \{0,\ldots,p-1\}} %
are such that \m{i=i_1+\cdots+i_{n-1}p^{n-2}}, then there is an index \m{w\in\{1,\ldots,n-1\}} with \m{i_w<p-1}.
Lemma~\ref{LEMMAX6primep} and the fact that \algn{\CONST{l}{1}_{-\TEMPORARYA}=\BIGO{\QUARK}} for all \m{\TEMPORARYA\in\ZZ_p} and all \m{l\in\{1,\ldots,p-2\}}  (due to part~\eqref{gammaprop3} of lemma~\ref{LEMMAX3.02}) imply that
\algn{\textstyle \sum_{\TEMPORARYA_w \in \FF_p} [\TEMPORARYA_w]^{i_w} \begin{Lmatrix}
1 & 0 \\ [\TEMPORARYA_w]p^w & 1 
\end{Lmatrix} \sqbrackets{}{\INZB}{\QQp}{\GEN_j} \EQUIVZEROIDEAL   \sum_{\TEMPORARYA_w\in \FF_p} [\TEMPORARYA_w]^{i_w}  \sqbrackets{}{\INZB}{\QQp}{\GEN_j} + \BIGO{\QUARK} = \BIGO{\QUARK}} for all \m{j\in\{1,\ldots,p-1\}}.
Therefore the reduction modulo \m{\maxideal} of this element is in \m{\INAAA  }, implying by equation~\eqref{eq3t094u43u984398} that  each element of \m{\ind_{KZ}^G I_{i,j}} is represented by an element of \m{\INAAA  }, and therefore also by an element of the kernel of the surjective homomorphism 
\m{\textstyle \QUOTIENTMAP : \ind_{\INZA}^G \QUOTIENT \xtwoheadrightarrow{} \overline\Theta_{k,a,\varepsilon}} given by lemma~\ref{LEMMAX3}. Therefore any factor of \m{\overline\Theta_{k,a,\varepsilon}} must be a factor of a representation in the set \Algn{\label{eq3k40og30ok3094kf3o0gfjgi5g45}\left\{\ind_{KZ}^G I_{p^{n-1}-1,j} \,|\, j\in\{1,\ldots,p-1\}\right\}.} 
For \m{j\in\{1,\ldots,p-1\}}, let 
\algn{\SUBMODULE(j) & \textstyle = \sigma_{\MORE{2j-\NEBENPOWEY}}(\NEBENPOWEX-s)\subset I_{p^{n-1}-1,j},\\
\QUOTIENTI(j) & \textstyle = \sigma_{\LESS{\NEBENPOWEX-2j}}(j)\cong I_{p^{n-1}-1,j}/\sigma_{\MORE{2j-\NEBENPOWEY}}(\NEBENPOWEX-j).}These two irreducible representations are the only factors of \m{I_{p^{n-1}-1,j}}, as shown in section~\ref{SECnot}. For \m{j\in \{1,\ldots,p-1\}},
equation~\eqref{eq3t094u43u984398} implies that   
\algn{\nonumber
       & \textstyle  \sum_{\TEMPORARYA ,\TEMPORARYA_1,\ldots,\TEMPORARYA_{n-1} \in \FF_p} \TEMPORARYA_1^{i_1}\cdots \TEMPORARYA_{n-1}^{i_{n-1}}
    f(-\TEMPORARYA,1)\begin{Lmatrix}
    1 & 0 \\ [\TEMPORARYA] + [\TEMPORARYA_1]p+\cdots+[\TEMPORARYA_{n-1}]p^{n-1} & 1
    \end{Lmatrix}  \sqbrackets{}{\INZB}{\FFp}{\GENN_j} \\ & \null\qquad + \textstyle \sum_{\TEMPORARYA_1,\ldots,\TEMPORARYA_{n-1} \in \FF_p} \TEMPORARYA_1^{i_1}\cdots \TEMPORARYA_{n-1}^{i_{n-1}}
    f(-1,0)\begin{Lmatrix}
 -[\TEMPORARYA_1]p-\cdots-[\TEMPORARYA_{n-1}]p^{n-1} & -1 \\ 1 & 0
\end{Lmatrix}  \sqbrackets{}{\INZB}{\FFp}{\GENN_j}
    } represents \m{1\sqbracketss{}{KZ}{\FFp}{f} \in \ind_{KZ}^G I_{p^{n-1}-1,j}} when \m{i_1=\cdots=i_{n-1}=p-1}. Moreover, when \m{i_w<p-1} for some index \m{w\in\{1,\ldots,n-1\}}, it represents \m{1\sqbracketss{}{KZ}{\FFp}{f} \in \ind_{KZ}^G I_{i,j}} for some \m{i<p^{n-1}-1}, and therefore it represents \m{0 \in \ind_{KZ}^G I_{p^{n-1}-1,j}}.
 So if the degree of each \m{z_j} in \m{h \in \FFp[z_1,\ldots,z_{n-1}]} is exactly \m{p-1} and the coefficient of \m{z_1^{p-1}\cdots z_{n-1}^{p-1}} is \m{1}, then 
 \algn{\nonumber
       & \textstyle  \sum_{\TEMPORARYA ,\TEMPORARYA_1,\ldots,\TEMPORARYA_{n-1} \in \FF_p} h(\TEMPORARYA_1,\ldots, \TEMPORARYA_{n-1})
    f(-\TEMPORARYA,1)\begin{Lmatrix}
    1 & 0 \\ [\TEMPORARYA] + [\TEMPORARYA_1]p+\cdots+[\TEMPORARYA_{n-1}]p^{n-1} & 1
    \end{Lmatrix}  \sqbrackets{}{\INZB}{\FFp}{\GENN_j} \\ & \null\qquad + \textstyle \sum_{\TEMPORARYA_1,\ldots,\TEMPORARYA_{n-1} \in \FF_p} h(\TEMPORARYA_1,\ldots, \TEMPORARYA_{n-1})
    f(-1,0)\begin{Lmatrix}
 -[\TEMPORARYA_1]p-\cdots-[\TEMPORARYA_{n-1}]p^{n-1} & -1 \\ 1 & 0
\end{Lmatrix}  \sqbrackets{}{\INZB}{\FFp}{\GENN_j}
    } represents \m{1\sqbracketss{}{KZ}{\FFp}{f} \in \ind_{KZ}^G I_{p^{n-1}-1,j}}.
 For \algn{h(z_1,\ldots,z_{n-1}) = \prod_{j=1}^{n-1}(z_j^{p-1}-1),} this implies that
 \Algn{\label{eq394j9jh39fj93jg4398gj49jg498jg}
       & \textstyle  \sum_{\TEMPORARYA  \in \FF_p} 
    f(-\TEMPORARYA,1)\begin{Lmatrix}
    1 & 0 \\ [\TEMPORARYA] & 1
    \end{Lmatrix}  \sqbrackets{}{\INZB}{\FFp}{\GENN_j} + \textstyle 
    f(-1,0)\begin{Lmatrix}
0 & -1 \\ 1 & 0
\end{Lmatrix}  \sqbrackets{}{\INZB}{\FFp}{\GENN_j}
    } represents \m{(-1)^{n-1}\sqbracketss{}{KZ}{\FFp}{f} \in \ind_{KZ}^G I_{p^{n-1}-1,j}}. Together with equation~\eqref{EQUATION1}, this gives an explicit description of \m{\ind_{KZ}^G\SUBMODULE(j)} and \m{\ind_{KZ}^G\QUOTIENTI(j)} for \m{j \in \{1,\ldots,p-1\}}.
 
\subsection{Theorem~\ref{MAINTHM1} follows from propositions~\ref{PROP1} \& \ref{PROP2}} In this section we show that theorem~\ref{MAINTHM1} follows from two propositions about the kernel of \m{\QUOTIENTMAP}.

\begin{proposition}\label{PROP1}
 The kernel of \m{\QUOTIENTMAP} contains a representative of a generator for each  representation in the set
 \algn{&\textstyle\left\{\ind_{KZ}^G I_{p^{n-1}-1,l} \,|\, l \in \{1,\ldots,\TEMPORARYI-1\}\right\} \cup \left\{\ind_{KZ}^G \SUBMODULE(l) \,|\, l \in \{\TEMPORARYI,\ldots,p-1\}\right\} \\ &\textstyle \null\,\,\,\quad \cup \left\{\ind_{KZ}^{G} \QUOTIENTI(l) \,|\, l \in \{1,\ldots,p-1\}\backslash \{\MORE{\NEBENPOWEX-w}\,|\,w\in \{0,\ldots,\TEMPORARYI-1\}\}\right\}.}
\end{proposition}

\begin{proposition}\label{PROP2}
 If \m{\MORE{\NEBENPOWEX}=2\TEMPORARYI-1} and \m{\varrho\in\ind_{KZ}^G\QUOTIENTI(\TEMPORARYI)}, then the kernel of \m{\QUOTIENTMAP} contains a representative of
 \algn{\left\{\begin{array}{rl} (T_{\sigma}^2 - \phi\, T_{\sigma}+1)(\varrho)%
    & \text{if }  a\in\CIRCLES_{\TEMPORARYI} \text{,} \\[3pt] T_{\sigma}(\varrho)  &  \text{if }  a\not\in\CIRCLES_{\TEMPORARYI} \text{,} \end{array}\right.} where \m{T_\sigma} is the Hecke operator in \m{\End_G(\ind_{KZ}^G\QUOTIENTI(\TEMPORARYI))} corresponding to the double coset of \m{\begin{Smatrix} p&0 \\ 0 &1\end{Smatrix}}, and \m{\phi} is the reduction modulo \m{\maxideal} of 
\algn{\textstyle  ((\TEMPORARYI-1)!)^{-1} (1-\ZETA)^{\TEMPORARYI-1} a^{-3} \left( \frac{1}{\TEMPORARYI!(\TEMPORARYI-1)!} (\ZETA-1)^{2\TEMPORARYI-1}-a^2\right) .}
\end{proposition}
 
 \begin{proof}{Proof that propositions~\ref{PROP1} \& \ref{PROP2} imply theorem~\ref{MAINTHM1}. } Theorem~\ref{MAINTHM1} is equivalent to theorem~\ref{MAINTHM2}, so it is enough to show that propositions~\ref{PROP1} and~\ref{PROP2} imply theorem~\ref{MAINTHM2}. Equation~\eqref{eq3k40og30ok3094kf3o0gfjgi5g45} and proposition~\ref{PROP1} imply  that  any  factor of \m{\overline\Theta_{k,a,\varepsilon}} must be a factor of a representation in the set
 \Algn{\label{eqf34kt309jg49jfoifjk3fj349ij}\textstyle \left\{\ind_{KZ}^G\QUOTIENTI(\MORE{\NEBENPOWEX-l}) \,|\, l\in\{0,\ldots,\TEMPORARYI-1\}\right\} \backslash \left\{\ind_{KZ}^G\QUOTIENTI(l) \,|\, l\in\{1,\ldots,\TEMPORARYI-1\}\right\}.} 
 Theorem~\ref{BERGERX1} and theorem~2.7.1 in~\cite{Bre03a} (which classifies  \m{\repnn(t,\lambda,\psi)}) show that 
 either \m{\overline\Theta_{k,a,\varepsilon}} is irreducible (and infinite-dimensional), or it has exactly two infinite-dimensional factors. If \m{\overline\Theta_{k,a,\varepsilon}} is irreducible and a factor of
 some \m{
 \ind_{KZ}^G\sigma_b(d)} that belongs to the set in equation~\eqref{eqf34kt309jg49jfoifjk3fj349ij}, then theorem~\ref{BERGERX1} and theorem~2.7.1 in~\cite{Bre03a} imply that it is \m{\repnn(b,0,\omega^{d})}, so it belongs to the set 
   \algn{   \left\{\BIRREDUCIBLE{k}{\NEBENPOWER}{\NEBENPOWEX-l} %
    \,|\, l \in \{0,\ldots,\TEMPORARYI-1\}\right\}\backslash \left\{\BIRREDUCIBLE{k}{\NEBENPOWER}{l} 
    \,|\, l \in \{1,\ldots,\TEMPORARYI-1\}\right\}.}
 If \m{\overline\Theta_{k,a,\varepsilon}} is reducible and has two infinite-dimensional factors \m{R_1} and \m{R_2}, and if, for \m{i\in \{1,2\}}, \m{R_i} is a  factor of
 some \m{\ind_{KZ}^G\QUOTIENTI(\MORE{\NEBENPOWEX-b_i}) \cong
 \ind_{KZ}^G\sigma_{{2b_i-\NEBENPOWEY}}(\NEBENPOWEX-b_i)}  that belongs to the set in equation~\eqref{eqf34kt309jg49jfoifjk3fj349ij}, then theorem~\ref{BERGERX1} and theorem~2.7.1 in~\cite{Bre03a} imply that
 \algn{%
  \NEBENPOWEX-b_2 \equiv (\NEBENPOWEX-b_1)+(2b_1-\NEBENPOWEY+1)\bmod{p-1} \Longleftrightarrow b_1+b_2\equiv \NEBENPOWEX-1\bmod{p-1}.} Since \m{b_1,b_2\in\{0,\ldots,\TEMPORARYI-1\}} and \m{\MORE{\NEBENPOWEX-b_i}\not\in \{1,\ldots,\TEMPORARYI\} } for \m{i\in\{1,2\}},  
it follows that
\algn{b_1 = \MORE{\NEBENPOWEX-b_2}-1\not\in\{0,\ldots,\TEMPORARYI-2\}, \text{ and }b_2 = \MORE{\NEBENPOWEX-b_1}-1\not\in\{0,\ldots,\TEMPORARYI-2\}.} Therefore \m{b_1=b_2=\TEMPORARYI-1} and \m{\MORE{\NEBENPOWEX}=2\TEMPORARYI-1}. So if \m{\overline\Theta_{k,a,\varepsilon}} is reducible then \m{\MORE{\NEBENPOWEX}=2\TEMPORARYI-1} and both of its infinite-dimensional factors are factors of \m{\QUOTIENTI(\TEMPORARYI)}. To summarize, either \m{\overline\Theta_{k,a,\varepsilon}}  belongs to the set  
   \algn{   \left\{\BIRREDUCIBLE{k}{\NEBENPOWER}{\NEBENPOWEX-l} %
    \,|\, l \in \{1,\ldots,\TEMPORARYI-1\}\right\}\backslash \left\{\BIRREDUCIBLE{k}{\NEBENPOWER}{l} 
    \,|\, l \in \{1,\ldots,\TEMPORARYI-1\}\right\},} or \m{\MORE{\NEBENPOWEX}=2\TEMPORARYI-1} and all infinite-dimensional factors of \m{\overline\Theta_{k,a,\varepsilon}} are factors of \m{\QUOTIENTI(\TEMPORARYI)}. In the latter case, theorem~\ref{BERGERX1} and proposition~\ref{PROP2} imply that \m{\overline\Theta_{k,a,\varepsilon}\cong \BPI}.
 \end{proof}
 
 Therefore our goal is to prove propositions~\ref{PROP1} and~\ref{PROP2}.
 
 \section{Proof of proposition~\ref{PROP1}}
 \label{secprop1}
 
 Equation~\eqref{eq394j9jh39fj93jg4398gj49jg498jg} implies that, for \m{w\in\{1,\ldots,p-1\}},
 \algn{
\textstyle
     \mathrm{gen}_1(w) = 1\sqbrackets{}{\INZB}{\FFp}{\GENN_w} & \textstyle \text{ represents a generator of } I_{p^{n-1}-1,w}, \\ \textstyle
   \mathrm{gen}_2(w)=  \sum_{\TEMPORARYA \in \FF_p} { \begin{Lmatrix}
    1 & 0 \\ [\TEMPORARYA] & 1
    \end{Lmatrix}}\sqbrackets{}{\INZB}{\FFp}{\GENN_w} & \textstyle \text{ represents a generator of } \SUBMODULE(w).
 }
 Let us prove the following three claims.

\emph{Claim 1. } If \m{w \in \{1,\ldots,\TEMPORARYI-1\}}, then \m{\mathrm{gen}_1(w) \in \INAAA}. %

     \emph{Claim 2. } If \m{w \in \{\TEMPORARYI,\ldots,p-1\}} and \m{\MORE{\NEBENPOWEX-w}\in \{\TEMPORARYI,\ldots,p-1\}}, then \m{\mathrm{gen}_2(w) \in \INAAA}. %
     
    \emph{Claim 
     3. } If \m{w \in \{\TEMPORARYI,\ldots,p-1\}} and \m{\MORE{\NEBENPOWEX-w}\in \{1,\ldots,\TEMPORARYI-1\}}, then  \algn{\mathrm{gen}_2(w) + \sum_{i=0}^{\TEMPORARYI-1}b_i\mathrm{gen}_1(i) \in \INAAA,} for some \m{b_0,\ldots,b_{\TEMPORARYI-1}\in\FFp[G]}.

 \emph{Proof of claim 1. } 
If we apply lemma~\ref{LEMMAX6} with \m{\alpha=w\in\{1,\ldots,p-1\}} and the set of lifts \m{\LCLASSTEICH}, and if we multiply both sides of the resulting equation by \m{a\QUARK^{-w}}, we get %
 \algn{\textstyle{1}\sqbrackets{}{\INZB}{\QQp}{\GEN_w} \EQUIVZEROIDEAL \BIGO{a\QUARK^{-w}}=  \BIGO{a\QUARK^{1-\TEMPORARYI}}.} The reduction modulo \m{\maxideal} of this equation implies that %
  \m{\mathrm{gen}_1(w) \in \INAAA}. \hspace{\stretch{1}} \m{\square}%
 
 \emph{Proof of claim 2. } 
Let %
 \m{l = \MORE{\NEBENPOWEX-w}}.
Since \m{0\not\in S_{l}} and therefore \m{\binom{0}{j}=0} for all \m{j\in S_{l}}, %
 the quadruple
\algn{(m,\DELTAAA,(u_i)_{1\LEQSLANT i \LEQSLANT m},(\LAMBDA_i)_{1\LEQSLANT i\LEQSLANT m}) = (1,l,(0),(1))} 
satisfies equation~\eqref{eqg3oigj98j45g9j58ghufieofji3u5}. %
  Lemma~\ref{LEMMAX9} then implies that %
\algn{
& \textstyle
\sum_{\TEMPORARYA \in \FF_p} { \begin{Lmatrix}
    1 & 0 \\ [\TEMPORARYA] & 1
    \end{Lmatrix}}\sqbrackets{}{\INZB}{\QQp}{\GEN_w} %
    + (-1)^{\NEBENPOWEX}\VERIFYIF{w\equiv\NEBENPOWEX \bmod{p-1}} \begin{Lmatrix}
0 & -1 \\ 1 & 0
\end{Lmatrix}\sqbrackets{}{\INZB}{\QQp}{\GEN_{0}} \EQUIVZEROIDEAL \BIGO{\QUARK^\TEMPORARYI a^{-1}}.} Lemma~\ref{LEMMAX3} implies that  \m{1\sqbrackets{}{\INZB}{\FFp}{\GENN_0} \in\INAAA},
so  \m{\INAAA} also contains \m{\mathrm{gen}_2(w)}. %
\hspace{\stretch{1}} \m{\square}
    
    \emph{Proof of claim 3. } 
    Let  \m{l = \MORE{\NEBENPOWEX-w}}.
  Then
  \algn{S_l = \{0,\ldots,l-1\} \cup \{p-\TEMPORARYI+l,\ldots,p-1\}.}
   The columns of
\algn{M = \MMM{l+1}{(0,p-w,\ldots,p-w+l-1)}{l}}
   indexed by 
  \m{\textstyle j\in \{p-\TEMPORARYI+l,\ldots,p-1\}} are zero because \algn{&\textstyle u_i < p-w+l\LEQSLANT p-\TEMPORARYI+l\LEQSLANT j \\&\null\qquad\textstyle\Longrightarrow \binom{u_i}{j}=0 \text{ for all } {i\in \{1,\ldots,l+1\}} \text{ and all } {j\in \{p-\TEMPORARYI+l,\ldots,p-1\}}.} 
 The \m{l\times l} submatrix of \m{M} consisting of the rows indexed by \m{i \in \{2,\ldots,l+1\}} and columns indexed by \m{j \in \{0,\ldots,l-1\}} is obtained from
  \algn{\textstyle \MATRIX_{l-1}(p-w,0) =\left(\binom{p-2-w+i}{j}\right)_{2\LEQSLANT i \LEQSLANT l+1,\, 0\LEQSLANT j \LEQSLANT l-1}} by reducing modulo \m{\maxideal}, so it is invertible by part~\eqref{vander2} of lemma~\ref{LEMMAX8}. This implies that the first row of \m{M} is a linear combination of the other rows of \m{M}, so there exist \m{\lambda_2,\ldots,\lambda_{l+1}\in\FF_p} such that the quadruple
\algn{%
(l+1,l,(0,p-w,\ldots,p-w+l-1),(1,\lambda_2,\ldots,\lambda_{l+1}))} satisfies equation~\eqref{eqg3oigj98j45g9j58ghufieofji3u5}. %
  Lemma~\ref{LEMMAX9} then implies that
 \algn{\textstyle \sum_{\TEMPORARYA\in\FF_p} {\begin{Lmatrix}
 1 & 0 \\ [\TEMPORARYA] & 1
 \end{Lmatrix}}\sqbrackets{}{\INZB}{\QQp}{\GEN_{w}} + \sum_{i=0}^{l} b_i \sqbrackets{}{\INZB}{\QQp}{\GEN_{i}} \EQUIVZEROIDEAL \BIGO{\QUARK^\TEMPORARYI a^{-1}},}
 for some \m{b_0,\ldots,b_l \in \ZZp[G]}. We get the desired claim after reducing modulo~\m{\maxideal}. \hspace{\stretch{1}} \m{\square}
The three claims we showed imply that the kernel of \m{\QUOTIENTMAP} contains a representative of a generator for each  representation in the set
 \algn{\left\{\ind_{KZ}^G I_{p^{n-1}-1,l} \,|\, l \in \{1,\ldots,\TEMPORARYI-1\}\right\} \cup \left\{\ind_{KZ}^G \SUBMODULE(l) \,|\, l \in \{\TEMPORARYI,\ldots,p-1\}\right\}.} In order to find in the kernel of \m{\QUOTIENTMAP} a representative of a generator of \algn{\ind_{KZ}^{G} \QUOTIENTI(w)\text{  for }w \in \{1,\ldots,p-1\}\backslash \{\MORE{\NEBENPOWEX-l}\,|\,l\in \{0,\ldots,\TEMPORARYI-1\}\},} we 
   use a different construction in each of the following five cases.

\begin{center}
\begin{tabular}{r|l}
    case 1 & \m{\TEMPORARYI\LEQSLANT\MORE{\NEBENPOWEX}} and \m{w \in (\MORE{\NEBENPOWEX},p-1]}. \NEWLINEA
    case 2 & \m{\TEMPORARYI\LEQSLANT\MORE{\NEBENPOWEX}} and \m{w \in [\TEMPORARYI+1,\MORE{\NEBENPOWEX}-\TEMPORARYI]}. \NEWLINEA
    case 3 & \m{\TEMPORARYI\LEQSLANT\MORE{\NEBENPOWEX}} and \m{w=\TEMPORARYI}.\NEWLINEA
    case 4 & \m{\MORE{\NEBENPOWEX} < \TEMPORARYI} and \m{w \in [\TEMPORARYI+1,p-\TEMPORARYI+\MORE{\NEBENPOWEX})}.  \NEWLINEA
    case 5 & \m{\MORE{\NEBENPOWEX} < \TEMPORARYI} and \m{w=\TEMPORARYI}. \NEWLINEA
\end{tabular}
\end{center}

\emph{Case 1. }\label{case1}
      Let \m{\delta=w}. Then
  \m{S_\delta = \{w-\TEMPORARYI+1,\ldots,w-1\}},
     so
\algn{M = \MMM{\TEMPORARYI}{(w-i+1)_{1\LEQSLANT i \LEQSLANT \TEMPORARYI}}{\delta} = \left(\binom{w-i+1}{w-j}\right)_{1\LEQSLANT i\LEQSLANT\TEMPORARYI,\, 1\LEQSLANT j \LEQSLANT\TEMPORARYI-1}.} The square submatrix obtained from \m{M} by discarding the first row is
     upper triangular with \m{1}'s on the diagonal and is therefore invertible. This means that the first row of \m{M} is a linear combination of the other rows of \m{M}, so there exist \m{\lambda_2,\ldots,\lambda_{\TEMPORARYI}\in\FF_p} such that
the quadruple
\algn{%
(\TEMPORARYI,\delta,(w-i+1)_{1\LEQSLANT i \LEQSLANT \TEMPORARYI},(1,\lambda_2,\ldots,\lambda_{\TEMPORARYI}))} satisfies equation~\eqref{eqg3oigj98j45g9j58ghufieofji3u5}. %
  Lemma~\ref{LEMMAX9} then implies that
 \Algn{\label{eqf43j93jr9348}\textstyle (-1)^{w+\NEBENPOWEX} {\begin{Lmatrix}
 0 & -1 \\ 1 & 0
 \end{Lmatrix}}\sqbrackets{}{\INZB}{\QQp}{\GEN_{w}} + \sum_{i=\MORE{\NEBENPOWEX}-\TEMPORARYI+1}^{\MORE{\NEBENPOWEX}} b_i \sqbrackets{}{\INZB}{\QQp}{\GEN_{i}} \EQUIVZEROIDEAL \BIGO{\QUARK^\TEMPORARYI a^{-1}},}
  for some \m{b_{\MORE{\NEBENPOWEX}-\TEMPORARYI+1},\ldots,b_{\MORE{\NEBENPOWEX}} \in \ZZp[G]}. %
 Since \m{w>\MORE{\NEBENPOWEX}}, the reduction modulo~\m{\maxideal} of the left side of equation~\eqref{eqf43j93jr9348} is in \m{\INAAA} and represents a generator of \m{\ind_{KZ}^G \QUOTIENTI(w)} by equation~\eqref{eq394j9jh39fj93jg4398gj49jg498jg}. \hspace{\stretch{1}} \m{\square}

\emph{Case 2. } \label{case2} %
 Let %
 \m{\delta=\MORE{\NEBENPOWEX-w}=\MORE{\NEBENPOWEX}-w\in\{\TEMPORARYI,\ldots,p-2-\TEMPORARYI\}}. Then \algn{S_\delta = \{\delta-\TEMPORARYI+1,\ldots,\delta-1\} \subseteq \{1,\ldots,p-3-\TEMPORARYI\},} so
 \algn{M = \MMM{w}{(p-i)_{1\LEQSLANT i \LEQSLANT w}}{\delta} = \left(\binom{p-i}{\delta-j}\right)_{1\LEQSLANT i\LEQSLANT w,\, 1\LEQSLANT j \LEQSLANT\TEMPORARYI-1}.}
 Let
  \algn{A %
  = \left(\binom{p-i}{\delta-j}\right)_{1\LEQSLANT i\LEQSLANT w,\, 0\LEQSLANT j \LEQSLANT\TEMPORARYI-1},} so that \m{M} is %
   obtained from \m{A} by removing the column indexed  \m{j=0}. Then the \m{\TEMPORARYI\times\TEMPORARYI} submatrix of \m{A} consisting of the rows indexed \m{i\in\{2,\ldots,\TEMPORARYI+1\}} is obtained from
     \algn{\textstyle\MATRIX_{\TEMPORARYI-1}(p-1-\TEMPORARYI,\delta+1-\TEMPORARYI)} by reducing modulo \m{\maxideal} and permuting the rows and columns, so it is invertible by part~\eqref{vander2} of lemma~\ref{LEMMAX8} because \m{\delta<p} and \m{p\nmid \binom{p-2}{\delta}}. %
 This implies that the first row of \m{A} is a linear combination of the other rows of \m{A}, so there exist \m{\lambda_1,\ldots,\lambda_{w}\in\FF_p} such that 
 \algn{((-1)^{p-1}\lambda_1,\ldots,(-1)^{p-w}\lambda_{w})\,A=0} and \m{\lambda_1=1}.
 Then the quadruple \algn{%
(w,\delta,(p-i)_{1\LEQSLANT i \LEQSLANT w},(\lambda_i)_{1\LEQSLANT i \LEQSLANT w})} satisfies equation~\eqref{eqg3oigj98j45g9j58ghufieofji3u5}, %
  so lemma~\ref{LEMMAX9} %
   implies that
  \Algn{\label{eqfi43jfoij394fj9348jf}\nonumber
 & \textstyle
\sum_{i=1}^w [\lambda_i] \sum_{\TEMPORARYA \in \FF_p} {[\TEMPORARYA]^{p-i} \begin{Lmatrix}
    1 & 0 \\ [\TEMPORARYA] & 1
    \end{Lmatrix}}\sqbrackets{}{\INZB}{\QQp}{\GEN_{{w-i+1}}} \\& \null\qquad\textstyle + \sum_{i=1}^w (-1)^{p-i+\NEBENPOWEX}[\lambda_i]  {\binom{p-i}{\delta}\begin{Lmatrix}
0 & -1 \\ 1 & 0
\end{Lmatrix}}\sqbrackets{}{\INZB}{\QQp}{\GEN_{\delta}} %
 \EQUIVZEROIDEAL \BIGO{\QUARK^\TEMPORARYI a^{-1}}.}
  The fact that \m{((-1)^{p-1}\lambda_1,\ldots,(-1)^{p-w}\lambda_{w})} is orthogonal to the column of \m{A} indexed \m{j=0} implies that \algn{\sum_{i=1}^w (-1)^{p-i}[\lambda_i]  \binom{p-i}{\delta} = \BIGO{p}, }
  so equation~\eqref{eqfi43jfoij394fj9348jf} implies that
    \Algn{\label{eqf934jf9834u98}
  \textstyle
\sum_{i=1}^w [\lambda_i] \sum_{\TEMPORARYA \in \FF_p} {[\TEMPORARYA]^{p-i} \begin{Lmatrix}
    1 & 0 \\ [\TEMPORARYA] & 1
     \end{Lmatrix}}\sqbrackets{}{\INZB}{\QQp}{\GEN_{{w-i+1}}}
 \EQUIVZEROIDEAL \BIGO{\QUARK^\TEMPORARYI a^{-1}}.}
 Since \m{\lambda_1=1}, the reduction modulo~\m{\maxideal} of the left side of equation~\eqref{eqf934jf9834u98} is in \m{\INAAA} and represents a generator of \m{\ind_{KZ}^G \QUOTIENTI(w)} by equation~\eqref{eq394j9jh39fj93jg4398gj49jg498jg}. \hspace{\stretch{1}} \m{\square} %

\emph{Case 3. } \label{case3} 
Let \m{\delta=\TEMPORARYI} and \m{u=p+2\TEMPORARYI-\MORE{\NEBENPOWEX}}.
 Then \m{S_\delta = \{1,\ldots,\TEMPORARYI-1\}}, so
 \algn{M = \MMM{\TEMPORARYI}{(u-i)_{1\LEQSLANT i \LEQSLANT \TEMPORARYI}}{\delta} = \left(\binom{u-i}{j}\right)_{1\LEQSLANT i\LEQSLANT \TEMPORARYI,\, 1\LEQSLANT j \LEQSLANT\TEMPORARYI-1}.}
As \m{\MORE{\NEBENPOWEX}-w=\MORE{\NEBENPOWEX}-\TEMPORARYI\not\in\{0,\ldots,\TEMPORARYI-1\}}, we have \m{\MORE{\NEBENPOWEX}\GEQSLANT 2\TEMPORARYI}, so \m{u\in\{2\TEMPORARYI+1,\ldots,p\}}. 
Let
\algn{A = \left(\binom{u-i}{j}-\VERIFYIF{(i,j)=(1,\TEMPORARYI)}\right)_{1\LEQSLANT i,j\LEQSLANT \TEMPORARYI},} so that \m{M} is %
  obtained from \m{A} by removing the column indexed  \m{j=\TEMPORARYI}.
 If \m{A^\prime} is obtained from  \m{\MATRIX_{\TEMPORARYI-1}(u-\TEMPORARYI,1)} by adding \m{-1} to the bottom-right entry, then
\algn{
    \textstyle \det A^\prime & \textstyle = \det \MATRIX_{\TEMPORARYI-1}(u-\TEMPORARYI,1) - \det \MATRIX_{\TEMPORARYI-2}(u-\TEMPORARYI,1) \\ 
    & \textstyle = \frac{(u-\TEMPORARYI)\cdots(u-1)}{\TEMPORARYI!}-\frac{(u-\TEMPORARYI)\cdots(u-2)}{(\TEMPORARYI-1)!} = \binom{u-2}{\TEMPORARYI}
}
when \m{\TEMPORARYI>1}, and \m{\det A^\prime=u-2} when \m{\TEMPORARYI=1}.
 Since \m{u-2\in \{2\TEMPORARYI-1,\ldots,p-2\}}, it follows that \m{A^\prime} is invertible over \m{\ZZ_p}. As \m{A} is obtained from \m{A^\prime} by reducing modulo \m{\maxideal} and permuting the rows, %
  \m{A} is invertible over \m{\FF_p}.  Therefore there exist \m{\lambda_1,\ldots,\lambda_{\TEMPORARYI}\in\FF_p} such that 
 \Algn{\label{eq398jug983j3ijfi3u4jfi3u4}((-1)^{u-1}\lambda_1,\ldots,(-1)^{u-\TEMPORARYI}\lambda_{\TEMPORARYI})\,A=(0,\ldots,0,(-1)^{\NEBENPOWEX}).}
  So the quadruple
\algn{%
(\TEMPORARYI,\delta,(u-i)_{1\LEQSLANT i \LEQSLANT \TEMPORARYI},(\lambda_i)_{1\LEQSLANT i \LEQSLANT \TEMPORARYI})} satisfies equation~\eqref{eqg3oigj98j45g9j58ghufieofji3u5}, %
  and lemma~\ref{LEMMAX9} implies that
    \Algn{\label{eq0g340f0943ir0349}\nonumber
 & \textstyle
\sum_{i=1}^\TEMPORARYI [\lambda_i] \sum_{\TEMPORARYA \in \FF_p} {[\TEMPORARYA]^{u-i} \begin{Lmatrix}
    1 & 0 \\ [\TEMPORARYA] & 1
    \end{Lmatrix}}\sqbrackets{}{\INZB}{\QQp}{\GEN_{{\TEMPORARYI-i+1}}} \\& \null\qquad\textstyle \vphantom{\null\qquad\textstyle + \sum_{i=1}^\TEMPORARYI (-1)^{u-i+\NEBENPOWEX}[\lambda_i]  {\binom{u-i}{\TEMPORARYI}\begin{Lmatrix}
0 & -1 \\ 1 & 0
\end{Lmatrix}}\sqbrackets{}{\INZB}{\QQp}{\GEN_{\TEMPORARYI}} %
 \EQUIVZEROIDEAL \BIGO{\QUARK^\TEMPORARYI a^{-1}}} + \sum_{i=1}^\TEMPORARYI (-1)^{u-i+\NEBENPOWEX}[\lambda_i]  {\binom{u-i}{\TEMPORARYI}\begin{Lmatrix}
0 & -1 \\ 1 & 0
\end{Lmatrix}}\sqbrackets{}{\INZB}{\QQp}{\GEN_{\TEMPORARYI}} %
 \EQUIVZEROIDEAL \BIGO{\QUARK^\TEMPORARYI a^{-1}}.} Equation~\eqref{eq398jug983j3ijfi3u4jfi3u4} implies that 
    \algn{\sum_{i=1}^\TEMPORARYI (-1)^{u-i+\NEBENPOWEX}[\lambda_i]  \binom{u-i}{\TEMPORARYI} = (-1)^{u-1+\NEBENPOWEX}[\lambda_1] + 1 + \BIGO{p} = [\lambda_1] + 1 + \BIGO{p}, }
  and therefore equation~\eqref{eq0g340f0943ir0349} implies that
     \Algn{\label{eq298ut9284ut298ur983ru32}\nonumber
 & \textstyle
\sum_{i=1}^\TEMPORARYI [\lambda_i] \left(\sum_{\TEMPORARYA \in \FF_p} {[\TEMPORARYA]^{u-i} \begin{Lmatrix}
    1 & 0 \\ [\TEMPORARYA] & 1
    \end{Lmatrix}}+\VERIFYIF{i=1}\begin{Lmatrix}
0 & -1 \\ 1 & 0
\end{Lmatrix}\right)\sqbrackets{}{\INZB}{\QQp}{\GEN_{{\TEMPORARYI-i+1}}} \\& \null\qquad\textstyle \vphantom{\null\qquad\textstyle + \sum_{i=1}^\TEMPORARYI (-1)^{u-i+\NEBENPOWEX}[\lambda_i]  {\binom{u-i}{\TEMPORARYI}\begin{Lmatrix}
0 & -1 \\ 1 & 0
\end{Lmatrix}}\sqbrackets{}{\INZB}{\QQp}{\GEN_{\TEMPORARYI}} %
 \EQUIVZEROIDEAL \BIGO{\QUARK^\TEMPORARYI a^{-1}}}+ \begin{Lmatrix}
0 & -1 \\ 1 & 0
\end{Lmatrix}\sqbrackets{}{\INZB}{\QQp}{\GEN_{\TEMPORARYI}} %
 \EQUIVZEROIDEAL \BIGO{\QUARK^\TEMPORARYI a^{-1}}.}
 By equation~\eqref{eq394j9jh39fj93jg4398gj49jg498jg},
 the term
 \algn{\left(\sum_{\TEMPORARYA \in \FF_p} {\TEMPORARYA^{u-i} \begin{Lmatrix}
    1 & 0 \\ [\TEMPORARYA] & 1
    \end{Lmatrix}}+\VERIFYIF{i=1}\begin{Lmatrix}
0 & -1 \\ 1 & 0
\end{Lmatrix}\right)\sqbrackets{}{\INZB}{\FFp}{\GENN_{{\TEMPORARYI-i+1}}}}
represents the trivial element of \m{\ind_{KZ}^G I_{p^{n-1}-1,\TEMPORARYI}} when \m{i\in\{2,\ldots,\TEMPORARYI\}}, and  it represents an element of \m{\ind_{KZ}^G \SUBMODULE(\TEMPORARYI)}  and therefore the trivial element of \m{\ind_{KZ}^{G}\QUOTIENTI(\TEMPORARYI)} when \m{i=1}. So  the reduction modulo~\m{\maxideal} of the left side of equation~\eqref{eq298ut9284ut298ur983ru32} is in \m{\INAAA} and represents a generator of \m{\ind_{KZ}^G \QUOTIENTI(w)}. \hspace{\stretch{1}} \m{\square}%

\emph{Case 4. }\label{case4} Let \m{\delta=w}.  Then \m{S_\delta = \{w-\TEMPORARYI+1,\ldots,w-1\}}, so
 \algn{M = \MMM{\TEMPORARYI}{(w-\MORE{\NEBENPOWEX}+i)_{1\LEQSLANT i \LEQSLANT \TEMPORARYI}}{\delta} = \left(\binom{w-\MORE{\NEBENPOWEX}+i}{w-j}\right)_{1\LEQSLANT i\LEQSLANT \TEMPORARYI,\, 1\LEQSLANT j \LEQSLANT\TEMPORARYI-1}.}
Let
\algn{A = \left(\binom{w-\MORE{\NEBENPOWEX}+i}{w-j}\right)_{1\LEQSLANT i\LEQSLANT \TEMPORARYI,\, 0\LEQSLANT j \LEQSLANT\TEMPORARYI-1},} so that \m{M} is %
  obtained from \m{A} by removing the column indexed  \m{j=0}.  Then \m{A} is obtained from
     \algn{\textstyle\MATRIX_{\TEMPORARYI-1}(w-\MORE{\NEBENPOWEX}+1,w-\TEMPORARYI+1)} by reducing modulo \m{\maxideal} and permuting the columns, so it is invertible by part~\eqref{vander2} of lemma~\ref{LEMMAX8} because \m{w<p} and \m{p\nmid \binom{w-\MORE{\NEBENPOWEX}+\TEMPORARYI}{w}}. So there exist \m{\lambda_1,\ldots,\lambda_{\TEMPORARYI}\in\FF_p} such that 
 \Algn{\label{eq34t43t53t4g4g5h56}((-1)^{w-\MORE{\NEBENPOWEX}+1}\lambda_1,\ldots,(-1)^{w-\MORE{\NEBENPOWEX}+\TEMPORARYI}\lambda_{\TEMPORARYI})\,A=((-1)^{\NEBENPOWEX},0,\ldots,0),}
  the quadruple
\algn{%
(\TEMPORARYI,\delta,(w-\MORE{\NEBENPOWEX}+i)_{1\LEQSLANT i \LEQSLANT \TEMPORARYI},(\lambda_i)_{1\LEQSLANT i \LEQSLANT \TEMPORARYI})} satisfies equation~\eqref{eqg3oigj98j45g9j58ghufieofji3u5}, %
  and lemma~\ref{LEMMAX9} implies that
    \Algn{\label{eqf0439r39rf3498jugu548y8uje}\nonumber
 & \textstyle
\sum_{i=1}^\TEMPORARYI [\lambda_i] \sum_{\TEMPORARYA \in \FF_p} {[\TEMPORARYA]^{w-\MORE{\NEBENPOWEX}+i} \begin{Lmatrix}
    1 & 0 \\ [\TEMPORARYA] & 1
    \end{Lmatrix}}\sqbrackets{}{\INZB}{\QQp}{\GEN_{i}} \\& \null\qquad\textstyle  + \sum_{i=1}^\TEMPORARYI (-1)^{w+i}[\lambda_i]  {\binom{w-\MORE{\NEBENPOWEX}+i}{w}\begin{Lmatrix}
0 & -1 \\ 1 & 0
\end{Lmatrix}}\sqbrackets{}{\INZB}{\QQp}{\GEN_{w}} %
 \EQUIVZEROIDEAL \BIGO{\QUARK^\TEMPORARYI a^{-1}}.} Equation~\eqref{eq34t43t53t4g4g5h56} implies that 
    \algn{\sum_{i=1}^\TEMPORARYI (-1)^{w+i}[\lambda_i]  \binom{w-\MORE{\NEBENPOWEX}+i}{w} = 1 + \BIGO{p}, }
  and therefore equation~\eqref{eqf0439r39rf3498jugu548y8uje} implies that
  \Algn{\label{eqf43j983u5yur9eiw}\begin{Lmatrix}
0 & -1 \\ 1 & 0
\end{Lmatrix}\sqbrackets{}{\INZB}{\QQp}{\GEN_{w}} %
+ \sum_{i=1}^{\TEMPORARYI} b_i \sqbrackets{}{\INZB}{\QQp}{\GEN_{i}} \EQUIVZEROIDEAL \BIGO{\QUARK^\TEMPORARYI a^{-1}},}
  for some \m{b_{1},\ldots,b_\TEMPORARYI \in \ZZp[G]}. Since \m{w>\TEMPORARYI}, the reduction modulo~\m{\maxideal} of the left side of equation~\eqref{eqf43j983u5yur9eiw} is in \m{\INAAA} and represents a generator of \m{\ind_{KZ}^G \QUOTIENTI(w)} by equation~\eqref{eq394j9jh39fj93jg4398gj49jg498jg}. \hspace{\stretch{1}} \m{\square}%

\emph{Case 5. } \label{case5} This case is nearly identical to case~3, the only difference being that we define 
 \m{u=2\TEMPORARYI-\MORE{\NEBENPOWEX}+1} instead of 
 \m{u=p+2\TEMPORARYI-\MORE{\NEBENPOWEX}}, so we omit the full details. \hspace{\stretch{1}} \m{\square}%

The five constructions in the five different cases show that the kernel of \m{\QUOTIENTMAP} contains a representative of a generator for each  representation in the set
 \algn{  \left\{\ind_{KZ}^{G} \QUOTIENTI(l) \,|\, l \in \{1,\ldots,p-1\}\backslash \{\MORE{\NEBENPOWEX-w}\,|\,w\in \{0,\ldots,\TEMPORARYI-1\}\}\right\},} completing the proof of proposition~\ref{PROP1}.
\null\hspace{\stretch{1}} {\BLACKSQUARE}

  \section{Proof of proposition~\ref{PROP2}}\label{secprop2}

Equation~\eqref{EQUATION1} implies that the quotient map \m{I_{p^{n-1}-1,\TEMPORARYI} \xtwoheadrightarrow{} \QUOTIENTI(\TEMPORARYI)} is given by
 \algn{
    \textstyle f \in I_{p^{n-1}-1,\TEMPORARYI} \xmapsto{ \ \ } \sum_{u,v} f(u,v) (vX-uY)^{p-2} \in \QUOTIENTI(\TEMPORARYI).}
    Here we view \m{\QUOTIENTI(\TEMPORARYI)\cong \sigma_{p-2}(\TEMPORARYI)} as the module of homogeneous polynomials in \m{X} and \m{Y} of total degree \m{p-2}. Equation~\eqref{eq394j9jh39fj93jg4398gj49jg498jg} implies that
    \Algn{\label{eqclass1}\nonumber
    \textstyle  \textstyle {(-1)^n}\sqbrackets{}{\INZB}{\FFp}{\GENN_\TEMPORARYI} & \textstyle  \text{ represents }{1}\sqbracketss{}{K Z}{\FFp}{X^{p-2}},
   \\ \textstyle {(-1)^n\begin{Lmatrix}
 0 & -1 \\ 1 & 0
 \end{Lmatrix}}\sqbrackets{}{\INZB}{\FFp}{\GENN_\TEMPORARYI} & \textstyle \text{ represents }{1}\sqbracketss{}{K Z}{\FFp}{Y^{p-2}}.
} Let \m{T_\sigma} be the Hecke operator in \m{\End_G(\ind_{KZ}^G\QUOTIENTI(\TEMPORARYI))} corresponding to the double coset of \m{\begin{Smatrix} p&0 \\ 0 &1\end{Smatrix}}. Since \m{p-2>0}, the formula for \m{T} given in section~2.1 of~\cite{Bre03b}
 implies that \algn{\nonumber T_\sigma ({1}\sqbracketss{}{K Z}{\FFp}{X^{p-2}}) & \textstyle = \sum_{\MU\in\FF_p} {\begin{Lmatrix}
   p & [\MU] \\ 0  & 1
    \end{Lmatrix} }{}\sqbracketss{}{K Z}{\FFp}{X^{p-2}},\\
    T_\sigma ({1}\sqbracketss{}{K Z}{\FFp}Y^{p-2}) & \textstyle = -\sum_{\MU\in\FF_p} [\MU]^{p-2}{\begin{Lmatrix}
   p & [\MU] \\ 0  & 1
    \end{Lmatrix} }{}\sqbracketss{}{K Z}{\FFp}{X^{p-2}} + {\begin{Lmatrix}
  1 & 0 \\ 0  & p
    \end{Lmatrix} }{}\sqbracketss{}{K Z}{\FFp}{Y^{p-2}}.
    }
    Therefore 
    \Algn{\label{eqclass2}\nonumber
    \textstyle  \textstyle (-1)^n\sum_{\MU\in\FF_p} {\begin{Lmatrix}
   p & [\MU] \\ 0  & 1
    \end{Lmatrix} }\sqbrackets{}{\INZB}{\FFp}{\GENN_\TEMPORARYI} & \textstyle  \text{ represents }T_\sigma({1}\sqbracketss{}{K Z}{\FFp}{X^{p-2}}),
   \\ \textstyle (-1)^{n-1}\left( \sum_{\MU\in\FF_p} [\MU]^{p-2}\begin{Lmatrix}
   p & [\MU] \\ 0  & 1
    \end{Lmatrix} + \begin{Lmatrix}
    0 & 1 \\ -p & 0
    \end{Lmatrix}\right) \sqbrackets{}{\INZB}{\FFp}{\GENN_\TEMPORARYI} & \textstyle \text{ represents }T_\sigma({1}\sqbracketss{}{K Z}{\FFp}{Y^{p-2}}).
}

 Note that \m{(-1)^\NEBENPOWEX = (-1)^{2\TEMPORARYI-1}=-1}. For  \m{u \in\ZZ_{\GEQSLANT 0}},  \m{s \in \{1,\ldots,p-1\}}, and  \m{\TEMPORARYA\in\FF_p}, let
\algn{c_\TEMPORARYA(u,s)\textstyle = (-1)^{u+\NEBENPOWEX} \sum_{\MU \in \FF_p} \sum_{j=0}^{u} \binom{u}{j}[\MU]^{u-j} [\TEMPORARYA]^{\LESS{-j+s-\NEBENPOWEY}} \begin{Lmatrix}
        0 & -1 \\ 1 & [\MU]
        \end{Lmatrix} \in \ZZ_p[G]}
    be the constant defined in lemma~\ref{LEMMAX5}.  For \m{\alpha\in\ZZ_{\GEQSLANT0}}, let \m{\CONSTANT{\alpha}=-\frac{(1-\ZETA)^{\alpha}}{\alpha!}}.  Let \Algn{\label{eq93jf93ujf98uf9843u}\TEMPORARYMU = -\frac{C_{\TEMPORARYI-1}}{a}-\frac{a}{C_{\TEMPORARYI}} & \textstyle= \frac{1}{aC_\TEMPORARYI}\left(\frac{1}{\TEMPORARYI!(\TEMPORARYI-1)!}(\ZETA-1)^{2\TEMPORARYI-1}-a^2\right) \in \QQp.}
    The set  \m{\CIRCLES_{\TEMPORARYI}} is designed precisely so that \Algn{\label{eqReduct}a \in \CIRCLES_{\TEMPORARYI} \Longleftrightarrow \TEMPORARYMU\in\ZZp.} If \m{\TEMPORARYMU\in\ZZp} then \m{\TADICVALUATION(a)=\TEMPORARYI-\frac12} and \m{\TADICVALUATION(a^2+C_\TEMPORARYI C_{\TEMPORARYI-1})\GEQSLANT 2\TEMPORARYI-\frac12}, so
    \Algn{\label{eq09g54ig9045i5irf094i5f0}\TEMPORARYMU = -\left(1+\BIGO{\QUARK^{1/2}}\right) \frac{C_{\TEMPORARYI-1}}{a^3}\left(\frac{1}{\TEMPORARYI!(\TEMPORARYI-1)!}(\ZETA-1)^{2\TEMPORARYI-1}-a^2\right),} and the reduction of \m{\TEMPORARYMU} modulo \m{\maxideal} is \m{\phi}.
    To begin with, let us prove the following equations.
\newcommand{\tempheightaerb}{\vphantom{ \textstyle \EQUIVZEROIDEAL  \sum_{\MU \in \FF_p} { \left(\frac{C_\TEMPORARYI}{a} +  (\ZETA-1)\TEMPORARYMU [\MU]^{p-1}\right)  \begin{Lmatrix}
       p & [\MU] \\ 0 & 1
 \end{Lmatrix} }\sqbrackets{}{\INZB}{\QQp}{\GEN_\TEMPORARYI}+  \BIGO{\QUARK} \cup \BIGO{a^2\QUARK^{2-2\TEMPORARYI}}}}   \Algn{\label{eq0f943if0943ik09f}\textstyle \textstyle\sum_{\TEMPORARYA \in \FF_p} {[\TEMPORARYA]^{p-1-\alpha} \begin{Lmatrix}
    1 & 0 \\ [\TEMPORARYA] & 1
    \end{Lmatrix}}\sqbrackets{}{\INZB}{\QQp}{\GEN}  &\textstyle\tempheightaerb\EQUIVZEROIDEAL \BIGO{\QUARK^\alpha a^{-1}} \text{ for } \alpha \in \{1,\ldots,p-1\}.\\ \textstyle\label{eq029if984ujf9jf9wjf} \textstyle 1\sqbrackets{}{\INZB}{\QQp}{\GEN_\alpha} & \textstyle\tempheightaerb\EQUIVZEROIDEAL \BIGO{a\QUARK^{-\alpha}}\text{ for } \alpha \in \{0,\ldots,p-1\}.
    \\
    \textstyle \label{eqfu3h4g8fh34hf3487fh387}\nonumber\textstyle \sum_{\TEMPORARYA,\MU \in \FF_p}  [\TEMPORARYA]^{p-1-\TEMPORARYI} \begin{Lmatrix}
        1 & [\MU] \\ 0 & 1
        \end{Lmatrix}\begin{Lmatrix}
    1 & 0 \\ [\TEMPORARYA] & 1
    \end{Lmatrix}\sqbrackets{}{\INZB}{\QQp}{\GEN} & \textstyle \tempheightaerb\EQUIVZEROIDEAL\begin{Lmatrix}
0 & -1 \\ 1 & 0
\end{Lmatrix} \sum_{\TEMPORARYA\in\FF_p}\begin{Lmatrix}
    1 & 0 \\ [\TEMPORARYA] & 1
    \end{Lmatrix}\sqbrackets{}{\INZB}{\QQp}{\GEN_{\MORE{\TEMPORARYI-1}}} + \ERRR^\prime,%
 \\ \textstyle%
 \text{where } \frac{1}{\ZETA-1}\ERRR ^\prime & \textstyle \tempheightaerb\EQUIVZEROIDEAL %
    \frac{1}{\TEMPORARYI}   %
     \sqbrackets{}{\INZB}{\QQp}{\GEN_\TEMPORARYI} + \BIGO{\QUARK^{\TEMPORARYI} a^{-1}}\cup \BIGO{a\QUARK^{1-\TEMPORARYI}}.}

\emph{Proof of equations~\eqref{eq0f943if0943ik09f} \&~\eqref{eq029if984ujf9jf9wjf}. }
If \m{\alpha\in\{1,\ldots,p-1\}} then both equations follow from lemma~\ref{LEMMAX6} applied to \m{\alpha} and  the set of lifts \m{\LCLASSTEICH}. If  \m{\alpha=0} then equation~\eqref{eq029if984ujf9jf9wjf}  follows from equation~\eqref{eq3498ut4398jg394j3498} and parts~\eqref{propeta2} and~\eqref{propeta5} of lemma~\ref{LEMMAX3.01}.
\hspace{\stretch{1}} \m{\square}

\emph{Proof of equation~\eqref{eqfu3h4g8fh34hf3487fh387}. }
 Since 
 \algn{\textstyle\textstyle c_\TEMPORARYA(0,\MORE{\TEMPORARYI-1}) = - \sum_{\MU \in \FF_p}  [\TEMPORARYA]^{p-1-\TEMPORARYI} \begin{Lmatrix}
        0 & -1 \\ 1 & [\MU]
        \end{Lmatrix},}
if we apply lemma~\ref{LEMMAX5} to \m{(u,s)= (0,\MORE{\TEMPORARYI-1})}
 and multiply both sides of the resulting equation by \m{\begin{Smatrix} 0 & -1 \\ 1 & 0
 \end{Smatrix}} we get 
\algn{\sum_{\TEMPORARYA,\MU \in \FF_p}  [\TEMPORARYA]^{p-1-\TEMPORARYI} \begin{Lmatrix}
        1 & [\MU] \\ 0 & 1
        \end{Lmatrix}\begin{Lmatrix}
    1 & 0 \\ [\TEMPORARYA] & 1
    \end{Lmatrix}\sqbrackets{}{\INZB}{\QQp}{\GEN} & \textstyle \tempheightaerb\EQUIVZEROIDEAL\begin{Lmatrix}
0 & -1 \\ 1 & 0
\end{Lmatrix} \sum_{\TEMPORARYA\in\FF_p}\begin{Lmatrix}
    1 & 0 \\ [\TEMPORARYA] & 1
    \end{Lmatrix}\sqbrackets{}{\INZB}{\QQp}{\GEN_{\MORE{\TEMPORARYI-1}}} + \ERRR^{\prime\prime},} and      we can use the formula for the error term given in lemma~\ref{LEMMAX5} with the sets of lifts \algn{(\LCLASS_\TEMPORARYAA)_{\TEMPORARYAA\in\FF_p^\times} = \left(\{-[\MU]-[\TEMPORARYAA]\,|\,\MU\in\FF_p\}\right)_{\TEMPORARYAA\in\FF_p^\times}} %
 to compute that
    \algn{\frac{1}{\ZETA-1}\ERRR^{\prime\prime}\EQUIVZEROIDEAL  \sum_{\TEMPORARYAA\in\FF_p^\times,\, \MU\in\FF_p} \sum_{j=1}^{p-1} 
        \frac{1}{j} ([\MU]+[\TEMPORARYAA])^{j} [\TEMPORARYAA]^{\TEMPORARYI-j} A_{\TEMPORARYAA,\MU}\sqbrackets{}{\INZB}{\QQp}{\GEN} + \BIGO{\QUARK},}
     where \m{%
      A_{\TEMPORARYAA,\MU} = \begin{Smatrix} 0 & -1 \\ 1 & 0
 \end{Smatrix}\begin{Smatrix}
            1 & 0 \\ -[\MU]-[\TEMPORARYAA] & 1 \end{Smatrix}  \begin{Smatrix}
            1/[\TEMPORARYAA] & 1 \\ 0 & [\TEMPORARYAA]
            \end{Smatrix} = \begin{Smatrix}
        1 & [\MU]
       \\ 0 & 1     \end{Smatrix}\begin{Smatrix}
             1 & 0 \\ 1/[\TEMPORARYAA] & 1
             \end{Smatrix}}. 
           By part~\eqref{comb10.3} of lemma~\ref{LEMMAX3.04},
         \Algn{\label{eqg9j59f8jw98fr983u398tu}\nonumber\frac{1}{\ZETA-1}\ERRR ^{\prime\prime} & \textstyle\EQUIVZEROIDEAL  \sum_{\TEMPORARYAA\in\FF_p^\times,\, \MU\in\FF_p}   \sum_{l=1}^{p-1}\frac{(-1)^l}{l} [\MU]^l [\TEMPORARYAA]^{\TEMPORARYI-l} \begin{Lmatrix}
        1 & [\MU] \\ 0 & 1
        \end{Lmatrix}\begin{Lmatrix}
        1 & 0 \\ 1/[\TEMPORARYAA] & 1  \end{Lmatrix}\sqbrackets{}{\INZB}{\QQp}{\GEN} + \BIGO{\QUARK}
     \\ \nonumber& \textstyle \EQUIVZEROIDEAL \sum_{\TEMPORARYA\in\FF_p^\times,\, \MU\in\FF_p}
    \sum_{l=1}^{p-1}\frac{(-1)^l}{l} [\MU]^l [\TEMPORARYA]^{l-\TEMPORARYI} %
      \begin{Lmatrix}
        1 & [\MU] \\ 0 & 1
        \end{Lmatrix}\begin{Lmatrix}
        1 & 0 \\ [\TEMPORARYA] & 1  \end{Lmatrix}\sqbrackets{}{\INZB}{\QQp}{\GEN} + \BIGO{\QUARK}
        \\ &\nonumber \textstyle \EQUIVZEROIDEAL \sum_{\TEMPORARYA, \MU\in\FF_p}
    \sum_{l=1}^{p-1}\frac{(-1)^l}{l} [\MU]^l [\TEMPORARYA]^{\LESS{l-\TEMPORARYI}} %
    \begin{Lmatrix}
        1 & [\MU] \\ 0 & 1
        \end{Lmatrix}\begin{Lmatrix}
        1 & 0 \\ [\TEMPORARYA] & 1  \end{Lmatrix}\sqbrackets{}{\INZB}{\QQp}{\GEN}
        \\ & \textstyle\null\qquad \phantom{\EQUIVZEROIDEAL} - \sum_{\MU\in\FF_p}
    \frac{(-1)^\TEMPORARYI}{\TEMPORARYI} [\MU]^\TEMPORARYI  %
     \begin{Lmatrix}
        1 & [\MU] \\ 0 & 1
        \end{Lmatrix}\sqbrackets{}{\INZB}{\QQp}{\GEN} + \BIGO{\QUARK}.
        } 
   For the second congruence we change the variable \m{\TEMPORARYAA} to \m{\TEMPORARYA=1/\TEMPORARYAA}, and for the third congruence we rewrite the sum ``\m{\sum_{\TEMPORARYA\in\FF_p^\times}}'' as ``\m{\sum_{\TEMPORARYA\in\FF_p}} minus the term with \m{\TEMPORARYA=0}''.
      For \m{l
    \in\{0\}\cup\{\TEMPORARYI,\ldots,p-1\}},
  \Algn{\label{eqf03j498rf3498u9834u}\sum_{\TEMPORARYA \in \FF_p} {[\TEMPORARYA]^{\LESS{l-\TEMPORARYI}} \begin{Lmatrix}
    1 & 0 \\ [\TEMPORARYA] & 1
    \end{Lmatrix}}\sqbrackets{}{\INZB}{\QQp}{\GEN} \EQUIVZEROIDEAL \BIGO{\QUARK^{\TEMPORARYI}a^{-1}},} %
 by equation~\eqref{eq0f943if0943ik09f} with \m{\alpha\in\{\TEMPORARYI,\ldots,p-1\}}.  If \m{\TEMPORARYI=1} then equation~\eqref{eqfu3h4g8fh34hf3487fh387} follows from equations~\eqref{eqg9j59f8jw98fr983u398tu} and~\eqref{eqf03j498rf3498u9834u}. %
 Suppose that \m{\TEMPORARYI>1}.  Let   \algn{c_\TEMPORARYA^\ast = \sum_{i=1}^{\TEMPORARYI-1} \frac{1}{i}\binom{\TEMPORARYI-1}{i}c_\TEMPORARYA(p-i,\TEMPORARYI-i).}
For \m{i\in\{1,\ldots,\TEMPORARYI-1\}},
         \algn{c_\TEMPORARYA(p-i,\TEMPORARYI-i) & \textstyle = (-1)^{i} \sum_{\MU \in \FF_p} \sum_{j=0}^{p-i} \binom{p-i}{j}[\MU]^{p-i-j} [\TEMPORARYA]^{\LESS{1-\TEMPORARYI-i-j}} \begin{Lmatrix}
        0 & -1 \\ 1 & [\MU]
        \end{Lmatrix} \\ & \textstyle = (-1)^{i} \sum_{\MU \in \FF_p} \sum_{l=0}^{p-i} \binom{p-i}{l}[\MU]^{l} [\TEMPORARYA]^{\LESS{l-\TEMPORARYI}} \begin{Lmatrix}
        0 & -1 \\ 1 & [\MU]
        \end{Lmatrix}. } For the second equality we change the variable \m{j} to \m{l=p-i-j}. So part~\eqref{comb10.4} of lemma~\ref{LEMMAX3.04} implies that \algn{c_\TEMPORARYA^\ast =  -\sum_{\MU \in \FF_p} \sum_{l=1}^{\TEMPORARYI-1} \frac{(-1)^{l}}{l} [\MU]^{l} [\TEMPORARYA]^{\LESS{l-\TEMPORARYI}} \begin{Lmatrix}
        0 & -1 \\ 1 & [\MU]
        \end{Lmatrix} + \sum_{l\in\{0\}\cup\{\TEMPORARYI,\ldots,p-1\}} b_l [\TEMPORARYA]^{\LESS{l-\TEMPORARYI}} + \BIGO{p},} for some \m{b_0,b_{\TEMPORARYI},\ldots,b_{p-1}\in\ZZp[G]}.
         Therefore equations~\eqref{eqg9j59f8jw98fr983u398tu} and~\eqref{eqf03j498rf3498u9834u} imply that \Algn{\label{eqg398u4g983uj94j5g984jg9845ujguj45g9u}\frac{1}{\ZETA-1}\ERRR ^{\prime\prime} \EQUIVZEROIDEAL  \begin{Lmatrix}
        0 & -1 \\ 1 & 0
        \end{Lmatrix} \sum_{\TEMPORARYA \in \FF_p} {c_\TEMPORARYA^\ast \begin{Lmatrix}
    1 & 0 \\ [\TEMPORARYA] & 1
    \end{Lmatrix}}\sqbrackets{}{\INZB}{\QQp}{\GEN}  -
    \frac{(-1)^\TEMPORARYI}{\TEMPORARYI}   %
     \sqbrackets{}{\INZB}{\QQp}{\GEN_\TEMPORARYI} +  \BIGO{\QUARK^{\TEMPORARYI}a^{-1}}.}
         Lemma~\ref{LEMMAX5} applied to \m{(u,s) =(p-i,\TEMPORARYI-i)} for \m{i\in\{1,\ldots,\TEMPORARYI-1\}} implies that %
    \Algn{\label{eqf054j948jg45jgt98ujt}\nonumber
& \textstyle\sum_{\TEMPORARYA \in \FF_p} {c_\TEMPORARYA^\ast \begin{Lmatrix}
    1 & 0 \\ [\TEMPORARYA] & 1
    \end{Lmatrix}}\sqbrackets{}{\INZB}{\QQp}{\GEN} \\ \textstyle & \textstyle\qquad  \EQUIVZEROIDEAL\nonumber
\sum_{i=1}^{\TEMPORARYI-1} \frac{1}{i}\binom{\TEMPORARYI-1}{i} \sum_{\TEMPORARYA \in \FF_p} {[\TEMPORARYA]^{p-i} \begin{Lmatrix}
    1 & 0 \\ [\TEMPORARYA] & 1
    \end{Lmatrix}}\sqbrackets{}{\INZB}{\QQp}{\GEN_{\TEMPORARYI-i}} \\& \null\qquad\qquad\textstyle + \sum_{i=1}^{\TEMPORARYI-1} \frac{(-1)^{i}}{i}\binom{\TEMPORARYI-1}{i}  {\binom{p-i}{\TEMPORARYI}\begin{Lmatrix}
0 & -1 \\ 1 & 0
\end{Lmatrix}}\sqbrackets{}{\INZB}{\QQp}{\GEN_{\TEMPORARYI}} + \BIGO{\QUARK}\nonumber%
 \\ \textstyle & \textstyle\qquad  \EQUIVZEROIDEAL
 -{\frac{1+(-1)^\TEMPORARYI}{\TEMPORARYI}\begin{Lmatrix}
0 & -1 \\ 1 & 0
\end{Lmatrix}}\sqbrackets{}{\INZB}{\QQp}{\GEN_{\TEMPORARYI}} + \BIGO{a\QUARK^{1-\TEMPORARYI}}.}
The second congruence is due to equation~\eqref{eq029if984ujf9jf9wjf} with \m{\alpha\in\{1,\ldots,\TEMPORARYI-1\}} and part~\eqref{comb10.4} of lemma~\ref{LEMMAX3.04}. Equations~\eqref{eqg398u4g983uj94j5g984jg9845ujguj45g9u} and~\eqref{eqf054j948jg45jgt98ujt} imply that equation~\eqref{eqfu3h4g8fh34hf3487fh387} is true in the case when \m{\TEMPORARYI>1} as well.
\hspace{\stretch{1}} \m{\square}

\subsection{\protect\boldmath Completing the proof for \texorpdfstring{\m{{\TEMPORARYI=1}}}{\textnu=1}} 
\label{subsection5.4} In this subsection we assume that \m{\TEMPORARYI=1}, so that  \m{\MORE{\NEBENPOWEX}=1}. %
     Let  \algn{\ERRR_1 = \sum_{\TEMPORARYA \in \FF_p} {(c_\TEMPORARYA(0,1)-c_\TEMPORARYA(p-1,1)) \begin{Lmatrix}
    1 & 0 \\ [\TEMPORARYA] & 1
    \end{Lmatrix}}\sqbrackets{}{\INZB}{\QQp}{\GEN}. }
    Since \algn{c_\TEMPORARYA(0,1)-c_\TEMPORARYA(p-1,1) = - \sum_{\MU \in \FF_p} [\MU]^{p-2}[\TEMPORARYA]^{p-2} \begin{Lmatrix}
          0 & -1 \\ 1 & [\MU]
  \end{Lmatrix}+\sum_{j=0}^{p-3} b_{j}[\TEMPORARYA]^j + \BIGO{p}}  for some \m{b_{0},\ldots,b_{p-3}\in\ZZ_p[G]}, it follows from equation~\eqref{eq0f943if0943ik09f} with \m{\alpha\in\{2,\ldots,p-1\}} and lemma~\ref{LEMMAX6} applied to \m{\alpha=1} and the set of lifts \m{\LCLASSTEICH} that
   \Algn{\label{eq3040394ir093ri0394i}\nonumber\textstyle  \ERRR_1 \nonumber & \textstyle\null \EQUIVZEROIDEAL - \sum_{\TEMPORARYA,\MU \in \FF_p} [\MU]^{p-2} [\TEMPORARYA]^{p-2} \begin{Lmatrix}
         0 & -1 \\ 1 & [\MU]
  \end{Lmatrix} \begin{Lmatrix}
     1 & 0 \\ [\TEMPORARYA] & 1
     \end{Lmatrix}\sqbrackets{}{\INZB}{\QQp}{\GEN} +  \BIGO{\QUARK^2 a^{-1}}
     \\ & \textstyle\null\EQUIVZEROIDEAL - %
     \COVERA{1}\sum_{\MU\in \FF_p} [\MU]^{p-2}  \begin{Lmatrix}
        0 & -1 \\ p & [\MU]
  \end{Lmatrix} \sqbrackets{}{\INZB}{\QQp}{\GEN_1} +  \BIGO{\QUARK^2 a^{-1}}.} 
  It follows from lemma~\ref{LEMMAX5} applied to \m{(u,s)=(0,1)} and \m{(u,s)=(p-1,1)}  that
\Algn{\label{eq3948ujt9438u3948}
\nonumber& \textstyle1\sqbrackets{}{\INZB}{\QQp}{\GEN_1} + \begin{Lmatrix}
0 & -1 \\ 1 & 0
\end{Lmatrix}\sqbrackets{}{\INZB}{\QQp}{\GEN_{p-1}} \\ \textstyle & \textstyle\qquad  \EQUIVZEROIDEAL
 \sum_{\TEMPORARYA \in \FF_p} {(c_\TEMPORARYA(0,1)-c_\TEMPORARYA(p-1,1)) \begin{Lmatrix}
    1 & 0 \\ [\TEMPORARYA] & 1
    \end{Lmatrix}}\sqbrackets{}{\INZB}{\QQp}{\GEN} + \BIGO{\QUARK} = \ERRR_1+\BIGO{\QUARK}.}
        Equation~\eqref{eq3498ut4398jg394j3498} and parts~\eqref{propeta2} and~\eqref{propeta5} of lemma~\ref{LEMMAX3.01}
    imply that
\Algn{\label{eqn66}%
\textstyle \begin{Lmatrix}
    p & 0 \\ 0 & 1
    \end{Lmatrix} \sqbrackets{}{\INZB}{\QQp}{\GEN_0} \EQUIVZEROIDEAL a  \sqbrackets{}{\INZB}{\QQp}{\GEN} + \BIGO{p},}
    and part~\eqref{propeta4} of lemma~\ref{LEMMAX3.01} then implies that
    \Algn{\label{eq2093ri2309ir2039ir2039ri}
    \nonumber    \textstyle 1
    \sqbrackets{}{\INZB}{\QQp}{\GEN_{p-1}} & \textstyle \EQUALZEROIDEAL 1
    \sqbrackets{}{\INZB}{\QQp}{(\GEN_0-\GEN)}
    \\ & \EQUIVZEROIDEAL \left(a\begin{Lmatrix}
    1 & 0 \\ 0 & p
    \end{Lmatrix}-1\right)
 \sqbrackets{}{\INZB}{\QQp}{\GEN} + \BIGO{p}. 
    }
        Equations~\eqref{eq3948ujt9438u3948} and~\eqref{eq2093ri2309ir2039ir2039ri} imply that
    \Algn{\label{eq3409it3049ir0934ri}\nonumber\textstyle {1}\sqbrackets{}{\INZB}{\QQp}{\GEN_{1}} & \textstyle \EQUIVZEROIDEAL -{\begin{Lmatrix}
    0 & -1 \\ 1 & 0
    \end{Lmatrix} }\sqbrackets{}{\INZB}{\QQp}{\GEN_{p-1}} + \ERRR_1+ \BIGO{\QUARK} \\ & \textstyle \EQUIVZEROIDEAL {\left(\begin{Lmatrix}
    0 & -1 \\ 1 & 0
    \end{Lmatrix}-a\begin{Lmatrix}
    0 & -p \\ 1 & 0
    \end{Lmatrix}\right) }\sqbrackets{}{\INZB}{\QQp}{\GEN} + \ERRR_1+ \BIGO{\QUARK}.}
       If we multiply equation~\eqref{eq3409it3049ir0934ri} by \m{a \sum_{\LAMBDA\in\FF_p} \begin{Smatrix}
            p & [\LAMBDA] \\ 0 & 1
            \end{Smatrix}} we get 
    \newcommand{\AUXERRebgefw}{\ERRR_2}\Algn{\label{eq9834ut938u9384t}\nonumber
                & \textstyle a \sum_{\LAMBDA\in\FF_p} {\begin{Lmatrix}
            p & [\LAMBDA] \\ 0 & 1
            \end{Lmatrix} }\sqbrackets{}{\INZB}{\QQp}{\GEN_1} \\ & \textstyle \null\qquad \EQUIVZEROIDEAL a \begin{Lmatrix}
            0 & -1 \\ 1 & 0
            \end{Lmatrix} \sum_{\LAMBDA\in\FF_p}{ \begin{Lmatrix}
             1 & 0 \\ -[\LAMBDA] & p
      \end{Lmatrix}}\sqbrackets{}{\INZB}{\QQp}{\GEN} + \AUXERRebgefw + \BIGO{a\QUARK}\nonumber
            \\ & \textstyle \null\qquad \EQUIVZEROIDEAL \begin{Lmatrix}
            0 & -1 \\ 1 & 0
            \end{Lmatrix}\sum_{\MU\in\FF_p}{\begin{Lmatrix}
             1 & 0 \\ [\MU] & 1
     \end{Lmatrix}}\sqbrackets{}{\INZB}{\QQp}{\GEN_0} + \AUXERRebgefw + \BIGO{a\QUARK}
          \nonumber        \\ & \textstyle \null\qquad \EQUIVZEROIDEAL \begin{Lmatrix}
            0 & -1 \\ 1 & 0
            \end{Lmatrix}\sum_{\MU\in\FF_p}{ \begin{Lmatrix}
             1 & 0 \\ [\MU] & 1
         \end{Lmatrix}}\sqbrackets{}{\INZB}{\QQp}{\GEN_{p-1}} + \AUXERRebgefw + \BIGO{a\QUARK}\cup \BIGO{\QUARK^{2} a^{-1}},\textstyle 
          } where
     \Algn{\label{eq3juf983uju9832}\nonumber\AUXERRebgefw  & \textstyle 
        \EQUIVZEROIDEAL - \sum_{\LAMBDA\in\FF_p} {\begin{Lmatrix}
             p & [\LAMBDA] \\ 0 & 1
            \end{Lmatrix} } \left( C_1 \sum_{\MU\in \FF_p} [\MU]^{p-2}  \begin{Lmatrix}
         0 & -1 \\ p & [\MU]
 \end{Lmatrix} \sqbrackets{}{\INZB}{\QQp}{\GEN_1} + a^2\begin{Lmatrix}
    0 & -p \\ 1 & 0
    \end{Lmatrix}\sqbrackets{}{\INZB}{\QQp}{\GEN}\right)
    \\ & \textstyle\nonumber
        \EQUIVZEROIDEAL - \sum_{\LAMBDA\in\FF_p} {\begin{Lmatrix}
             [\LAMBDA] & -p \\ 1 & 0
            \end{Lmatrix} } \left( C_1 \sum_{\MU\in \FF_p} [\MU]^{p-2}  \begin{Lmatrix}
         p & [\MU] \\ 0 & 1
 \end{Lmatrix} + a^2\begin{Lmatrix}
    0 & 1 \\ -p & 0
    \end{Lmatrix}\right)\sqbrackets{}{\INZB}{\QQp}{\GEN_1}\\&\null\hspace{7.5cm} \qquad+\BIGO{a^3}\cup\BIGO{a\QUARK}.}
    For the second congruence in equation~\eqref{eq9834ut938u9384t} we change the variable \m{\LAMBDA} to \m{\MU=-\LAMBDA} and replace \m{\GEN} with \m{\GEN_0} by using equation~\eqref{eqn66}. For the third congruence in equation~\eqref{eq9834ut938u9384t} we replace \m{\GEN_{0}} with \m{\GEN_{p-1}} by using part~\eqref{propeta4} of lemma~\ref{LEMMAX3.01} and subtracting \algn{\textstyle \begin{Lmatrix}
            0 & -1 \\ 1 & 0
            \end{Lmatrix} \sum_{\MU\in\FF_p}{\begin{Lmatrix} 1 & 0 \\ [\MU] & 1
     \end{Lmatrix}}\sqbrackets{}{\INZB}{\QQp}{\GEN},} which is in \m{\BIGO{\QUARK^{p-1} a^{-1}}} by equation~\eqref{eq0f943if0943ik09f} with \m{\alpha=p-1}.
     The second congruence in equation~\eqref{eq3juf983uju9832} follows from equation~\eqref{eq3409it3049ir0934ri}.
 The first line of equation~\eqref{eq3juf983uju9832} %
    implies that
    \m{\label{eqoifjio3j43jf93jg983gj}\ERRR_2\in\BIGO{a^2}\cup\BIGO{\QUARK}}, as   \m{\ERRR_1=\BIGO{\QUARK a^{-1}}}  due to equation~\eqref{eq3040394ir093ri0394i}. Equation~\eqref{eqfu3h4g8fh34hf3487fh387} implies that
 \Algn{\label{eq20fj9283jr982j3urf23}\nonumber
 \textstyle
\begin{Lmatrix}
0 & -1 \\ 1 & 0
\end{Lmatrix} \sum_{\MU\in\FF_p}{ \begin{Lmatrix}
             1 & 0 \\ [\MU] & 1
         \end{Lmatrix}} \sqbrackets{}{\INZB}{\QQp}{\GEN_{p-1}} 
 & \textstyle \nonumber\null \EQUIVZEROIDEAL  \sum_{\TEMPORARYA ,\LAMBDA\in \FF_p} { [\TEMPORARYA]^{p-2}\begin{Lmatrix}
1 & [\LAMBDA] \\ 0 & 1 
\end{Lmatrix} \begin{Lmatrix}
    1 & 0 \\ [\TEMPORARYA] & 1
    \end{Lmatrix}}\sqbrackets{}{\INZB}{\QQp}{\GEN} 
-  \ERRR^\prime
    \\ & \textstyle \null \EQUIVZEROIDEAL  %
    \COVERA{1}
    \sum_{\LAMBDA\in\FF_p} {\begin{Lmatrix} p & [\LAMBDA] \\ 0 & 1
      \end{Lmatrix}}\sqbrackets{}{\INZB}{\QQp}{\GEN_1}
- \ERRR^\prime.
}
            The second congruence follows from lemma~\ref{LEMMAX6} applied to \m{\alpha=1} and the set of lifts \m{\LCLASSTEICH}. Equations~\eqref{eq9834ut938u9384t} and~\eqref{eq20fj9283jr982j3urf23} imply that 
            \Algn{\label{eq2093ri209ri209rf0fj093}\nonumber
                & \textstyle (1-\COVERASQ{1}) \sum_{\LAMBDA\in\FF_p} {\begin{Lmatrix}
            p & [\LAMBDA] \\ 0 & 1
            \end{Lmatrix} }\sqbrackets{}{\INZB}{\QQp}{\GEN_1} \\ & \textstyle\null\qquad \EQUIVZEROIDEAL a^{-1} \left(\ERRR_2-  \ERRR^\prime\right)+ \BIGO{\QUARK}\cup\BIGO{\QUARK^2 a^{-2}} \EQUIVZEROIDEAL \BIGO{a}\cup\BIGO{\QUARK a^{-1}}.
            } 
            By equation~\eqref{eqReduct},
            \algn{1 - \COVERASQ{1} \in  \BIGO{a}\cup\BIGO{\QUARK a^{-1}} \Longleftrightarrow \TADICVALUATION\left(a^2-\CONSTANT{1}\right)\GEQSLANT \frac32 \Longleftrightarrow a\in \CIRCLES_{1} \Longleftrightarrow \TEMPORARYMU\in\ZZp.}
           If   \m{a\not\in\CIRCLES_1} then after dividing equation~\eqref{eq2093ri209ri209rf0fj093}  by \m{1-\COVERASQ{1}} and reducing modulo \m{\maxideal}, we get by equation~\eqref{eqclass2} that \m{T_\sigma({(-1)^n}\sqbracketss{}{K Z}{\FFp}{X^{p-2}})} is represented by an element of \m{\INAAA  } and therefore  by an element of  the kernel of \m{\QUOTIENTMAP}. If \m{a\in\CIRCLES_1} then \m{\TADICVALUATION(a)=\frac12} and \m{\TEMPORARYMU\in\ZZp}, and equation~\eqref{eq2093ri209ri209rf0fj093} implies that
     \Algn{\label{eq340k309f943f3984f3948uuj9834uj}
     \textstyle \TEMPORARYMU \sum_{\LAMBDA\in\FF_p} {\begin{Lmatrix}
            p & [\LAMBDA] \\ 0 & 1
            \end{Lmatrix} }\sqbrackets{}{\INZB}{\QQp}{\GEN_1} %
            \EQUIVZEROIDEAL a^{-2} \left( \ERRR^\prime-\ERRR_2\right)+ \BIGO{a}.
            }
          As 
          \m{\TEMPORARYMU\in\ZZp} implies that \m{\ZETA-1 = a^2+\BIGO{a^3}}, by equations~\eqref{eqclass1}, \eqref{eqclass2}, \eqref{eq09g54ig9045i5irf094i5f0}, \eqref{eqfu3h4g8fh34hf3487fh387}, and~\eqref{eq3juf983uju9832}, the reductions of 
\algn{\TEMPORARYMU\sum_{\LAMBDA\in\FF_p}\begin{Lmatrix}
                p & [\LAMBDA] \\ 0 & 1
              \end{Lmatrix}\sqbrackets{}{\INZB}{\QQp}{\GEN_{1}} ,\,\,
               a^{-2}\ERRR^\prime ,\,\, -a^{-2}\ERRR_2}
modulo \m{\maxideal}  are
\algn{\phi\, {\textstyle T_\sigma({(-1)^n}\sqbracketss{}{K Z}{\FFp}{X^{p-2}})},\,\, {(-1)^{n}}\sqbracketss{}{K Z}{\FFp}{X^{p-2}},\,\, {\textstyle T_\sigma^2({(-1)^n}\sqbracketss{}{K Z}{\FFp}{X^{p-2}})},} respectively.
 So equation~\eqref{eq340k309f943f3984f3948uuj9834uj} implies that the kernel of \m{\QUOTIENTMAP} contains a representative of \m{(T_\sigma^2-\phi\, T_\sigma+1)({(-1)^n}\sqbracketss{}{K Z}{\FFp}{X^{p-2}})}, which completes the proof of proposition~\ref{PROP2} in the case when  \m{\TEMPORARYI=1} since \m{{(-1)^{n}}\sqbracketss{}{K Z}{\FFp}{X^{p-2}}} generates \m{\ind_{KZ}^G \QUOTIENTI(1)}.

\subsection{\protect\boldmath Completing the proof for \texorpdfstring{\m{{\TEMPORARYI>1}}}{\textnu>1}}
\label{subsection5.5} In this subsection we assume that \m{\TEMPORARYI>1}. %
  Let \m{\delta = \TEMPORARYI}. Then \m{S_\delta = \{1,\ldots,\TEMPORARYI-1\}}, so
\algn{M = \MMM{\TEMPORARYI}{(p-1-\TEMPORARYI+i)_{1\LEQSLANT i \LEQSLANT \TEMPORARYI}}{\delta} = \left(\binom{p-1-\TEMPORARYI+i}{j}\right)_{1\LEQSLANT i\LEQSLANT \TEMPORARYI,\, 1\LEQSLANT j \LEQSLANT\TEMPORARYI-1}.}
By part~\eqref{comb10.5} of lemma~\ref{LEMMAX3.04} we have
 \algn{\left((-1)^{p-1-\TEMPORARYI+i}\binom{\TEMPORARYI}{i-1}\right)_{1\LEQSLANT i\LEQSLANT \TEMPORARYI} \, M = 0,}
 so  the quadruple
\algn{%
\Big(\TEMPORARYI,\delta,(p-1-\TEMPORARYI+i)_{1\LEQSLANT i \LEQSLANT \TEMPORARYI},\big(\binom{\TEMPORARYI}{i-1}\big)_{1\LEQSLANT i \LEQSLANT \TEMPORARYI}\Big)} satisfies equation~\eqref{eqg3oigj98j45g9j58ghufieofji3u5}, %
  and lemma~\ref{LEMMAX9} implies that  
   \Algn{\label{eq3j49j98jt945jt59jtu4}\nonumber
 & \textstyle
 \sum_{\TEMPORARYA \in \FF_p} {[\TEMPORARYA]^{p-\TEMPORARYI} \begin{Lmatrix}
    1 & 0 \\ [\TEMPORARYA] & 1
    \end{Lmatrix}}\sqbrackets{}{\INZB}{\QQp}{\GEN_{p-1}} +  \TEMPORARYI\sum_{\TEMPORARYA \in \FF_p} {[\TEMPORARYA]^{p-1} \begin{Lmatrix}
    1 & 0 \\ [\TEMPORARYA] & 1
    \end{Lmatrix}}\sqbrackets{}{\INZB}{\QQp}{\GEN_{\TEMPORARYI-1}}  +  \sum_{i=1}^{\TEMPORARYI-2} b_i {}\sqbrackets{}{\INZB}{\QQp}{\GEN_i} \\& \null\qquad\textstyle + \sum_{i=1}^\TEMPORARYI (-1)^{p-\TEMPORARYI+i}\binom{\TEMPORARYI}{i-1} {\binom{p-1-\TEMPORARYI+i}{\TEMPORARYI}\begin{Lmatrix}
0 & -1 \\ 1 & 0
\end{Lmatrix}}\sqbrackets{}{\INZB}{\QQp}{\GEN_{\TEMPORARYI}}\nonumber %
  \\ & \null\qquad\qquad \textstyle \EQUIVZEROIDEAL \sum_{\TEMPORARYA \in \FF_p} {A_\TEMPORARYA \begin{Lmatrix}
    1 & 0 \\ [\TEMPORARYA] & 1
    \end{Lmatrix}}\sqbrackets{}{\INZB}{\QQp}{\GEN} + \BIGO{\QUARK},}
    where, for \m{\TEMPORARYA\in\FF_p},
    \algn{A_\TEMPORARYA & \textstyle= \sum_{i=1}^\TEMPORARYI (-1)^{p-\TEMPORARYI+i} \binom{\TEMPORARYI}{i-1} \sum_{\MU \in \FF_p} \sum_{j\equiv 0\bmod{p-1}} \binom{p-1-\TEMPORARYI+i}{j}[\MU]^{j} [\TEMPORARYA]^{p-1-\TEMPORARYI} \begin{Lmatrix}
        0 & -1 \\ 1 & [\MU]
 \end{Lmatrix} \\ & \textstyle=  \sum_{\MU \in \FF_p} \sum_{i=1}^\TEMPORARYI (-1)^{p-\TEMPORARYI+i} \binom{\TEMPORARYI}{i-1} (1+\VERIFYIF{i=\TEMPORARYI}[\MU]^{p-1}) [\TEMPORARYA]^{p-1-\TEMPORARYI} \begin{Lmatrix}
        0 & -1 \\ 1 & [\MU]
 \end{Lmatrix}
,}
  and \m{b_1,\ldots,b_{\TEMPORARYI-2}\in\ZZp[G]}. 
  We can use equation~\eqref{eq029if984ujf9jf9wjf} with \m{\alpha\in\{1,\ldots,\TEMPORARYI-2\}} and part~\eqref{comb10.5} of lemma~\ref{LEMMAX3.04} to conclude that \algn{1\sqbrackets{}{\INZB}{\QQp}{\GEN_i} & \textstyle{}\in \BIGO{a\QUARK^{2-\alpha}}  \textstyle \textstyle \text{ for } i\in\{1,\ldots,\TEMPORARYI-2\}, \\ \textstyle \sum_{i=1}^\TEMPORARYI (-1)^{p-\TEMPORARYI+i}\binom{\TEMPORARYI}{i-1} \binom{p-1-\TEMPORARYI+i}{\TEMPORARYI} & \textstyle = 1 + \BIGO{p},} and therefore simplify equation~\eqref{eq3j49j98jt945jt59jtu4} as
 \Algn{\label{eq0394ikf3409rf4309fk49gj4j5}\nonumber
 & \textstyle
 \sum_{\TEMPORARYA \in \FF_p} {[\TEMPORARYA]^{p-\TEMPORARYI} \begin{Lmatrix}
    1 & 0 \\ [\TEMPORARYA] & 1
    \end{Lmatrix}}\sqbrackets{}{\INZB}{\QQp}{\GEN_{p-1}} +  \TEMPORARYI\sum_{\TEMPORARYA \in \FF_p} {[\TEMPORARYA]^{p-1} \begin{Lmatrix}
    1 & 0 \\ [\TEMPORARYA] & 1
    \end{Lmatrix}}\sqbrackets{}{\INZB}{\QQp}{\GEN_{\TEMPORARYI-1}}  %
    + \begin{Lmatrix}
0 & -1 \\ 1 & 0
\end{Lmatrix}\sqbrackets{}{\INZB}{\QQp}{\GEN_{\TEMPORARYI}}\nonumber %
  \\ & \null\qquad \textstyle \EQUIVZEROIDEAL -\sum_{\TEMPORARYA,\MU \in \FF_p} { (1+\TEMPORARYI[\MU]^{p-1}) [\TEMPORARYA]^{p-1-\TEMPORARYI} \begin{Lmatrix}
        0 & -1 \\ 1 & [\MU]
 \end{Lmatrix} \begin{Lmatrix}
    1 & 0 \\ [\TEMPORARYA] & 1
    \end{Lmatrix}}\sqbrackets{}{\INZB}{\QQp}{\GEN} + \BIGO{\QUARK}.}
    In particular, equation~\eqref{eq0f943if0943ik09f} with \m{\alpha=\TEMPORARYI} and equation~\eqref{eq029if984ujf9jf9wjf} with \m{\alpha=\TEMPORARYI-1} imply that
     \Algn{\label{eqf049fi4309ir0943ri09}
 & \textstyle
 \sum_{\TEMPORARYA \in \FF_p} {[\TEMPORARYA]^{p-\TEMPORARYI} \begin{Lmatrix}
    1 & 0 \\ [\TEMPORARYA] & 1
    \end{Lmatrix}}\sqbrackets{}{\INZB}{\QQp}{\GEN_{p-1}}   %
    + \begin{Lmatrix}
0 & -1 \\ 1 & 0
\end{Lmatrix}\sqbrackets{}{\INZB}{\QQp}{\GEN_{\TEMPORARYI}} %
   \EQUIVZEROIDEAL %
   \BIGO{\QUARK^{\TEMPORARYI}a^{-1}} \cup \BIGO{a\QUARK^{1-\TEMPORARYI}}.}
   If we multiply equation~\eqref{eq0394ikf3409rf4309fk49gj4j5} by \m{\begin{Smatrix}0 & -1 \\ 1 & 0\end{Smatrix}} we get
   \Algn{\label{eq03j4g9j4u8hrfeiwjf9i5jg954ujg984ugj94} \begin{Lmatrix}
0 & -1 \\ 1 & 0
\end{Lmatrix} \sum_{\TEMPORARYA \in \FF_p} {[\TEMPORARYA]^{p-\TEMPORARYI} \begin{Lmatrix}
    1 & 0 \\ [\TEMPORARYA] & 1
    \end{Lmatrix}}\sqbrackets{}{\INZB}{\QQp}{\GEN_{p-1}} - 1\sqbrackets{}{\INZB}{\QQp}{\GEN_{\TEMPORARYI}} \EQUIVZEROIDEAL \ERRR_3,}
    where
 \Algn{\label{eqf034jf9jf948jgh94jf94jf}\nonumber
  \textstyle
 \ERRR_3
\nonumber& \null  \textstyle \EQUIVZEROIDEAL \sum_{\TEMPORARYA,\MU \in \FF_p} { (1+\TEMPORARYI[\MU]^{p-1}) [\TEMPORARYA]^{p-1-\TEMPORARYI} \begin{Lmatrix}
        1 & [\MU] \\ 0 & 1
 \end{Lmatrix} \begin{Lmatrix}
    1 & 0 \\ [\TEMPORARYA] & 1
    \end{Lmatrix}}\sqbrackets{}{\INZB}{\QQp}{\GEN}\nonumber \\& \null\qquad\phantom{\EQUIVZEROIDEAL} \textstyle -\TEMPORARYI\begin{Lmatrix}
0 & -1 \\ 1 & 0
\end{Lmatrix}\sum_{\TEMPORARYA \in \FF_p} {[\TEMPORARYA]^{p-1} \begin{Lmatrix}
    1 & 0 \\ [\TEMPORARYA] & 1
    \end{Lmatrix}}\sqbrackets{}{\INZB}{\QQp}{\GEN_{\TEMPORARYI-1}} + \BIGO{\QUARK}\nonumber
    \\ & \null   \textstyle \EQUIVZEROIDEAL \frac{C_\TEMPORARYI}{a}\sum_{\MU \in \FF_p} { (1+\TEMPORARYI[\MU]^{p-1})  \begin{Lmatrix}
        p & [\MU] \\ 0 & 1
 \end{Lmatrix} }\sqbrackets{}{\INZB}{\QQp}{\GEN_\TEMPORARYI}\nonumber \\&\nonumber \null\qquad\phantom{\EQUIVZEROIDEAL} \textstyle -  \frac{\TEMPORARYI a} {C_{\TEMPORARYI-1}}\begin{Lmatrix}
0 & -1 \\ 1 & 0
\end{Lmatrix}\sum_{\TEMPORARYA ,\MU\in \FF_p} {[\TEMPORARYA]^{p-1}[\MU]^{p-\TEMPORARYI} \begin{Lmatrix}
    1 & 0 \\ [\TEMPORARYA] & p
    \end{Lmatrix}\begin{Lmatrix}
    1 & 0 \\ [\MU] & 1
    \end{Lmatrix}}\sqbrackets{}{\INZB}{\QQp}{\GEN} + \BIGO{\QUARK}
     \\ & \null   \textstyle \EQUIVZEROIDEAL \frac{C_\TEMPORARYI}{a}\sum_{\MU \in \FF_p} { (1+\TEMPORARYI[\MU]^{p-1})  \begin{Lmatrix}
         p & [\MU] \\ 0 & 1
 \end{Lmatrix} }\sqbrackets{}{\INZB}{\QQp}{\GEN_\TEMPORARYI}\nonumber \\& \nonumber\null \qquad\phantom{\EQUIVZEROIDEAL} \textstyle - \frac{\TEMPORARYI a}{C_{\TEMPORARYI-1}} \begin{Lmatrix}
0 & -1 \\ 1 & 0
\end{Lmatrix}\sum_{\TEMPORARYA \in \FF_p} {[\TEMPORARYA]^{p-1} \begin{Lmatrix}
    1 & 0 \\ [\TEMPORARYA] & p
    \end{Lmatrix} \begin{Lmatrix}
    0 & -1 \\ 1 & 0
    \end{Lmatrix}}\sqbrackets{}{\INZB}{\QQp}{\GEN_\TEMPORARYI} +  \BIGO{\QUARK} \cup \BIGO{a^2\QUARK^{2-2\TEMPORARYI}}
    \\& \null \EQUIVZEROIDEAL \textstyle \sum_{\MU \in \FF_p} { \left(\frac{C_\TEMPORARYI}{a}+ \frac{a^2 +  C_\TEMPORARYI C_{\TEMPORARYI-1}}{aC_{\TEMPORARYI-1}}\TEMPORARYI[\MU]^{p-1}\right)  \begin{Lmatrix}
       p & [\MU] \\ 0 & 1
 \end{Lmatrix} }\sqbrackets{}{\INZB}{\QQp}{\GEN_\TEMPORARYI} +  \BIGO{\QUARK} \cup \BIGO{a^2\QUARK^{2-2\TEMPORARYI}}.
   } 
The second congruence follows from lemma~\ref{LEMMAX6} applied to \m{\alpha\in\{\TEMPORARYI-1,\TEMPORARYI\}} and the set of lifts \m{\LCLASSTEICH}. %
      The third congruence follows from equation~\eqref{eq029if984ujf9jf9wjf} with \m{\alpha=0}, part~\eqref{propeta4} of lemma~\ref{LEMMAX3.01}, and equation~\eqref{eqf049fi4309ir0943ri09}. Therefore we have %
\newcommand{\tempheightaer}{\textstyle\vphantom{\EQUIVZEROIDEAL -{\frac{a^2}{\CONSTANT{\TEMPORARYI}\CONSTANT{\TEMPORARYI-1}}\binom{2\TEMPORARYI-1}{\TEMPORARYI} \sum_{\LAMBDA\in\FF_p} \begin{Lmatrix}
            p & [\LAMBDA] \\ 0 & 1
         \end{Lmatrix}}\sqbrackets{}{\INZB}{\QQp}{\GEN_{\TEMPORARYI}} + \BIGO{a\QUARK^{1-\TEMPORARYI}} \cup \BIGO
        {a^3\QUARK^{2-3\TEMPORARYI}}.}} \Algn{\label{eq04f93jfk349f}\nonumber
                & \textstyle \textstyle {\sum_{\MU\in\FF_p}\begin{Lmatrix}
                p & [\MU] \\ 0 & 1
         \end{Lmatrix}}\sqbrackets{}{\INZB}{\QQp}{\GEN_{\TEMPORARYI}}
            \\ & \nonumber \textstyle \null\qquad\tempheightaer  
            \EQUIVZEROIDEAL {\frac{a}{\CONSTANT{\TEMPORARYI}} \sum_{\TEMPORARYA,\MU\in\FF_p} [\TEMPORARYA]^{p-1-\TEMPORARYI} \begin{Lmatrix}
            1 & [\MU] \\ 0 & 1
         \end{Lmatrix}\begin{Lmatrix}
            1 & 0 \\ [\TEMPORARYA] & 1
            \end{Lmatrix} }\sqbrackets{}{\INZB}{\QQp}{\GEN} + \BIGO{a\QUARK^{1-\TEMPORARYI}}
    \\ & \nonumber \textstyle \null\qquad\tempheightaer 
    \EQUIVZEROIDEAL {\frac{a}{\CONSTANT{\TEMPORARYI}} \sum_{\MU\in\FF_p} \begin{Lmatrix}
            0 & -1 \\ 1 & 0
         \end{Lmatrix}\begin{Lmatrix}
            1 & 0 \\ [\MU] & 1
         \end{Lmatrix}}\sqbrackets{}{\INZB}{\QQp}{\GEN_{\TEMPORARYI-1}} + \ERRR_4 + \BIGO{a\QUARK^{1-\TEMPORARYI}}
        \\ & \nonumber \textstyle \null\qquad\tempheightaer 
        \EQUIVZEROIDEAL {\frac{a^2}{\CONSTANT{\TEMPORARYI}\CONSTANT{\TEMPORARYI-1}} \sum_{\TEMPORARYA,\MU\in\FF_p}   \begin{Lmatrix}
            0 & -1 \\ 1 & 0
            \end{Lmatrix} \begin{Lmatrix}
            1 & 0 \\ [\MU] & p
            \end{Lmatrix} [\TEMPORARYA]^{p-\TEMPORARYI}\begin{Lmatrix}
            1 & 0 \\ [\TEMPORARYA] & 1
         \end{Lmatrix}}\sqbrackets{}{\INZB}{\QQp}{\GEN} + \ERRR_4 + \BIGO{a\QUARK^{1-\TEMPORARYI}}%
           \\ & \nonumber \textstyle \null\qquad\tempheightaer 
        \EQUALZEROIDEAL {\frac{a^2}{\CONSTANT{\TEMPORARYI}\CONSTANT{\TEMPORARYI-1}} \sum_{\TEMPORARYA,\LAMBDA\in\FF_p} \begin{Lmatrix}
            p & [\LAMBDA] \\ 0 & 1
            \end{Lmatrix}  \begin{Lmatrix}
            0 & -1 \\ 1 & 0
            \end{Lmatrix} [\TEMPORARYA]^{p-\TEMPORARYI}\begin{Lmatrix}
            1 & 0 \\ [\TEMPORARYA] & 1
         \end{Lmatrix}}\sqbrackets{}{\INZB}{\QQp}{\GEN} + \ERRR_4 + \BIGO{a\QUARK^{1-\TEMPORARYI}}%
             \\ & \nonumber \textstyle \null\qquad\tempheightaer 
        \EQUIVZEROIDEAL -{\frac{a^2}{\CONSTANT{\TEMPORARYI}\CONSTANT{\TEMPORARYI-1}} \sum_{\TEMPORARYA,\LAMBDA\in\FF_p} \begin{Lmatrix}
            p & [\LAMBDA] \\ 0 & 1
            \end{Lmatrix}  \begin{Lmatrix}
            0 & -1 \\ 1 & 0
            \end{Lmatrix} [\TEMPORARYA]^{p-\TEMPORARYI}\begin{Lmatrix}
            1 & 0 \\ [\TEMPORARYA] & 1
         \end{Lmatrix}}\sqbrackets{}{\INZB}{\QQp}{\GEN_{p-1}} + \ERRR_4 + \BIGO{a\QUARK^{1-\TEMPORARYI}}%
         \\ & \textstyle \null\qquad\tempheightaer 
        \EQUIVZEROIDEAL -{\frac{a^2}{\CONSTANT{\TEMPORARYI}\CONSTANT{\TEMPORARYI-1}} %
        \sum_{\LAMBDA\in\FF_p} \begin{Lmatrix}
            p & [\LAMBDA] \\ 0 & 1
         \end{Lmatrix}}\sqbrackets{}{\INZB}{\QQp}{\GEN_{\TEMPORARYI}} + \ERRR_4 + \ERRR_5 + \BIGO{a\QUARK^{1-\TEMPORARYI}}%
         ,
        }
        where \Algn{\label{eqdk93fj923jf9283jd}
\textstyle \ERRR_4 %
 \EQUIVZEROIDEAL \frac{a}{\CONSTANT{\TEMPORARYI}}\ERRR^\prime,
   \text{ and } \textstyle \ERRR_5  \null  \textstyle \EQUIVZEROIDEAL -  \frac{a^2}{C_\TEMPORARYI C_{\TEMPORARYI-1}}\sum_{\LAMBDA \in \FF_p}  \begin{Lmatrix}
       p & [\LAMBDA] \\ 0 & 1
 \end{Lmatrix}\ERRR_3.}
The first and third congruences %
 follow from lemma~\ref{LEMMAX6} applied to \m{\alpha\in\{\TEMPORARYI-1,\TEMPORARYI\}} and the set of lifts \m{\LCLASSTEICH}. %
  The second congruence follows from equation~\eqref{eqfu3h4g8fh34hf3487fh387}. 
   For the equality we  change  the variable \m{\MU} to \m{\LAMBDA=-\MU}. The fourth congruence  follows from equation~\eqref{eq029if984ujf9jf9wjf} with \m{\alpha=0} and part~\eqref{propeta4} of lemma~\ref{LEMMAX3.01}. The fifth congruence 
 follows from equation~\eqref{eq03j4g9j4u8hrfeiwjf9i5jg954ujg984ugj94}.
 Then equation~\eqref{eq04f93jfk349f} implies that
\Algn{\label{eqdj2398i9834j98f3}\textstyle {\left(1 +\frac{a^2}{\CONSTANT{\TEMPORARYI}\CONSTANT{\TEMPORARYI-1}} \right)\sum_{\LAMBDA\in\FF_p}\begin{Lmatrix}
                p & [\LAMBDA] \\ 0 & 1
              \end{Lmatrix}}\sqbrackets{}{\INZB}{\QQp}{\GEN_{\TEMPORARYI}} \EQUIVZEROIDEAL \BIGO{\QUARK^\TEMPORARYI a^{-1}}\cup\BIGO{a\QUARK^{1-\TEMPORARYI}}.}
By equation~\eqref{eqReduct},
            \algn{1 +\frac{a^2}{\CONSTANT{\TEMPORARYI}\CONSTANT{\TEMPORARYI-1}} \in  \BIGO{\QUARK^\TEMPORARYI a^{-1}}\cup\BIGO{\QUARK^{1-\TEMPORARYI} a} %
            \Longleftrightarrow a\in \CIRCLES_{\TEMPORARYI} \Longleftrightarrow \TEMPORARYMU\in\ZZp.}
           If   \m{a\not\in\CIRCLES_\TEMPORARYI} then after dividing equation~\eqref{eqdj2398i9834j98f3}  by \m{1 +\frac{a^2}{\CONSTANT{\TEMPORARYI}\CONSTANT{\TEMPORARYI-1}}} and reducing modulo \m{\maxideal}, we get by equation~\eqref{eqclass2} that \m{T_\sigma({(-1)^n}\sqbracketss{}{K Z}{\FFp}{X^{p-2}})} is represented by an element of \m{\INAAA  } and therefore by an element of  the kernel of \m{\QUOTIENTMAP}. If \m{a\in\CIRCLES_\TEMPORARYI} then \m{\TADICVALUATION(a)=\TEMPORARYI-\frac12} and \m{\TEMPORARYMU\in\ZZp}, and equation~\eqref{eq04f93jfk349f} implies that
     \Algn{\label{eq958tugtue09uhy67u4e09gut9jrwe}
     \textstyle \TEMPORARYMU \sum_{\LAMBDA\in\FF_p} {\begin{Lmatrix}
            p & [\LAMBDA] \\ 0 & 1
            \end{Lmatrix} }\sqbrackets{}{\INZB}{\QQp}{\GEN_\TEMPORARYI} %
            \EQUIVZEROIDEAL -\frac{C_{\TEMPORARYI-1}}{a} \left( \ERRR_4+\ERRR_5\right)+ \BIGO{\QUARK^\TEMPORARYI a^{-1}} \cup \BIGO{a \QUARK^{1-\TEMPORARYI}}.
            }
           By equations~\eqref{eqclass1}, \eqref{eqclass2}, \eqref{eq09g54ig9045i5irf094i5f0}, \eqref{eqfu3h4g8fh34hf3487fh387}, \eqref{eqf034jf9jf948jgh94jf94jf}, and~\eqref{eqdk93fj923jf9283jd}, the reductions of 
\algn{\TEMPORARYMU\sum_{\LAMBDA\in\FF_p}\begin{Lmatrix}
                p & [\LAMBDA] \\ 0 & 1
              \end{Lmatrix}\sqbrackets{}{\INZB}{\QQp}{\GEN_{\TEMPORARYI}} ,\,\,
               -\frac{C_{\TEMPORARYI-1}}{a}  \ERRR_4 ,\,\,-\frac{C_{\TEMPORARYI-1}}{a}  \ERRR_5}
modulo \m{\maxideal}  are
\algn{\phi\, {\textstyle T_\sigma({(-1)^{n}}\sqbracketss{}{K Z}{\FFp}{X^{p-2}})},\,\, {(-1)^{n}}\sqbracketss{}{K Z}{\FFp}{X^{p-2}},\,\, {\textstyle T_\sigma^2({(-1)^{n}}\sqbracketss{}{K Z}{\FFp}{X^{p-2}})},} respectively.
 So equation~\eqref{eq958tugtue09uhy67u4e09gut9jrwe} implies that the kernel of \m{\QUOTIENTMAP} contains a representative of \m{(T_\sigma^2-\phi\, T_\sigma+1)({(-1)^{n}}\sqbracketss{}{K Z}{\FFp}{X^{p-2}})}, which completes the proof of proposition~\ref{PROP2} in the case when  \m{\TEMPORARYI>1} since \m{{(-1)^{n}}\sqbracketss{}{K Z}{\FFp}{X^{p-2}}} generates \m{\ind_{KZ}^G \QUOTIENTI(\TEMPORARYI)}. \BLACKSQUARE

\end{document}